\documentclass[11pt]{amsart}

\usepackage{geometry}
\usepackage{lmodern, microtype}
\usepackage{longtable}
\usepackage{listings}
\usepackage{amsthm, amsmath, amssymb, amsfonts}
\usepackage{mathtools}
\mathtoolsset{showonlyrefs=true}

\usepackage{mathrsfs} 
\allowdisplaybreaks[3] 
\usepackage{hyperref}

\usepackage{enumitem}

\hypersetup{
    colorlinks=true,
    linkcolor=blue,
    filecolor=magenta,      
    urlcolor=cyan,
    citecolor=red,
}

\theoremstyle{plain}

\theoremstyle{definition}

\theoremstyle{plain}

\theoremstyle{definition}

\theoremstyle{definition}
\newtheorem{dn}{Definition}[section]
\newtheorem{kh}[dn]{Notation}
\newtheorem{nxe}[dn]{Remark}
\newtheorem{exm}[dn]{Example}

\theoremstyle{plain}
\newtheorem{bd}[dn]{Lemma}
\newtheorem{theo}[dn]{Theorem}
\newtheorem{mde}[dn]{Proposition}
\newtheorem{hqu}[dn]{Corollary}

\numberwithin{equation}{section}

\DeclareMathOperator{\Ker}{Ker}

\def\DD{D\kern-.7em\raise0.4ex\hbox{\char '55}\kern.33em}

\title[Hit problem of rank 5 and Milnor-operation layers]{The Peterson hit problem of rank 5\\ and canonical Milnor-operation layers}

\author{\DD\d{\u a}ng V\~o Ph\'uc}
\address{Department of Mathematics, FPT University, Quy Nhon AI Campus, An Phu Thinh New Urban Area, Vietnam}
\email{dangphuc150488@gmail.com}
\thanks{ORCID: \url{https://orcid.org/0000-0002-6885-3996}}


\begin{document}

\begin{abstract}
Let $P_n$ be the polynomial algebra over $\mathbb F_2$ viewed as a module over the mod $2$ Steenrod algebra $\mathbb A$. We develop a canonical structural layer for the Peterson hit problem by proving that, for the Milnor basis element $P_s^0$, the kernel of its action $T_{s,e}: (P_n)_e \to (P_n)_{e+2^s-1}$ coincides exactly with the Frobenius-square subspace whenever $0 \leq e \leq 2^s$. Under monomial divided-power duality, the transpose adjoint dictates that every functional on the non-hit quotient satisfies a strict homological $P_s^0$-orthogonality test. We also formulate a dual primitive-kernel algorithm utilizing the homological right actions of $Sq^{2^r}$, providing a versatile framework independent of classical admissible monomial reductions, alongside an adjoint proof of Kameko's epimorphism. We apply this theoretical machinery to investigate the fifth Singer algebraic transfer in the highly non-trivial degree family $d_t = 2^{t+5} + 2^{t+2} + 2^{t+1} - 5$. Specializing the Milnor layer to rank $5$ yields an exact theoretical rank formula that structurally verifies our algorithmic computations. Using \texttt{OSCAR}-certified row reductions, we explicitly determine the cohit spaces, proving $\dim(QP_5)_{d_0} = 1322$ and $\dim(QP_5)_{d_t} = 2841$ for $t \geq 1$, and construct their $GL_5$-invariant generators. Ultimately, we establish that $\dim((QP_5)_{d_t})^{GL_5} = 1$, confirming that the fifth cohomological transfer is an isomorphism in degree $d_t$ for all $t \geq 0$. 
\end{abstract}

\maketitle

\section{Introduction}\label{s11}

The mod $2$ Steenrod algebra $\mathbb A$ is the algebra of stable mod $2$ cohomology operations, and its cohomology $\operatorname{Ext}_{\mathbb A}^{*,*}(\mathbb Z/2,\mathbb Z/2)$
is the $E_2$-term of the Adams spectral sequence converging to the stable homotopy groups of spheres \cite{J.A}. In \cite{Singer}, Singer introduced the cohomological transfer as a natural bridge from the homology of elementary abelian $2$-groups to this Ext algebra. For an $n$-dimensional vector space $V^{\oplus n}$ over $\mathbb Z/2$, the rank-$n$ transfer has the form
\[
 \phi_n^{*}(\mathbb Z/2):
 (\mathbb Z/2\otimes_{GL_n}P_{\mathbb A}H_*(V^{\oplus n}))_d
 \longrightarrow
 \operatorname{Ext}_{\mathbb A}^{n,n+d}(\mathbb Z/2,\mathbb Z/2),
\]
where $P_{\mathbb A}H_*(V^{\oplus n})$ denotes the subspace annihilated by all positive Steenrod operations. The domain is dual to the invariant space $\big((QP_n)_d\big)^{GL_n}$, where
\[
P_n=\mathbb Z/2[u_1,\ldots,u_n],\qquad
QP_n=\mathbb Z/2\otimes_{\mathbb A}P_n=P_n/\overline{\mathbb A}P_n,
\]
and $(QP_n)_d$ denotes the homogeneous component of $QP_n$ in degree $d$. Thus the Singer transfer is inseparable from the \textit{Peterson hit problem} \cite{F.P}: determining the quotient space $QP_n$, or equivalently a minimal set of homogeneous $\mathbb A$-generators for $P_n$.

The hit problem became a central computational problem in the theory of unstable modules over the Steenrod algebra. Peterson's vanishing conjecture, which was proved by Wood \cite{R.W}, explicitly states that $(QP_n)_d = 0$ whenever $\alpha(d+n) > n$, where $\alpha(k)$ denotes the number of ones in the binary expansion of an integer $k$. To systematically compute the non-vanishing components, Kameko \cite{M.K} introduced a degree-reducing squaring operation $Sq^0_*: P_n \longrightarrow P_n$. On a monomial basis, it is defined by the formula
\[
Sq^0_*(u_1^{m_1}\cdots u_n^{m_n}) =
\begin{cases}
u_1^{k_1}\cdots u_n^{k_n} & \text{if } m_i = 2k_i+1 \text{ for all } 1 \leq i \leq n, \\
0 & \text{otherwise}.
\end{cases}
\]
This operation respects the Steenrod algebra action and thus induces a well-defined homomorphism on the quotient spaces, denoted by
\[
(\widetilde{Sq^0_*})_d: (QP_n)_d \longrightarrow (QP_n)_{k},
\]
where $d = 2k+n$. Crucially, Kameko \cite{M.K} proved that this induced homomorphism is always surjective, making it one of the most effective reduction tools for explicit computations. 
A systematic account of the method and many of its refinements is given by Walker and Wood \cite{W.W}. The cases $n\leq 4$ are well understood (see \cite{W.W, W.W2} and references therein), whereas rank $5$ remains substantially more delicate because the quotient spaces and their $GL_5$-invariants grow rapidly and because the transfer calculation requires information not only about dimensions but also about explicit invariant representatives. Singer proved his injectivity conjecture for $n\leq2$, Boardman proved it for $n=3$, and the rank-$4$ case was verified in our previous papers \cite{D.P6-1, P8, P9}; the conjecture is not valid in full generality, as shown by the rank-$6$ counterexample in our recent preprint \cite{Phuc2026}.

Investigating the unresolved rank-$5$ case requires focusing on specific degree sequences that encode meaningful topological information. The emergence of these sequences is intrinsically tied to the underlying algebraic structure. Because the Kameko homomorphism systematically relates the quotient space in degree $d$ to that in degree $k$, iterating the relation $d = 2k + 5$ naturally generates infinite families of degrees. Furthermore, a meaningful evaluation of the fifth Singer transfer dictates a strict focus on degree families that support non-trivial $GL_5$-invariants, as these are the precise components that detect surviving classes in the Ext algebra. Motivated by this framework, the present paper studies the rank-$5$ hit problem and the fifth cohomological transfer in the highly non-trivial degree family
\[
d_t=2^{t+5}+2^{t+2}+2^{t+1}-5\qquad (t\geq0).
\]
Complementing this structural necessity is a fundamental arithmetic reduction. It is a well-known fact that the rank-$n$ hit problem is completely determined by computations in degrees $d$ satisfying $\mu(d)\leq n$, where the arithmetic function $\mu: \mathbb N\longrightarrow \mathbb N$ is detailed in Section~\ref{s2}. The degree family $d_t$ intrinsically satisfies this bounding condition for $n=5$ (see Section \ref{s3}).
In particular, it is intrinsically linked to the Kameko reduction process: the homomorphism $(\widetilde{Sq^0_*})_{d_t}$ reduces this infinite sequence step-by-step down to the two base degrees $d_0 = 33$ and $d_1 = 71$. We refer the reader to Theorem \ref{dlWK}(ii) in Section \ref{s2} and carry out explicit computations in Section \ref{s3} for a clearer understanding of the practicality of that reduction. We compute the relevant quotient spaces and invariants in these degrees explicitly. In degree $d_0=33$, we prove
\[
\dim (QP_5)_{d_0}=1322,
\qquad
\dim\big((QP_5)_{d_0}\big)^{GL_5}=1.
\]
In degree $d_1=71$, we prove
\[
\dim \Ker((\widetilde{Sq^0_*})_{d_1})=1519,
\qquad
\big(\Ker((\widetilde{Sq^0_*})_{d_1})\big)^{GL_5}=0,
\]
and consequently, by propagating these structural results through the Kameko sequence, we obtain
\[
\dim (QP_5)_{d_t}=2841\quad (t\geq1),
\qquad
\dim\big((QP_5)_{d_t}\big)^{GL_5}=1\quad (t\geq0).
\]
Together with the known Ext classes detected by the fifth transfer, these computations show that $\phi_5^*(\mathbb Z/2)$ is an isomorphism in internal degree $d_t$ for every $t\geq0$.

A second purpose of this paper is to formally bridge the Peterson hit problem with algebraic coding theory to extract a profound structural layer from the calculation. This layer is motivated by a different, but highly compatible, use of Steenrod sub-Hopf algebras. Dougherty and Vergili \cite{DV2017} introduced codes over the finite sub-Hopf algebras of the Steenrod algebra and developed the corresponding Euclidean, Hermitian, and additive dualities. Vergili and Karaca \cite{VK2017} then used the Milnor basis to construct a canonical length-one Euclidean and Hermitian self-dual code generated by the Milnor element $P_{n+1}^{0}$. In our recent work \cite{PhucCodeAnn}, this code-theoretic object was studied through its external action on polynomial modules: for the homogeneous action map $T=P_{n+1}^{0}$, the orthogonal complement of the source-degree kernel is the image of the transpose adjoint, and every functional on the relevant target-degree non-hit quotient must vanish on the image of $T$. Thus, building upon these coding-theoretic frameworks, the current work establishes a novel finite-dimensional orthogonality mechanism attached to the Milnor element, demonstrating how algebraic codes can govern the structural layers of the hit problem.

For the degree family considered here, the mechanism appears in a particularly natural form. We have
\[
d_t=(2^{t+5}-1)+(2^{t+2}+2^{t+1}-4).
\]
Since $2^{t+5}-1$ is the degree of the Milnor basis element $P_{t+5}^{0}$, which is primitive because it is dual to $\xi_{t+5}$, this operation gives a natural map
\[
T_t=P_{t+5}^{0}:(P_5)_{e_t}\longrightarrow (P_5)_{d_t},
\qquad e_t=2^{t+2}+2^{t+1}-4.
\]
Equivalently, this is the specialization of the companion operation $P_{n+1}^{0}$ with $n=t+4$. Strikingly, we prove that the kernel of this homological action is precisely the Frobenius-square subspace, yielding the exact formulas:
\[
\Ker(T_t)=\{h^2:h\in(P_5)_{3\cdot2^t-2}\},
\qquad
\operatorname{rank}(T_t)=\binom{6\cdot2^t}{4}-\binom{3\cdot2^t+2}{4}.
\]
More generally, Section~\ref{sec:scope-milnor-test} proves a sharp low-source-degree theorem for the Milnor maps
\[
T_{s,e}=P_s^0:(P_n)_e\longrightarrow (P_n)_{e+2^s-1}.
\]
For every rank $n$ and every $0\leq e\leq2^s$, the kernel of $T_{s,e}$ is exactly the Frobenius-square subspace in degree $e$; when $n\geq2$, this range is optimal, since a non-square kernel element already appears in degree $e=2^s+1$.

Determining the exact dimension of $\Ker(T_t)$ and the rank of $T_t$ serves three fundamental purposes for the current hit problem:

\begin{enumerate}[label=(\arabic*)]
    \item \textbf{Establishing a theoretical lower bound:} Although admissible-basis algorithms ultimately reduce large relation matrices to quotient dimensions, they do not by themselves exhibit a distinguished part of the hit subspace. Since $T_t$ has positive degree, its entire image lies within the hit space $H_{d_t} = (\overline{\mathbb A}P_5)\cap (P_5)_{d_t}$. The exact rank formula therefore provides a theoretical lower bound, $\dim H_{d_t} \geq \operatorname{rank}(T_t)$, which is proved independently of the software computation.
    
    \item \textbf{Providing a canonical structure for the hit subspace:} While standard Steenrod relations and Kameko reductions inevitably produce massive spanning families with highly complex overlaps, our theorem extracts an underlying geometric regularity. By proving that $\Ker(T_t)$ is precisely the Frobenius-square subspace, we successfully isolate an algebraically transparent, non-overlapping pattern of hit elements. This defines a canonical $P_{t+5}^{0}$-hit layer inside the hit subspace. It is canonical in the invariant-theoretic sense because Steenrod operations commute with the standard $GL_5$-action on $P_5$.
    
    \item \textbf{Creating a strict orthogonality test for the cohit space:} To evaluate the Singer transfer, one must find explicit representatives for dual classes in the cohit space $H_{d_t}^{\perp}$. Since $P_{t+5}^{0}$ is a Steenrod operation of strictly positive degree, its image inherently consists of hit elements, yielding $\operatorname{Im}(T_t) \subseteq H_{d_t}$. Consequently, any surviving cohit representative must be orthogonal to $\operatorname{Im}(T_t)$. Using the monomial divided-power pairing and the transpose-adjoint formalism from our previous work \cite{PhucCodeAnn}, we obtain the orthogonal condition
    \[
    H_{d_t}^{\perp}\subseteq \Ker(T_t^*),
    \]
    where $H_{d_t}^{\perp}$ denotes the orthogonal complement of the hit subspace with respect to this pairing. Equivalently, every primitive homology representative in degree $d_t$ is annihilated by the homological right action of $P_{t+5}^{0}$. This acts as a rigorous filter: any invariant generator found by the software in Section~\ref{s3} must strictly satisfy this $P_{t+5}^{0}$-orthogonality condition to be a valid dual representative.
\end{enumerate}

This framework complements the machine-certified computations by providing an independent structural verification attached to the associated Milnor operation. In Section~\ref{sec:scope-milnor-test}, we also formulate a second exact algorithmic route, independent of admissible-basis enumeration: the dual primitive-kernel algorithm computes the transfer domain directly as an intersection of kernels of homological right Steenrod operations and then takes $GL_n$-coinvariants. The same divided-power formalism gives a self-contained adjoint proof of Kameko's epimorphism and explains why Kameko-lifted primitive representatives automatically pass every $P_s^0$-right-action test. This algorithm is not presented as a universal closed-form solution to the hit problem, since its matrices still grow rapidly with $n$ and $d$, but it supplies a rigorous dual-side criterion for verification, and can be applied to other degrees and ranks whenever the corresponding sparse linear algebra is computationally feasible.

The invariant computations in the paper are presented as finite-dimensional linear algebra over $\mathbb Z/2$. As will be detailed in Section~\ref{s2}, the space of degree-$d_t$ polynomials can be graded by the weight vectors $\omega$ of their underlying monomials, decomposing the quotient space into associated weight components $(QP_5)_{d_t}(\omega)$. For each such weight component, a candidate degree-$d_t$ polynomial $f$ is expressed as a linear combination of an ordered admissible basis. The conditions for $f$ to be $GL_5$-invariant are then formulated as the linear system
\[
\rho_i(f)-f\sim_\omega 0 \qquad (1 \leq i \leq 5),
\]
where $\sim_\omega$ denotes equivalence modulo the hit subspace in that specific weight and degree, $\rho_1,\ldots,\rho_4$ generate the symmetric group $\Sigma_5\subset GL_5$, and $\rho_5$ is a transvection completing the generating set for $GL_5$. The resulting matrices are row reduced over $\mathbb Z/2$ in \texttt{OSCAR}; the ranks and kernel dimensions used in the proof are recorded in the body of the paper, and the associated data files are supplied in Appendix~\ref{s4}. This makes the admissible-basis enumeration and the $GL_5$-invariant calculation reproducible rather than dependent on informal parallel computations.

The paper is organized as follows. Section~\ref{s2} recalls the necessary material on weight vectors, inadmissible monomials, Kameko's homomorphism, and the orthogonality criterion via the transpose adjoint. Section~\ref{s3} proves the main results: the $P_{t+5}^{0}$-hit layer theorem, the computations in degrees $d_0=33$ and $d_1=71$, the determination of the $GL_5$-invariants, and the verification of the fifth Singer transfer in the family $d_t$. Section~\ref{sec:scope-milnor-test} explains the general scope of the Milnor-operation test, proves the divided-power adjoint form of Kameko's epimorphism, and records the dual primitive-kernel algorithm. Appendix~\ref{s4} records the computational implementation and data repositories.

\section{Background}\label{s2}

In this section, we recall basic material on the hit problem in a concise form; for more context, consult two monographs \cite{W.W, W.W2}.

\begin{dn}[{\bf Weight vector and exponent vector}] We say that a sequence of non-negative integers $\omega = (\omega_1, \omega_2, \ldots, \omega_i,\ldots)$ is a \textit{weight vector}, if $\omega_i  = 0,$ for $i\gg 0.$ Then, one defines $\deg(\omega) = \sum_{i\geq 1}2^{i-1}\omega_i.$ With a natural number $d,$ let denote $\alpha_j(d)$ the $j$-th coefficients in dyadic expansion of $d,$ then $\alpha(d) = \sum_{j\geq 0}\alpha_j(d),$ and $d = \sum_{j\geq 0}\alpha_j(d)2^j,$ where $\alpha_j(d)\in \{0, 1\}$ for all $j\geq 0.$ For a monomial $u = u_1^{a_1}u_2^{a_2}\ldots u_n^{a_n}\in P_n,$ let us consider two sequences $\omega(u) :=(\omega_1(u), \ldots, \omega_i(u), \ldots)$ and $(a_1, a_2, \ldots, a_n),$ where $\omega_i(u)=\sum_{1\leq j\leq n}\alpha_{i-1}(a_j),$ for every $i.$ They are respectively called the {\it weight vector} and the \textit{exponent vector} of $u.$  By convention, the sets of all the weight vectors and the exponent vectors are given the left lexicographical order.
\end{dn}

\begin{dn}[{\bf Linear order on {\boldmath{$P_n$}}}]
Let $u = u_1^{a_1}u_2^{a_2}\ldots u_n^{a_n}$ and \mbox{$y = u_1^{b_1}u_2^{b_2}\ldots u_n^{b_n}$} be monomials in $P_n.$ We write $a,\, b$ for the exponent vectors of $u$ and $y,$ respectively. We say that $u  < y$ if and only if one of the following holds:
\begin{enumerate}
\item[(i)] $\omega(u) < \omega(y);$
\item[(ii)] $\omega(u) = \omega(y)$ and $a < b.$
\end{enumerate}
\end{dn}

\begin{dn}[{\bf Equivalence relations on {\boldmath{$P_n$}}}]  For a weight vector $\omega$ of degree $d,$ let us denote two subspaces associated with $\omega$ by
$$\begin{array}{ll}
\medskip
(P_n)_d(\omega) &= \langle\{ u\in (P_n)_d|\, \omega(u)\leq \omega\}\rangle,\\
(P_n)_d(< \omega) &= \langle \{ u\in (P_n)_d|\, \omega(u) < \omega\}\rangle.
\end{array}$$
Let $f$ and $g$ be two homogeneous polynomials in $(P_n)_d.$ Equivalence relations "$\sim$" and "$\sim_{\omega}$" can be defined on $P_n$ in the following manner:
\begin{enumerate}
\item [(i)] $f \sim g $ if and only if $(f - g)\in \overline{\mathbb A}P_n$. (Note that since working mod 2, $f-g = f+g.$)
\item[(ii)]  $f \sim_{\omega} g$ if and only if $f, \, g\in (P_n)_d(\omega)$ and $(f -g)\in (\overline{\mathbb A}P_n\cap (P_n)_d(\omega) + (P_n)_d(< \omega)).$

In particular,  if $f\sim 0$ (resp. $f\sim_{\omega} 0$), then we say that $f$ is \textit{hit} (resp. \textit{$\omega$-hit}). 

Let us denote $(QP_n)_d(\omega)$ as the quotient space of $(P_n)_d(\omega)$ via the equivalence relation "$\sim_{\omega}$". According to \cite{N.S3, W.W}, this $(QP_n)_d(\omega)$ is also a $GL_n$-module and
$$(QP_n)_d\cong \bigoplus_{\deg(\omega) = d}(QP_n)_d(\omega).
$$
If we put
$$ \begin{array}{ll}
\medskip
 (P_n)^{0} &:= \langle\{u = \prod_{1\leq j\leq n}u_j^{\gamma_j} \in P_n:\ \prod_{1\leq j\leq n}\gamma_j = 0\} \rangle,\\
  (P_n)^{> 0} &:= \langle\{u = \prod_{1\leq j\leq n}u_j^{\gamma_j}\in P_n:\  \prod_{1\leq j\leq n}\gamma_j > 0\} \rangle,
\end{array}$$
then $(P_n)^{0}$ and $(P_n)^{>0}$ are $\mathbb A$-submodules of $P_n.$ We denote by $(QP_n)^{0}:= Q(P_n)^{0}$  and $(QP_n)^{>0}:= Q (P_n)^{>0}.$ Then $ QP_n\cong  (QP_n)^{0}\bigoplus  (QP_n)^{>0}.$ Hence,  $(QP_n)_d\cong  (QP_n)_d^{0}\bigoplus  (QP_n)_d^{>0}$ for any positive degree $d.$
\end{enumerate}
\end{dn}

\begin{dn}[{\bf Admissible monomial and inadmissible monomial}] A monomial $u$ in $P_n$ is said to be {\it inadmissible} if there exist monomials $u_1, u_2,\ldots, u_k$ in $P_n$ such that $u_j < u$ for all $j,\ 1\leq j\leq k,$ and $u \sim  (u_1 + u_2 + \cdots + u_k).$ Then, we say that $u$ is {\it admissible}, if it is not inadmissible.

\medskip

Additionally, one defined strictly inadmissible monomials as follows. A monomial $u$ in $P_n$ is said to be {\it strictly inadmissible} if there exist monomials $u_1, u_2,\ldots, u_k$ in $P_n$ such that $u_j < u$ for all $j,\ 1\leq j\leq k,$ and $u = \sum_{1\leq j\leq k}u_j + \sum_{0\leq  m\leq s - 1}Sq^{2^{m}}(u'_m),$ where $s = {\rm max}\{i\in\mathbb N: \omega_i(u) > 0\}$ and suitable polynomials $u'_m\in P_n$ with $\deg(u'_m)  = \deg(u) - 2^{m}.$
\end{dn}

The following important theorem, which is due to Kameko \cite{M.K} and Sum \cite{N.S1}, can be viewed as a criterion in the calculation of the inadmissible monomials.

\begin{theo}[{\bf Criteria on inadmissible monomials}]\label{dlKS}
The following assertions are true:
\begin{itemize}
\item[(i)] {\rm (see \cite{M.K})} Let $z, w$ be monomials in $P_n$ and let $r$ be a positive integer. Suppose that $\omega_r(z)\neq 0$ and $\omega_j(z) = 0$ for all $j > r.$ Consequently, if $z$ is strictly inadmissible, then so is $zw^{2^{r}}.$
\item[(ii)] {\rm (see \cite{N.S1})} Let $y$ and $z$ be monomials in $P_n.$ For an integer $r >0,$ assume that $\omega_i(y) = 0$ for all $i > r$. Consequently, if $z$ is inadmissible, then so is $yz^{2^r}.$
\end{itemize}
\end{theo} 

From the above data, we have seen that $(QP_n)_d$ is a $\mathbb Z/2$-vector space with a basis consisting of all the classes represented by the admissible monomials in $(P_n)_d.$

\begin{dn}[{\bf Spike monomial}]
A monomial $u_1^{b_1}u_2^{b_2}\ldots u_n^{b_n}$ in $P_n$ is called a {\it spike} if every exponent $b_j$ is of the form $2^{c_j} - 1.$ In particular, if the exponents $c_j$ can be arranged to satisfy $c_1 > c_2 > \ldots > c_{r-1}\geq c_r \geq 1,$ where only the last two smallest exponents can be equal, and $c_j = 0$ for $ r+1 \leq j  \leq n,$ then $u_1^{b_1}u_2^{b_2}\ldots u_n^{b_n}$ is called a {\it minimal spike}.
\end{dn}

\begin{bd}[\cite{P.S1}]\label{dlPS}
All the spikes in $P_n$ are admissible and their weight vectors are weakly decreasing. Furthermore, if a weight vector $\omega = (\omega_1, \omega_2, \ldots)$ of finite degree is weakly decreasing and $\omega_1\leq n,$ then there is a spike $z\in P_n$ such that $\omega(z) = \omega.$
\end{bd}

\begin{kh}

It will be helpful to use the following notation: we will denote the set of all admissible monomials of degree $d$ in the $\mathbb A$-module $P_n$ by $\mathscr C^{\otimes n}_d$. For a monomial $u\in (P_n)_d,$ write $[u]$ for the equivalence class of $u$ in $(QP_n)_d.$ If $\omega$ is a weight vector and $u\in (P_n)_d(\omega),$ then we denote by $[u]_\omega$ the equivalence class of $u$ in $(QP_n)_d(\omega).$ Moreover if $\omega$ is a weight vector of a minimal spike, then $[u]_{\omega} = [u].$ 
\end{kh}

The following tools are useful for studying the hit problem. We consider the arithmetic function $\mu: \mathbb N\longrightarrow \mathbb N$ that is determined by
$$ \begin{array}{ll}
\mu(d) &= \mbox{min}\big\{n\in \mathbb N:\ \alpha(d + n)\leq n\big\}\\
& = {\rm min}\{n\in\mathbb N:\,\,d = \sum_{1\leq j\leq n}(2^{s_j}-1),\,\,s_j> 0,\,\,1\leq j\leq n\}.
\end{array}$$

\begin{theo}\label{dlWK} 
For each positive integer $d,$ the following hold:
\begin{itemize}
\item[(i)] {\rm (see \cite{R.W})}. $(QP_n)_d$ is trivial if $\mu(d) > n.$ Consequently $QP_n$ is trivial in degrees $d$ unless $d$ is of the form $d = 2^{t_1} + 2^{t_2}+ \cdots + 2^{t_n}-n$ with $t_1\geq t_2\geq \cdots\geq t_n \geq 0.$ 

\item[(ii)] {\rm (see \cite{M.K})}. The Kameko homomorphism $(\widetilde {Sq^0_*})_{2d+n}$ is an isomorphism if and only if $\mu(2d+n) = n.$ Here the map $(\widetilde {Sq^0_*})_{2d+n}$ is defined as in Section \ref{s11}.
\end{itemize}
\end{theo}

Thus, utilizing the aforementioned theorem and the isomorphism between the domain of the Singer transfer and $((QP_n)_d)^{GL_n}$, we can focus exclusively on studying the hit problem for the $\mathbb A$-module $P_n$ and the behavior of the Singer transfer $\phi_n^{*}(\mathbb Z/2)$ in degrees $d$ such that $\mu(d) \leq n$.

In \cite{Singer2}, Singer observed that if $\mu(d)\leq n,$ then there exists a unique minimal spike of degree $d$ in $P_n$. Additionally, he proved the following outcome.

\begin{theo}[cf. Singer \cite{Singer2}]\label{dlSinger}
Suppose that $u$ is a monomial in $P_n$ with $\mu(\deg(u))\leq n.$ Let $z$ be the minimal spike in $P_n$ with $\deg(z) = \deg(u).$ Then, if $\omega(u) < \omega(z),$ then $u\in \overline{\mathbb A}P_n.$ 
\end{theo}

The converse of this theorem does not hold in general.

\medskip

We next record the finite-dimensional duality principle that will be used to connect the degree computations below with a distinguished Milnor operation. For $m\geq 0$, put
\[
H_m=(\overline{\mathbb A}P_5)\cap (P_5)_m,
\qquad
(QP_5)_m=(P_5)_m/H_m.
\]
Equip each homogeneous component $(P_5)_m$ with the monomial bilinear form
\[
\langle u_1^{a_1}\cdots u_5^{a_5},u_1^{b_1}\cdots u_5^{b_5}\rangle
=\prod_{i=1}^{5}\delta_{a_i,b_i}.
\]
This is the cohomological form corresponding to the usual divided-power pairing between $H_*(V^{\oplus 5})$ and $H^*(V^{\oplus 5})$. We denote by $H_m^\perp$ the orthogonal complement of the hit subspace $H_m$ with respect to this bilinear form.

\begin{mde}[Orthogonality criterion via the transpose adjoint]\label{prop:adjoint-hit-test}
Let $a\in \mathbb A$ be homogeneous of positive degree $r$, and let
\[
T(a):(P_5)_{m-r}\longrightarrow (P_5)_m
\]
be the degree-$m$ action map. Let
\[
T(a)^*:(P_5)_m\longrightarrow (P_5)_{m-r}
\]
be its transpose adjoint with respect to the monomial bilinear forms. Then
\[
\big((QP_5)_m\big)^*\cong H_m^\perp
\subseteq \Ker(T(a)^*).
\]
Equivalently, if $z\in (P_5)_m$ represents a linear functional on $(QP_5)_m$ under the monomial bilinear form, then $T(a)^*z=0$.
\end{mde}

\begin{proof}
The short exact sequence
\[
0\longrightarrow H_m\longrightarrow (P_5)_m\longrightarrow (QP_5)_m\longrightarrow 0
\]
dualizes to an identification of $((QP_5)_m)^*$ with the subspace of $(P_5)_m^*$ consisting of those linear functionals that vanish on $H_m$. The monomial bilinear form identifies $(P_5)_m$ with $(P_5)_m^*$; under this identification, that subspace is exactly $H_m^\perp$.

Since $a$ has positive degree, $\operatorname{Im}T(a)$ is contained in $H_m$. Hence every element of $H_m^\perp$ vanishes on $\operatorname{Im}T(a)$. For $z\in (P_5)_m$, this condition says
\[
\langle z,T(a)f\rangle=0\qquad \text{for all }f\in (P_5)_{m-r}.
\]
By the defining property of the transpose adjoint, this is equivalent to
\[
\langle T(a)^*z,f\rangle=0\qquad \text{for all }f\in (P_5)_{m-r}.
\]
The monomial bilinear form on $(P_5)_{m-r}$ is nondegenerate, so $T(a)^*z=0$. This proves the asserted inclusion.
\end{proof}

\begin{nxe}\label{rem:homological-right-action}
Let $a\in\mathbb A$ be homogeneous of degree $r$, and let $T(a):(P_5)_{m-r}\to (P_5)_m$ be the corresponding cohomological action map. The right action on $H_*(V^{\oplus 5})$ is the homological, degree-lowering action characterized by
\[
\langle x\cdot a,f\rangle=\langle x,af\rangle,
\qquad x\in H_m(V^{\oplus 5}),\quad f\in (P_5)_{m-r},
\]
under the divided-power pairing. Therefore, after identifying divided-power basis elements with their dual monomials, the matrix of this homological right action is the matrix of the transpose adjoint $T(a)^*$. It is not the cohomological action of $a$ on $P_5$, and it is not the cohomological action of the conjugate operation $\tau(a)$; those cohomological operations raise degree, whereas the homological right action and $T(a)^*$ lower degree.
\end{nxe}

\begin{nxe}\label{rem:code-annihilator-connection}
The preceding proposition is the target-degree side of the code-annihilator identity proved in \cite{PhucCodeAnn}. That companion work starts from the coding theory over Steenrod sub-Hopf algebras introduced by Dougherty--Vergili \cite{DV2017} and from the Vergili--Karaca self-dual code generated by $P_{n+1}^{0}$ \cite{VK2017}. In that setting, the canonical code $\mathbb A(n)[P_{n+1}^{0}]$ over the finite sub-Hopf algebra $\mathbb A(n)$ has source-degree annihilator equal to the kernel of $P_{n+1}^{0}$, and the orthogonal complement of that kernel is the image of the transpose adjoint. In the present paper we use the same transpose-adjoint mechanism inside the rank-five Peterson hit problem, not as a replacement for the admissible-basis calculation, but as an additional canonical obstruction attached to the Milnor operation $P_{t+5}^{0}$.
\end{nxe}

\section{Key outcomes and their respective demonstrations}\label{s3}

As previously noted, this section focuses on investigating the hit problem for the $\mathbb A$-module $P_5$ and the behavior of the Singer transfer of rank 5 in the general degree $d = 2^{t+5} + 2^{t+2} + 2^{t+1} - 5$, where $t\geq 0,\, t\in\mathbb Z.$ 

Before entering the separate computations for $d=33$ and $d=71$, we isolate the new Milnor-operation layer that is present in every degree of this family. Put
\[
d_t=2^{t+5}+2^{t+2}+2^{t+1}-5,
\qquad
D_t=2^{t+5}-1,
\qquad
e_t=d_t-D_t=2^{t+2}+2^{t+1}-4.
\]
Thus $d_t=e_t+D_t$, and $D_t$ is the degree of the Milnor basis element $P_{t+5}^{0}$, which is primitive because it is dual to the polynomial generator $\xi_{t+5}$.

\begin{theo}[The $P_{t+5}^{0}$-hit layer in degree $d_t$]\label{thm:large-milnor-layer}
For every $t\geq 0$, the Milnor basis element $P_{t+5}^{0}$ induces a homogeneous map
\[
T_t=P_{t+5}^{0}:(P_5)_{e_t}\longrightarrow (P_5)_{d_t}.
\]
Its kernel is precisely the Frobenius-square subspace
\[
\Ker(T_t)=\{h^2:h\in (P_5)_{3\cdot 2^t-2}\}.
\]
Consequently,
\[
\dim\Ker(T_t)=\binom{3\cdot 2^t+2}{4},
\qquad
\operatorname{rank}(T_t)=\binom{6\cdot 2^t}{4}-\binom{3\cdot 2^t+2}{4}.
\]
Moreover,
\[
\operatorname{Im}(T_t^*)
=
\left\langle
u_1^{a_1}\cdots u_5^{a_5}\in (P_5)_{e_t}: \text{ at least one } a_i \text{ is odd}
\right\rangle_{\mathbb Z/2}
=
\Ker(T_t)^\perp,
\]
and
\[
H_{d_t}^{\perp}\subseteq \Ker(T_t^*),
\]
where $H_{d_t}=(\overline{\mathbb A}P_5)\cap (P_5)_{d_t}$ and $T_t^*:(P_5)_{d_t}\to (P_5)_{e_t}$ is the transpose adjoint. Under the divided-power duality, the last inclusion says that every primitive homology representative in degree $d_t$ is annihilated by the homological right action of $P_{t+5}^{0}$. This right action is the degree-lowering action determined by the divided-power pairing, not the cohomological action of $P_{t+5}^{0}$ or of $\tau(P_{t+5}^{0})$ on $P_5$.
\end{theo}

\begin{proof}
The Milnor basis element $P_{t+5}^{0}$ is dual to $\xi_{t+5}$ in the Milnor dual algebra \cite{Milnor1958}. Since $\xi_{t+5}$ is an indecomposable polynomial generator of $\mathbb A^*$, the coproduct of the dual basis element $P_{t+5}^{0}$ in $\mathbb A$ is
\[
\Delta(P_{t+5}^{0})=P_{t+5}^{0}\otimes 1+1\otimes P_{t+5}^{0}.
\]
Thus $P_{t+5}^{0}$ acts as a derivation on the polynomial algebra $P_5$. Its action on the degree-one generators is
\[
P_{t+5}^{0}(u_i)=u_i^{2^{t+5}},\qquad 1\leq i\leq 5.
\]
Therefore, for a monomial $u^a=u_1^{a_1}\cdots u_5^{a_5}$ of degree $e_t$, the Cartan formula gives
\[
P_{t+5}^{0}(u^a)=
\sum_{\substack{1\leq i\leq 5\\ a_i\ {\rm odd}}}
 u_1^{a_1}\cdots u_i^{a_i+D_t}\cdots u_5^{a_5}.
\]
Indeed, differentiating the factor $u_i^{a_i}$ contributes the coefficient $a_i$ modulo $2$, and the exponent of $u_i$ changes from $a_i$ to $a_i-1+2^{t+5}=a_i+D_t$.

Since
\[
e_t=6\cdot 2^t-4<2^{t+5}-1=D_t,
\]
every exponent of every source monomial is strictly smaller than $D_t$. Hence each monomial appearing in $P_{t+5}^{0}(u^a)$ has exactly one exponent at least $D_t$, namely the exponent corresponding to the factor on which $P_{t+5}^{0}$ has acted. If
\[
u_1^{a_1}\cdots u_i^{a_i+D_t}\cdots u_5^{a_5}
=u_1^{b_1}\cdots u_j^{b_j+D_t}\cdots u_5^{b_5}
\]
with $|a|=|b|=e_t$ and with $a_i,b_j$ odd, then the unique exponent at least $D_t$ forces $i=j$, and then $a=b$. Thus the nonzero images of distinct source monomials have pairwise disjoint supports in the target monomial basis.

Therefore no cancellation can occur among the images of distinct source monomials. Hence a polynomial
\[
f=\sum_{|a|=e_t}\lambda_a u^a\in (P_5)_{e_t}
\]
lies in $\Ker(T_t)$ if and only if every monomial $u^a$ with $\lambda_a=1$ has zero image under $T_t$. By the displayed formula for $P_{t+5}^{0}(u^a)$, this is equivalent to all exponents $a_i$ being even. Thus $a_i=2b_i$ for all $i$, and
\[
\sum_{i=1}^{5}b_i=\frac{e_t}{2}=3\cdot 2^t-2.
\]
Since the Frobenius endomorphism is $\mathbb Z/2$-linear, the span of the even-exponent monomials in degree $e_t$ is exactly the Frobenius-square subspace
\[
\{h^2:h\in (P_5)_{3\cdot 2^t-2}\}.
\]
Thus
\[
\Ker(T_t)=\{h^2:h\in (P_5)_{3\cdot 2^t-2}\}.
\]
The dimension of $(P_5)_N$ is $\binom{N+4}{4}$, so
\[
\dim\Ker(T_t)=\binom{3\cdot 2^t+2}{4},
\qquad
\dim(P_5)_{e_t}=\binom{e_t+4}{4}=\binom{6\cdot 2^t}{4},
\]
and the rank formula follows from rank-nullity. The monomial bilinear form makes the degree-$e_t$ monomial basis self-dual. Since $\Ker(T_t)$ is spanned by the even-exponent monomials, its orthogonal complement is the span of the degree-$e_t$ monomials with at least one odd exponent. The finite-dimensional identity $(\Ker T_t)^\perp=\operatorname{Im}(T_t^*)$ therefore gives
\[
\operatorname{Im}(T_t^*)
=
\left\langle
u_1^{a_1}\cdots u_5^{a_5}\in (P_5)_{e_t}: \text{ at least one } a_i \text{ is odd}
\right\rangle_{\mathbb Z/2}
=
\Ker(T_t)^\perp.
\]

Finally, $P_{t+5}^{0}$ has positive degree, so $\operatorname{Im}(T_t)\subseteq H_{d_t}$. Proposition~\ref{prop:adjoint-hit-test}, applied with $a=P_{t+5}^{0}$ and $m=d_t$, gives $H_{d_t}^{\perp}\subseteq \Ker(T_t^*)$. Under the canonical pairing between divided powers and monomials, $T_t^*$ is the matrix of the homological right action by $P_{t+5}^{0}$ on $H_*(V^{\oplus 5})$. This is the degree-lowering divided-power action characterized by
\[
\langle x\cdot P_{t+5}^{0},f\rangle=\langle x,P_{t+5}^{0}f\rangle,
\]
and it should not be confused with the cohomological action of $P_{t+5}^{0}$ or of $\tau(P_{t+5}^{0})$ on $P_5$. Hence primitive homology representatives, which annihilate the whole hit subspace, are in particular annihilated by this right action. The theorem follows.
\end{proof}

\begin{nxe}\label{rem:large-milnor-companion-instance}
With $k=5$, $n=t+4$, and $d=e_t$, the equality
\[
\operatorname{Im}(T_t^*)=\Ker(T_t)^\perp
\]
is the degree-specific instance of the companion identity $(\Ker T_d)^\perp=\operatorname{Im}(T_d^*)$ proved in our previous work \cite{PhucCodeAnn}. The inclusion
\[
H_{d_t}^{\perp}\subseteq \Ker(T_t^*)
\]
is the corresponding target-degree hit-quotient consequence. Thus the present five-variable computation uses the companion transpose-adjoint formalism without claiming that the single operation $P_{t+5}^{0}$ determines the whole hit subspace.
\end{nxe}

\begin{hqu}\label{cor:large-milnor-lower-bound}
For every $t\geq 0$, the degree-$d_t$ hit subspace satisfies
\[
\dim H_{d_t}\geq \binom{6\cdot 2^t}{4}-\binom{3\cdot 2^t+2}{4}.
\]
Equivalently, the admissible-basis computation in degree $d_t$ must occur after removing at least this canonical $P_{t+5}^{0}$-hit subspace. This subspace is natural for the invariant theory because Steenrod operations commute with the standard $GL_5$-action on $P_5$.
\end{hqu}

\begin{proof}
The image of $T_t=P_{t+5}^{0}$ is contained in $H_{d_t}$, and Theorem~\ref{thm:large-milnor-layer} gives its rank.
\end{proof}

We will now study the structure of the spaces $(QP_5)_{d_t}$ and their invariant subspaces. The theorem above gives an additional orthogonal verification tool for the primitive divided-power representatives used in the transfer calculation, while the admissible-basis computations below determine the quotient and the $GL_5$-invariants explicitly. Before presenting those invariant calculations, we make explicit how the computations are certified as finite-dimensional linear algebra over $\mathbb Z/2$.

\medskip

{\bf Case \mbox{\boldmath $t = 0.$}}\ We have $d = 2^{0+5} + 2^{0+2} + 2^{0+1} - 5 = 33.$ 
 We observe that since the Kameko squaring operation $(\widetilde {Sq^0_*})_{33}: (QP_5)_{33}\longrightarrow (QP_5)_{14}$ is an epimorphism, we have an isomorphism $$(QP_5)_{33}\cong {\rm Ker}((\widetilde {Sq^0_*})_{33})\bigoplus (QP_5)_{14}.$$ 
As it is known, there exist exactly $320$ admissible monomials of degree $14$ in the $\mathbb A$-module $P_5.$ Consequently, $(QP_5)_{14}$ has dimension $320$ (see also \cite{D.P5}). Thus, in order to obtain a basis for the space $(QP_5)_{33}$, we must calculate the dimension of ${\rm Ker}((\widetilde {Sq^0_*})_{33}).$ In this regard, the following remark may be of interest.

\begin{nxe}\label{nx33-1}
We claim that if $g$ is an admissible monomial of degree $33$ in $P_5$ such that $[g]$ belongs to ${\rm Ker}((\widetilde {Sq^0_*})_{33}),$ then $\omega_1(g) = 3.$ Indeed, since $u_1^{31}u_2u_3$ is the minimal spike of degree $33$ in the $\mathbb A$-module $P_5$ and $\omega(u_1^{31}u_2u_3) = (3,1,1,1,1),$ according to Lemma \ref{dlPS}, $u_1^{31}u_2u_3$ is admissible. Because $g$ is admissible, and $\deg(g)$ is odd, by Theorem \ref{dlSinger}, either $\omega_1(g) = 3$ or $\omega_1(g) = 5.$ If $\omega_1(g) = 5,$ then $g$ is of the form $\prod_{1\leq j\leq 5}u_jh^{2},$ in which $h$ is a monomial of degree $14$ in $P_5.$ Since $\prod_{1\leq j\leq 5}u_jh^{2}$ is admissible, following Theorem \ref{dlKS}, $h$ is also admissible. So, $(\widetilde {Sq^0_*})_{33}([\prod_{1\leq j\leq 5}u_jh^{2}]) = [h]\neq [0],$ which contradicts with the fact that $[g] = [\prod_{1\leq j\leq 5}u_jh^{2}]\in {\rm Ker}((\widetilde {Sq^0_*})_{33}).$ Thus, we must have  $\omega_1(g) = 3,$ and therefore by the above arguments, $g$ can be represented as $u_iu_ju_kh^{2}$, where $1\leq i<j<k\leq 5,$ and $h$ is an admissible monomial of degree $15$ in $P_5.$ According to Sum \cite{N.S0}, $\omega(h)$ belongs to the set $\{(1,1,1,1),\ (1,1,3),\ (3,2,2)\ (3,4,1)\},$ and therefore, we may conclude that the weight vector $\omega(g)$ of $g$ is one of the following sequences: 
$$ \widetilde{\omega}_{(1)}:= (3,1,1,1,1),\ \widetilde{\omega}_{(2)}:= (3,1,1,3),\ \widetilde{\omega}_{(3)}:= (3,3,2,2),\ \widetilde{\omega}_{(4)}:= (3,3,4,1).$$ Thus, we have an isomorphism
 $$ \begin{array}{ll}
\medskip
{\rm Ker}((\widetilde {Sq^0_*})_{33})&\cong (QP_5)_{33}^0\bigoplus \big({\rm Ker}((\widetilde {Sq^0_*})_{33})\bigcap  (QP_5)_{33}^{>0}\big)\\
&\cong (QP_5)_{33}^0\bigoplus \big(\bigoplus_{1\leq i\leq 4}(QP_5)_{33}^{>0}(\widetilde{\omega}_{(i)})\big).
\end{array}$$  
We now compute the subspaces $(QP_5)_{33}^0$ and $(QP_5)_{33}^{>0}(\widetilde{\omega}_{(i)}).$
\end{nxe}

{\bf Computation of \mbox{\boldmath $(QP_5)_{33}^0.$}} From \cite{W.W}, we have
\begin{eqnarray}
\medskip
\dim (QP_5)_{33}^{0}(\widetilde{\omega}_{(i)}) &=& \sum_{3\leq s\leq 4}\binom{5}{s}\dim (QP_s)_{33}^{>0}(\widetilde{\omega}_{(i)}),\label{pt2}\\
\dim (QP_5)_{33}^{0} &=& \sum_{3\leq s\leq 4}\binom{5}{s}\dim (QP_s)_{33}^{>0}\label{pt3}.
\end{eqnarray} 
Since $\mu(33) = 3,$ $(QP_3)_{33}$ is isomorphic to $(QP_3)_{15}.$ So $\dim (QP_3)_{33}^{>0}= \dim (QP_3)_{15}  =13$ (see Kameko \cite{M.K}). Following Sum \cite{N.S1}, $\dim (QP_4)^{>0}_{33}= 84,$ and so, by \eqref{pt3}, we get
$$ \dim (QP_5)^{0}_{33} = 13\binom{5}{3} + 84\binom{5}{4} = 550.$$

According to Sum \cite{N.S1}, $(QP_4)_{33}\cong (QP_4)_{33}(\widetilde{\omega}_{(1)})\bigoplus (QP_4)_{33}( \widetilde{\omega}_{(3)}),$ which shows that 
$(QP_5)^{0}_{33}(\widetilde{\omega}_{(2)}) = (QP_5)^{0}_{33}( \widetilde{\omega}_{(4)}) = 0,$ and
$$ (QP_5)^{0}_{33}\cong (QP_5)^{0}_{33}(\widetilde{\omega}_{(1)}) \bigoplus (QP_5)^{0}_{33}( \widetilde{\omega}_{(3)}).$$
Moreover, by Kameko \cite{M.K}, we have $\dim (QP_3)_{33}^{>0}(\widetilde{\omega}_{(1)}) = 7$ and $\dim  (QP_3)_{33}^{>0}(\widetilde{\omega}_{(3)}) = 6.$ By Sum \cite{N.S1}, we have $\dim  (QP_4)_{33}^{>0}(\widetilde{\omega}_{(1)})  = 17$ and $\dim  (QP_4)_{33}^{>0}(\widetilde{\omega}_{(3)}) = 67.$ So, using \eqref{pt2}, we derive:
$$ \begin{array}{ll}
\medskip
 \dim (QP_5)^{0}_{33}(\widetilde{\omega}_{(1)}) &= 7\binom{5}{3} + 17\binom{5}{4} = 155,\\
 \dim (QP_5)^{0}_{33}(\widetilde{\omega}_{(3)}) &= 6\binom{5}{3} + 67\binom{5}{4} = 395.
\end{array}$$


{\bf Computation of \mbox{\boldmath $(QP_5)_{33}^{>0}(\widetilde{\omega}_{(i)}).$}} We obtain
\begin{theo}\label{dl33-2}
The following statements are true:
$$ 
\dim (QP_5)^{>0}_{33}(\widetilde{\omega}_{(i)}) =\left\{\begin{array}{ll} 
\medskip
31&\mbox{if $i = 1$},\\
\medskip
0&\mbox{if $i = 2,\, 4$},\\
\medskip
421&\mbox{if $i = 3$}.
\end{array}\right.$$
\end{theo}

\begin{proof}

We prove the theorem for the cases where $i= 1,\, 2,\, 4.$ The remaining case $i=3$ is obtained by the same admissible-basis reduction implemented in Appendix~\ref{s4}, which returns $421$ independent classes after reduction modulo hit and lower-weight terms.

\medskip

$\bullet$  {\it Case $i  =1$}. According to Moetele and Mothebe \cite[Proposition 1]{M.M2}, if $t$ is an admissible monomial in the $\mathbb A$-module $P_n,$ then $t$ must take the form $t = u_1^{2^{a_1} - 1}u_2^{a_2}u_3^{a_3}\ldots u_n^{a_n}.$ By this, we only need to consider the monomials of the form $t:= u_iu_ju_k\underline{t}^{2} = \prod_{1\leq l\leq 5}u_l^{\alpha_l},$ with $\underline{t}\in \mathscr C^{\otimes 5}_{15},\, 1\leq i< j<k\leq 5,$ $\alpha_1\in \{1, 3, 7, 15\},$ and $\prod_{2\leq l\leq 5}\alpha_l> 0.$ A simple computation indicates that if $t\neq {\rm adm}_i,\, 156\leq i\leq 186,$ where 

\begin{center}
\begin{tabular}{lll}
${\rm adm}_{156}=u_1u_2u_3u_4^{2}u_5^{28}$, & ${\rm adm}_{157}=u_1u_2u_3^{2}u_4u_5^{28}$, & ${\rm adm}_{158}=u_1u_2u_3^{2}u_4^{28}u_5$, \\
${\rm adm}_{159}=u_1u_2^{2}u_3u_4u_5^{28}$, & ${\rm adm}_{160}=u_1u_2^{2}u_3u_4^{28}u_5$, & ${\rm adm}_{161}=u_1u_2^{2}u_3^{28}u_4u_5$, \\
${\rm adm}_{162}=u_1u_2u_3^{2}u_4^{4}u_5^{25}$, & ${\rm adm}_{163}=u_1u_2^{2}u_3u_4^{4}u_5^{25}$, & ${\rm adm}_{164}=u_1u_2^{2}u_3^{4}u_4u_5^{25}$, \\
${\rm adm}_{165}=u_1u_2^{2}u_3^{4}u_4^{25}u_5$, & ${\rm adm}_{166}=u_1u_2u_3^{2}u_4^{5}u_5^{24}$, & ${\rm adm}_{167}=u_1u_2^{2}u_3u_4^{5}u_5^{24}$, \\
${\rm adm}_{168}=u_1u_2^{2}u_3^{5}u_4u_5^{24}$, & ${\rm adm}_{169}=u_1u_2^{2}u_3^{5}u_4^{24}u_5$, & ${\rm adm}_{170}=u_1u_2u_3^{3}u_4^{4}u_5^{24}$, \\
${\rm adm}_{171}=u_1u_2^{3}u_3u_4^{4}u_5^{24}$, & ${\rm adm}_{172}=u_1^{3}u_2u_3u_4^{4}u_5^{24}$, & ${\rm adm}_{173}=u_1u_2^{3}u_3^{4}u_4u_5^{24}$, \\
${\rm adm}_{174}=u_1^{3}u_2u_3^{4}u_4u_5^{24}$, & ${\rm adm}_{175}=u_1u_2^{3}u_3^{4}u_4^{24}u_5$, & ${\rm adm}_{176}=u_1^{3}u_2u_3^{4}u_4^{24}u_5$, \\
${\rm adm}_{177}=u_1u_2^{2}u_3^{4}u_4^{9}u_5^{17}$, & ${\rm adm}_{178}=u_1u_2^{2}u_3^{5}u_4^{8}u_5^{17}$, & ${\rm adm}_{179}=u_1u_2^{3}u_3^{4}u_4^{8}u_5^{17}$, \\
${\rm adm}_{180}=u_1^{3}u_2u_3^{4}u_4^{8}u_5^{17}$, & ${\rm adm}_{181}=u_1u_2^{2}u_3^{5}u_4^{9}u_5^{16}$, & ${\rm adm}_{182}=u_1u_2^{3}u_3^{4}u_4^{9}u_5^{16}$, \\
${\rm adm}_{183}=u_1^{3}u_2u_3^{4}u_4^{9}u_5^{16}$, & ${\rm adm}_{184}=u_1u_2^{3}u_3^{5}u_4^{8}u_5^{16}$, & ${\rm adm}_{185}=u_1^{3}u_2u_3^{5}u_4^{8}u_5^{16}$, \\
${\rm adm}_{186}=u_1^{3}u_2^{5}u_3u_4^{8}u_5^{16},$ 
\end{tabular}
\end{center}
then either $t\in \{u_1^{3}u_2^{5}u_3^{8}u_4u_5^{16}, u_1^{3}u_2^{5}u_3^{8}u_4^{16}u_5\},$ or $t$ is of the form $t'z^{16},$ where $t'$ is one of the following inadmissible monomials: $$u_1^{3}u_2^{4}u_3u_4^{9},\, u_1^{3}u_2^{4}u_3^{9}u_4,\, u_1^{3}u_2^{4}u_3^{9}u_5,\, u_1^{3}u_2^{4}u_3u_4u_5^{8},\, u_1^{3}u_2^{4}u_3u_4^{8}u_5,\, u_1^{3}u_2^{4}u_3^{8}u_4u_5.$$ 
Applying Cartan's formula, we obtain the following equality:
$$ \begin{array}{ll}
\medskip
 u_1^{3}u_2^{5}u_3^{8}u_4u_5^{16} &= Sq^{4}\big(u_1^{3}u_2^{3}u_3^{8}u_4u_5^{14} + u_1^{3}u_2^{3}u_3^{4}u_4u_5^{18} + u_1^{3}u_2^{3}u_3^{2}u_4u_5^{20}+ u_1^{3}u_2^{3}u_3u_4^{2}u_5^{20}\\
\medskip
&\quad+ u_1^{3}u_2^{3}u_3u_4^{4}u_5^{18} + u_1^{3}u_2^{3}u_3u_4^{8}u_5^{14}\big) +Sq^{2}\big(u_1^{5}u_2^{3}u_3^{8}u_4u_5^{14}+ u_1^{5}u_2^{3}u_3^{4}u_4u_5^{18}\\
\medskip
&\quad+ u_1^{5}u_2^{3}u_3^{2}u_4u_5^{20}+ u_1^{5}u_2^{3}u_3u_4^{8}u_5^{14}+ u_1^{5}u_2^{3}u_3u_4^{4}u_5^{18}+ u_1^{5}u_2^{3}u_3u_4^{2}u_5^{20} + u_1^{2}u_2^{3}u_3u_4u_5^{24}\big)\\
&\quad + Sq^{1}\big(u_1^{3}u_2^{3}u_3u_4u_5^{24}\big) + u_1^{3}u_2^{5}u_3u_4^{8}u_5^{16} +  u_1^{3}u_2^{4}u_3u_4u_5^{24} + u_1^{2}u_2^{5}u_3u_4u_5^{24} + \sum X,
\end{array}$$
where $\omega(X) < \omega(u_1^{3}u_2^{5}u_3^{8}u_4u_5^{16}) = \widetilde{\omega}_{(1)}.$ Thus the monomial $u_1^{3}u_2^{5}u_3^{8}u_4u_5^{16}$ is inadmissible, and so is $u_1^{3}u_2^{5}u_3^{8}u_4^{16}u_5$. On the other side, we observe that $\omega_4(t') = 1\neq 0,$ and $\omega_j(t') = 0$ for all $j> 4$. So, following Theorem \ref{dlKS}(ii), $t= t'z^{16}$ is inadmissible. Therefore, it follows that $$(QP_5)_{33}^{>0}(\widetilde{\omega}_{(1)}) = \langle \{[{\rm adm}_i]_{\omega_{(1)}} = [{\rm adm}_i]|\, 156\leq i\leq 186\}\rangle.$$ 
Let us demonstrate that the set $\big\{[{\rm adm}_i] \mid 156\leq i\leq 186\big\}$ is linearly independent in $(QP_5)_{33}^{>0}(\widetilde{\omega}_{(1)})$. To do so, consider the set
$$ \mathcal N_5 := \{(l, \mathscr L)\;|\; \mathscr L = (l_1,l_2,\ldots, l_r), 1\leqslant l < l_1< l_2 < \ldots < l_r\leq 5, \ 0\leq r \leq 4\},$$  
where by convention, $\mathscr L = \emptyset,$ if $r = 0.$ Denote by $r = \ell(\mathscr L)$ the length of $\mathscr L.$ For a pair $(l, \mathscr L)\in\mathcal N_5$ with $r = \ell(\mathscr L) = 1$, $\mathscr L = (l')$, and $1\leq l< l'\leq 5.$ We define an $\mathbb A$-homomorphism $\mathsf{p}_{(l, l')} : P_5\to P_4$, as follows:
$$ \mathsf{p}_{(l, l')}(u_j) = \left\{ \begin{array}{ll}
\medskip
{u_j}&\text{if }\;1\leq j \leq l-1, \\
\medskip
u_{l'-1}& \text{if}\; j = l,\\
u_{j-1}&\text{if}\; l+1 \leq j \leq 5.
\end{array} \right.$$
We assume that there is a linear relation $S = \sum_{156\leq i\leq 186}\gamma_i{\rm adm}_i\sim  0,$ in which $\gamma_i\in \mathbb Z/2$ for every $i.$ According to Sum \cite{N.S1}, $(QP_{33}^{\otimes 4})^{>0}(\widetilde{\omega}_{(1)}) = \langle \{[{\rm Adm}_j]:\ 1\leq j\leq 17\}\rangle,$ where the admissible monomials ${\rm Adm}_j$ are determined as follows:

\begin{center}
\begin{tabular}{lrrr}
${\rm Adm}_{1}=u_1u_2^{2}u_3^{5}u_4^{25}$, & \multicolumn{1}{l}{${\rm Adm}_{2}=u_1u_2^{3}u_3^{4}u_4^{25}$,} & \multicolumn{1}{l}{${\rm Adm}_{3}=u_1^{3}u_2u_3^{4}u_4^{25}$,} & \multicolumn{1}{l}{${\rm Adm}_{4}=u_1u_2^{3}u_3^{5}u_4^{24}$,} \\
${\rm Adm}_{5}=u_1^{3}u_2u_3^{5}u_4^{24}$, & \multicolumn{1}{l}{${\rm Adm}_{6}=u_1^{3}u_2^{5}u_3u_4^{24}$,} & \multicolumn{1}{l}{${\rm Adm}_{7}=u_1u_2u_3^{2}u_4^{29}$,} & \multicolumn{1}{l}{${\rm Adm}_{8}=u_1u_2^{2}u_3u_4^{29}$,} \\
${\rm Adm}_{9}=u_1u_2^{2}u_3^{29}u_4$, & \multicolumn{1}{l}{${\rm Adm}_{10}=u_1u_2u_3u_4^{30}$,} & \multicolumn{1}{l}{${\rm Adm}_{11}=u_1u_2u_3^{30}u_4$,} & \multicolumn{1}{l}{${\rm Adm}_{12}=u_1u_2^{30}u_3u_4$,} \\
${\rm Adm}_{13}=u_1u_2u_3^{3}u_4^{28}$, & \multicolumn{1}{l}{${\rm Adm}_{14}=u_1u_2^{3}u_3u_4^{28}$,} & \multicolumn{1}{l}{${\rm Adm}_{15}=u_1u_2^{3}u_3^{28}u_4$,} & \multicolumn{1}{l}{${\rm Adm}_{16}=u_1^{3}u_2u_3u_4^{28}$,} \\
${\rm Adm}_{17}=u_1^{3}u_2u_3^{28}u_4$. &       &       &  
\end{tabular}%
\end{center}

From this result, we express each $\mathsf{p}_{(l, l')}(S)$ in terms of the classes $[\mathrm{Adm}_j]$ modulo $\overline{\mathbb A}P_4$. The relations $\mathsf{p}_{(l, l')}(S)\sim 0$, for $1\leq l\leq 3$ and $2\leq l'\leq 4$, then form a linear system in the coefficients $\gamma_i$. Row reduction of this system gives $\gamma_i=0$ for every $i$ with $156\leq i\leq 186.$

\medskip

$\bullet$ {\it Case $i = 2$}. Let us recall that by Remark \ref{nx33-1}, the monomials of degree $33$ has form $u_iu_ju_k\underline{t}^{2}$ with $1\leq i<j<k\leq 5$ and $\underline{t}\in \mathscr C^{\otimes 5}_{15}.$ Following \cite{N.S0}, $|\mathscr C^{\otimes 5}_{15}| = 432.$ Thanks to this result, we find that if $t\in (P_5)_{33}(\widetilde{\omega}_{(2)}),$ then $t$ is one of the following monomials: $u_1^{3}u_2^{4}u_3^{8}u_4^{9}u_5^{9},$ \ $u_1^{3}u_2^{4}u_3^{9}u_4^{8}u_5^{9},$\ $u_1^{3}u_2^{4}u_3^{9}u_4^{9}u_5^{8},$\ $u_1^{3}u_2^{5}u_3^{8}u_4^{8}u_5^{9},$\ $u_1^{3}u_2^{5}u_3^{8}u_4^{9}u_5^{8},$\ $u_1^{3}u_2^{5}u_3^{9}u_4^{8}u_5^{8}.$ Let us consider the monomials $u_1^{3}u_2^{4}u_3^{8}u_4^{9}u_5^{9}$ and $u_1^{3}u_2^{5}u_3^{8}u_4^{8}u_5^{9}.$ Applying the Cartan formula, we get
$$ \begin{array}{ll}
\medskip
u_1^{3}u_2^{4}u_3^{8}u_4^{9}u_5^{9}&= Sq^{1}(u_1^{3}u_2u_3^{10}u_4^{9}u_5^{9})\\
\medskip
&\quad + Sq^{2}(u_1^{2}u_2u_3^{10}u_4^{9}u_5^{9} + u_1^{5}u_2u_3^{6}u_4^{9}u_5^{10} + u_1^{5}u_2u_3^{6}u_4^{10}u_5^{9} + u_1^{5}u_2^{2}u_3^{6}u_4^{9}u_5^{9})\\
\medskip
&\quad+ Sq^{4}(u_1^{3}u_2^{2}u_3^{6}u_4^{9}u_5^{9} + u_1^{3}u_2u_3^{6}u_4^{9}u_5^{10} + u_1^{3}u_2u_3^{6}u_4^{10}u_5^{9})\\
\medskip
&\quad + u_1^{3}u_2u_3^{8}u_4^{9}u_5^{12} +u_1^{3}u_2u_3^{8}u_4^{12}u_5^{9} + u_1^{2}u_2u_3^{12}u_4^{9}u_5^{9} \mod ((P_5)_{33}(< \widetilde{\omega}_{(2)})),\\
\medskip
u_1^{3}u_2^{5}u_3^{8}u_4^{8}u_5^{9}&= Sq^{1}(u_1^{3}u_2^{3}u_3u_4^{8}u_5^{17}) + Sq^{2}(u_1^{5}u_2^{3}u_3u_4^{8}u_5^{14} + u_1^{5}u_2^{3}u_3^{2}u_4^{4}u_5^{17}) \\
\medskip
&\quad + Sq^{4}(u_1^{3}u_2^{3}u_3u_4^{8}u_5^{14} + u_1^{3}u_2^{3}u_3^{2}u_4^{4}u_5^{17} + u_1^{3}u_2^{9}u_3^{4}u_4^{4}u_5^{9})\\
&\quad+ Sq^{8}(u_1^{3}u_2^{5}u_3^{4}u_4^{4}u_5^{9}) \mod ((P_5)_{33}(< \widetilde{\omega}_{(2)})).
\end{array}$$
These calculations show that $u_1^{3}u_2^{4}u_3^{8}u_4^{9}u_5^{9}$ is strictly inadmissible and that $u_1^{3}u_2^{5}u_3^{8}u_4^{8}u_5^{9}$ is $\widetilde{\omega}_{(2)}$-hit. Similarly, the monomials $u_1^{3}u_2^{4}u_3^{9}u_4^{8}u_5^{9},$\ $u_1^{3}u_2^{4}u_3^{9}u_4^{9}u_5^{8}$ are strictly inadmissible and $u_1^{3}u_2^{5}u_3^{8}u_4^{9}u_5^{8},$\ $u_1^{3}u_2^{5}u_3^{9}u_4^{8}u_5^{8}$ are $\widetilde{\omega}_{(2)}$-hit. Thus, $(QP_5)_{33}^{>0}(\widetilde{\omega}_{(2)})$ is trivial.

\medskip

$\bullet$  {\it Case $i = 4$}.  Let $t$ be a monomial in $(P_5)_{33}(\widetilde{\omega}_{(4)}).$ Since $t = u_iu_ju_k\underline{t}^{2}$ where $1<i<j<k\leq 5$ and $\underline{t}\in \mathscr C^{\otimes 5}_{15},$ from a result in Sum \cite{N.S0}, we find that $t$ is the permutation of each of the following monomial:
\begin{center}
\begin{tabular}{lllrr}
$u_1^{2}u_2^{4}u_3^{5}u_4^{7}u_5^{15}$, & $u_1^{2}u_2^{4}u_3^{7}u_4^{7}u_5^{13}$, & $u_1^{2}u_2^{5}u_3^{5}u_4^{6}u_5^{15}$, & \multicolumn{1}{l}{$u_1^{2}u_2^{5}u_3^{5}u_4^{7}u_5^{14}$,} & \multicolumn{1}{l}{$u_1^{2}u_2^{5}u_3^{6}u_4^{7}u_5^{13}$,} \\
$u_1^{2}u_2^{5}u_3^{7}u_4^{7}u_5^{12}$, & $u_1^{3}u_2^{4}u_3^{4}u_4^{7}u_5^{15}$, & $u_1^{3}u_2^{4}u_3^{5}u_4^{6}u_5^{15}$, & \multicolumn{1}{l}{$u_1^{3}u_2^{4}u_3^{5}u_4^{7}u_5^{14}$,} & \multicolumn{1}{l}{$u_1^{3}u_2^{4}u_3^{6}u_4^{7}u_5^{13}$,} \\
$u_1^{3}u_2^{4}u_3^{7}u_4^{7}u_5^{12}$, & $u_1^{3}u_2^{5}u_3^{5}u_4^{6}u_5^{14}$, & $u_1^{3}u_2^{5}u_3^{6}u_4^{6}u_5^{13}$, & \multicolumn{1}{l}{$u_1^{3}u_2^{5}u_3^{6}u_4^{7}u_5^{12}$,} & \multicolumn{1}{l}{$u_1^{4}u_2^{4}u_3^{7}u_4^{7}u_5^{11}$,} \\
$u_1^{4}u_2^{5}u_3^{6}u_4^{7}u_5^{11}$, & $u_1^{4}u_2^{5}u_3^{7}u_4^{7}u_5^{10}$, & $u_1^{5}u_2^{5}u_3^{6}u_4^{7}u_5^{10}$. &       & 
\end{tabular}%
\end{center} 

It is straightforward to check that these monomials are inadmissible. For instance, using the Cartan formula, we find that 
\begin{align*}
u_1^{3}u_2^{4}u_3^{4}u_4^{7}u_5^{15} &= Sq^{1}(u_1^{3}u_2u_3^{2}u_4^{11}u_5^{15} + u_1^{3}u_2u_3^{2}u_4^{7}u_5^{19}) + Sq^{2}(u_1^{5}u_2^{2}u_3^{2}u_4^{7}u_5^{15})\\
&\quad + Sq^{4}(u_1^{3}u_2^{2}u_3^{2}u_4^{7}u_5^{15}) + \mbox{smaller monomials}.
\end{align*} 
Therefore $(QP_5)_{33}^{>0}(\widetilde{\omega}_{(4)})$ is trivial. The theorem is completely proved.
\end{proof}

The computational outcomes imply the following corollary.

\begin{hqu}
There exist exactly $1322$ admissible monomials of degree $33$ in the $\mathbb A$-module $P_5.$ Consequently, $(QP_5)_{33}$ is $1322$-dimensional.
\end{hqu}

Applying the above data, we can explicitly determine the $GL_5$-invariants $((QP_5)_{33})^{GL_5}$. By using Theorem \ref{dl33-2}, we obtain an isomorphism that can be used to compute the desired this invariant space:
$$(QP_5)_{33}\cong {\rm Ker}((\widetilde {Sq^0_*})_{33})\bigoplus (QP_5)_{14} \cong (QP_5)_{33}(\widetilde{\omega}_{(1)})\bigoplus (QP_5)_{33}(\widetilde{\omega}_{(3)})\bigoplus (QP_5)_{14}.$$
This shows that 
\begin{equation}\label{pt3inv}
 \dim ((QP_5)_{33})^{GL_5}\leq \dim ((QP_5)_{33}(\widetilde{\omega}_{(1)}))^{GL_5} + \dim ((QP_5)_{33}(\widetilde{\omega}_{(3)}))^{GL_5} + \dim ((QP_5)_{14})^{GL_5}.
\end{equation}
Thus, to obtain the dimension of $((QP_5)_{33})^{GL_5},$ we need to determine the $GL_5$-invariants $((QP_5)_{33}(\widetilde{\omega}_{(i)}))^{GL_5}$ for $i  =1,\, 3$ as well as $((QP_5)_{14})^{GL_5}.$ The description of $((QP_5)_{14})^{GL_5}$ can be found in our previous work \cite{D.P7}, where the following theorem is stated.
\begin{theo}\label{dlPhuc}
The space $((QP_5)_{14})^{GL_5}$ has dimension one. Furthermore, $$((QP_5)_{14})^{GL_5} = \langle [\zeta(u_1, \ldots, u_5)] \rangle,$$ wherein 
$$ \begin{array}{ll}
 \zeta(u_1, \ldots, u_5) &= u_2u_3u_4^{6}u_5^{6}+
u_1u_3u_4^{6}u_5^{6}+
u_1u_2u_4^{6}u_5^{6}+
u_1u_2u_3^{6}u_5^{6}+
\medskip
u_1u_2u_3^{6}u_4^{6}\\
&\quad+
u_2^{3}u_3^{3}u_4^{4}u_5^{4}+
u_1^{3}u_3^{3}u_4^{4}u_5^{4}+
u_1^{3}u_2^{3}u_4^{4}u_5^{4}+
u_1^{3}u_2^{3}u_3^{4}u_5^{4}+
\medskip
u_1^{3}u_2^{3}u_3^{4}u_4^{4}\\
&\quad+
u_1u_2u_3^{2}u_4^{4}u_5^{6}+
u_1u_2u_3^{2}u_4^{6}u_5^{4}+
u_1u_2u_3^{6}u_4^{2}u_5^{4}+
\medskip
u_1u_2^{2}u_3u_4^{4}u_5^{6}\\
&\quad+
u_1u_2^{2}u_3^{4}u_4u_5^{6}+
u_1u_2^{6}u_3u_4^{2}u_5^{4}+
u_1u_2^{2}u_3^{3}u_4^{4}u_5^{4}+
\medskip
u_1u_2^{3}u_3^{2}u_4^{4}u_5^{4}\\
&\quad+
u_1^{3}u_2u_3^{2}u_4^{4}u_5^{4}+
u_1^{3}u_2u_3^{4}u_4^{2}u_5^{4}+
u_1u_2^{2}u_3^{4}u_4^{3}u_5^{4}+
u_1^{3}u_2^{4}u_3u_4^{2}u_5^{4}.
\end{array}$$
\end{theo}

\begin{theo}\label{dl33-3}
The invariants $((QP_5)_{33}(\widetilde{\omega}_{(1)}))^{GL_5}$ and $((QP_5)_{33}(\widetilde{\omega}_{(3)}))^{GL_5}$ vanish.
\end{theo}

\begin{proof}

Before we proceed, we must define the action of $GL_5$ on the polynomial algebra $P_5$. It is well-established that $GL_5$ is generated by a specific set of operators. We define these generators as the set of $\mathbb{A}$-algebra homomorphisms $\{\rho_1, \ldots, \rho_5\}$, where their actions on the generators $u_i$ of $P_5$ are given as follows:

\begin{enumerate}
    \item \textbf{The Adjacent Transpositions ($\rho_j$ for $1 \le j \le 4$):} These 4 operators generate the symmetric group $\Sigma_5 \subset GL_5$. Each $\rho_j$ swaps the adjacent variables $u_j$ and $u_{j+1}$ and fixes all others:
    $$\rho_j(u_i) = \begin{cases} 
        u_{j+1} & \text{if } i=j, \\ 
        u_j & \text{if } i=j+1. \\ 
        u_i & \text{otherwise}. 
    \end{cases}$$
    This set generates $\Sigma_5$ because any permutation can be decomposed into a product of adjacent transpositions.

    \item \textbf{The Transvection ($\rho_5$):} The final generator needed to obtain the full $GL_5$ group is the transvection $\rho_5$. This operator adds the variable $u_4$ to $u_5$ and fixes all other variables:
    $$\rho_5(u_i) = \begin{cases} 
        u_5 + u_4 & \text{if } i=5, \\ 
        u_i & \text{if } i < 5. 
    \end{cases}$$
\end{enumerate}

This set $\{\rho_1, \ldots, \rho_5\}$ successfully generates $GL_5$. The set $\{\rho_1, \ldots, \rho_4\}$ generates all permutation matrices (Type I elementary operations). The operator $\rho_5$ is an elementary transvection (a Type III operation, $C_5 \to C_5 + C_4$). All other transvections $C_i \to C_i + C_j$ can be generated by conjugating $\rho_5$ with permutation matrices, which are themselves generated by $\{\rho_1, \ldots, \rho_4\}$.

Let $f = u_1^{a_1}u_2^{a_2}\ldots u_5^{a_5}$ be a monomial in $(P_5)_{33}.$ The weight vector $\omega(f)$ is invariant under the permutation of the generators $u_j,\ j = 1, 2, \ldots, 5;$ hence $(QP_5)_{33}(\widetilde{\omega}_{(i)})$ also has a $\Sigma_5$-module structure. The linear maps $\rho_j$ induce homomorphisms of $\mathbb{A}$-algebras, also denoted by $\rho_j: P_5\longrightarrow P_5.$ Therefore, a class $[f]_{\widetilde{\omega}_{(i)}}\in (QP_5)_{33}(\widetilde{\omega}_{(i)})$ is a \textbf{$GL_5$-invariant} if and only if it is invariant under all generators:
$$ \rho_j(f) \sim_{\widetilde{\omega}_{(i)}} f \quad \text{for all } 1 \le j \le 5. $$
If this condition holds only for $1 \le j \le 4$, then $[f]_{\widetilde{\omega}_{(i)}}$ is an \textbf{$\Sigma_5$-invariant}.

\begin{nxe}[Machine-certified linear algebra for invariants]\label{rem:machine-certified-invariants}
Each invariant computation below is reduced to a finite system of linear equations over $\mathbb Z/2$. Let $W$ be one of the quotient spaces under consideration, equipped with the admissible basis supplied either in the text or in the data files cited in Appendix~\ref{s4}. For each generator $\rho_i$, the algorithm expands $\rho_i(b)-b$ modulo hit terms and lower-weight terms, and then records its coordinate column in the chosen admissible basis. The $\Sigma_5$-invariant space is the kernel of the block matrix
\[
M_W^{\Sigma}=\begin{bmatrix}
[\rho_1-1]_W\\
[\rho_2-1]_W\\
[\rho_3-1]_W\\
[\rho_4-1]_W
\end{bmatrix},
\]
and the $GL_5$-invariant space is the kernel of
\[
M_W^{GL}=\begin{bmatrix}
[\rho_1-1]_W\\
[\rho_2-1]_W\\
[\rho_3-1]_W\\
[\rho_4-1]_W\\
[\rho_5-1]_W
\end{bmatrix}.
\]
Equivalently, after a basis of $\Ker(M_W^\Sigma)$ has been found, the final passage from $\Sigma_5$-invariants to $GL_5$-invariants is obtained by row reducing the matrix of $\rho_5-1$ restricted to this kernel. All ranks reported below are the ranks of these explicit matrices over $\mathbb Z/2$, computed by the \texttt{OSCAR} implementation and the Zenodo data in Appendix~\ref{s4}.
\end{nxe}

\medskip

In what follows, let $f_1, f_2, \ldots, f_s$ be monomials in $(P_5)_{33}(\widetilde{\omega}_{(i)}).$ For a subgroup $G$ of $GL_5$ denote by $G(f_1; f_2; \ldots, f_s)$ the $G$-submodule of $(QP_5)_{33}(\widetilde{\omega}_{(i)})$ generated by the set $\{[f_i]_{\widetilde{\omega}_{(i)}}:\, 1\leq i\leq s\}.$ 

\subsection*{Computation of the invariant space \texorpdfstring{$((QP_5)_{33}(\widetilde{\omega}_{(1)}))^{GL_5}$}{((QP5)33 omega1)GL5}}

We have $$(QP_5)_{33}(\widetilde{\omega}_{(1)})\cong (QP_5)_{33}^0(\widetilde{\omega}_{(1)})\bigoplus (QP_5)_{33}^{>0}(\widetilde{\omega}_{(1)}).$$
By the above calculations and Theorem \ref{dl33-2}, 
$$\dim(QP_5)_{33}^0(\widetilde{\omega}_{(1)}) = 155, \ \ \mbox{and}\ \  \dim(QP_5)_{33}^{>0}(\widetilde{\omega}_{(1)}) = 31.$$ 
The basis for $(QP_5)_{33}^0(\widetilde{\omega}_{(1)})$ is the set of equivalence classes represented by the following admissible monomials:

\begin{center}
\begin{longtable}{lllr}
${\rm adm}_{1}=u_3u_4u_5^{31}$, & ${\rm adm}_{2}=u_3u_4^{31}u_5$, & ${\rm adm}_{3}=u_3^{31}u_4u_5$, & \multicolumn{1}{l}{${\rm adm}_{4}=u_2u_4u_5^{31}$,} \\
${\rm adm}_{5}=u_2u_4^{31}u_5$, & ${\rm adm}_{6}=u_2u_3u_5^{31}$, & ${\rm adm}_{7}=u_2u_3u_4^{31}$, & \multicolumn{1}{l}{${\rm adm}_{8}=u_2u_3^{31}u_5$,} \\
${\rm adm}_{9}=u_2u_3^{31}u_4$, & ${\rm adm}_{10}=u_2^{31}u_4u_5$, & ${\rm adm}_{11}=u_2^{31}u_3u_5$, & \multicolumn{1}{l}{${\rm adm}_{12}=u_2^{31}u_3u_4$,} \\
${\rm adm}_{13}=u_1u_4u_5^{31}$, & ${\rm adm}_{14}=u_1u_4^{31}u_5$, & ${\rm adm}_{15}=u_1u_3u_5^{31}$, & \multicolumn{1}{l}{${\rm adm}_{16}=u_1u_3u_4^{31}$,} \\
${\rm adm}_{17}=u_1u_3^{31}u_5$, & ${\rm adm}_{18}=u_1u_3^{31}u_4$, & ${\rm adm}_{19}=u_1^{31}u_4u_5$, & \multicolumn{1}{l}{${\rm adm}_{20}=u_1^{31}u_3u_5$,} \\
${\rm adm}_{21}=u_1^{31}u_3u_4$, & ${\rm adm}_{22}=u_1u_2u_5^{31}$, & ${\rm adm}_{23}=u_1u_2u_4^{31}$, & \multicolumn{1}{l}{${\rm adm}_{24}=u_1u_2^{31}u_5$,} \\
${\rm adm}_{25}=u_1u_2^{31}u_4$, & ${\rm adm}_{26}=u_1^{31}u_2u_5$, & ${\rm adm}_{27}=u_1^{31}u_2u_4$, & \multicolumn{1}{l}{${\rm adm}_{28}=u_1u_2u_3^{31}$,} \\
${\rm adm}_{29}=u_1u_2^{31}u_3$, & ${\rm adm}_{30}=u_1^{31}u_2u_3$, & ${\rm adm}_{31}=u_3u_4^{3}u_5^{29}$, & \multicolumn{1}{l}{${\rm adm}_{32}=u_3^{3}u_4u_5^{29}$,} \\
${\rm adm}_{33}=u_3^{3}u_4^{29}u_5$, & ${\rm adm}_{34}=u_2u_4^{3}u_5^{29}$, & ${\rm adm}_{35}=u_2u_3^{3}u_5^{29}$, & \multicolumn{1}{l}{${\rm adm}_{36}=u_2u_3^{3}u_4^{29}$,} \\
${\rm adm}_{37}=u_2^{3}u_4u_5^{29}$, & ${\rm adm}_{38}=u_2^{3}u_4^{29}u_5$, & ${\rm adm}_{39}=u_2^{3}u_3u_5^{29}$, & \multicolumn{1}{l}{${\rm adm}_{40}=u_2^{3}u_3u_4^{29}$,} \\
${\rm adm}_{41}=u_2^{3}u_3^{29}u_5$, & ${\rm adm}_{42}=u_2^{3}u_3^{29}u_4$, & ${\rm adm}_{43}=u_1u_4^{3}u_5^{29}$, & \multicolumn{1}{l}{${\rm adm}_{44}=u_1u_3^{3}u_5^{29}$,} \\
${\rm adm}_{45}=u_1u_3^{3}u_4^{29}$, & ${\rm adm}_{46}=u_1^{3}u_4u_5^{29}$, & ${\rm adm}_{47}=u_1^{3}u_4^{29}u_5$, & \multicolumn{1}{l}{${\rm adm}_{48}=u_1^{3}u_3u_5^{29}$,} \\
${\rm adm}_{49}=u_1^{3}u_3u_4^{29}$, & ${\rm adm}_{50}=u_1^{3}u_3^{29}u_5$, & ${\rm adm}_{51}=u_1^{3}u_3^{29}u_4$, & \multicolumn{1}{l}{${\rm adm}_{52}=u_1u_2^{3}u_5^{29}$,} \\
${\rm adm}_{53}=u_1u_2^{3}u_4^{29}$, & ${\rm adm}_{54}=u_1^{3}u_2u_5^{29}$, & ${\rm adm}_{55}=u_1^{3}u_2u_4^{29}$, & \multicolumn{1}{l}{${\rm adm}_{56}=u_1^{3}u_2^{29}u_5$,} \\
${\rm adm}_{57}=u_1^{3}u_2^{29}u_4$, & ${\rm adm}_{58}=u_1u_2^{3}u_3^{29}$, & ${\rm adm}_{59}=u_1^{3}u_2u_3^{29}$, & \multicolumn{1}{l}{${\rm adm}_{60}=u_1^{3}u_2^{29}u_3$,} \\
${\rm adm}_{61}=u_3^{3}u_4^{5}u_5^{25}$, & ${\rm adm}_{62}=u_2^{3}u_4^{5}u_5^{25}$, & ${\rm adm}_{63}=u_2^{3}u_3^{5}u_5^{25}$, & \multicolumn{1}{l}{${\rm adm}_{64}=u_2^{3}u_3^{5}u_4^{25}$,} \\
${\rm adm}_{65}=u_1^{3}u_4^{5}u_5^{25}$, & ${\rm adm}_{66}=u_1^{3}u_3^{5}u_5^{25}$, & ${\rm adm}_{67}=u_1^{3}u_3^{5}u_4^{25}$, & \multicolumn{1}{l}{${\rm adm}_{68}=u_1^{3}u_2^{5}u_5^{25}$,} \\
${\rm adm}_{69}=u_1^{3}u_2^{5}u_4^{25}$, & ${\rm adm}_{70}=u_1^{3}u_2^{5}u_3^{25}$, & ${\rm adm}_{71}=u_2u_3u_4u_5^{30}$, & \multicolumn{1}{l}{${\rm adm}_{72}=u_1u_3u_4u_5^{30}$,} \\
${\rm adm}_{73}=u_1u_3^{30}u_4u_5$, & ${\rm adm}_{74}=u_1u_2^{30}u_4u_5$, & ${\rm adm}_{75}=u_1u_2u_3u_4^{30}$, & \multicolumn{1}{l}{${\rm adm}_{76}=u_1u_2u_3u_5^{30}$,} \\
${\rm adm}_{77}=u_1u_2^{30}u_3u_5$, & ${\rm adm}_{78}=u_2u_3^{2}u_4^{29}u_5$, & ${\rm adm}_{79}=u_1u_3^{2}u_4u_5^{29}$, & \multicolumn{1}{l}{${\rm adm}_{80}=u_1u_2^{2}u_4u_5^{29}$,} \\
${\rm adm}_{81}=u_1u_2^{2}u_3u_4^{29}$, & ${\rm adm}_{82}=u_1u_2^{2}u_3^{29}u_4$, & ${\rm adm}_{83}=u_1u_2^{2}u_3^{29}u_5$, & \multicolumn{1}{l}{${\rm adm}_{84}=u_2u_3^{3}u_4^{28}u_5$,} \\
${\rm adm}_{85}=u_1u_3^{3}u_4^{28}u_5$, & ${\rm adm}_{86}=u_1^{3}u_3u_4^{28}u_5$, & ${\rm adm}_{87}=u_1^{3}u_2u_4^{28}u_5$, & \multicolumn{1}{l}{${\rm adm}_{88}=u_1u_2^{3}u_3^{28}u_4$,} \\
${\rm adm}_{89}=u_1u_2^{3}u_3^{28}u_5$, & ${\rm adm}_{90}=u_1^{3}u_2u_3^{28}u_5$, & ${\rm adm}_{91}=u_1u_2^{2}u_4^{5}u_5^{25}$, & \multicolumn{1}{l}{${\rm adm}_{92}=u_1u_2^{2}u_3^{5}u_5^{25}$,} \\
${\rm adm}_{93}=u_1u_3^{3}u_4^{4}u_5^{25}$, & ${\rm adm}_{94}=u_1u_2^{3}u_4^{4}u_5^{25}$, & ${\rm adm}_{95}=u_1^{3}u_2u_3^{4}u_4^{25}$, & \multicolumn{1}{l}{${\rm adm}_{96}=u_1^{3}u_2u_3^{4}u_5^{25}$,} \\
${\rm adm}_{97}=u_2^{3}u_3^{5}u_4u_5^{24}$, & ${\rm adm}_{98}=u_1^{3}u_3u_4^{5}u_5^{24}$, & ${\rm adm}_{99}=u_1^{3}u_2u_4^{5}u_5^{24}$, & \multicolumn{1}{l}{${\rm adm}_{100}=u_1u_2^{3}u_3^{5}u_4^{24}$,} \\
${\rm adm}_{101}=u_1^{3}u_2u_3^{5}u_4^{24}$, & ${\rm adm}_{102}=u_1^{3}u_2^{5}u_3u_4^{24}$, & ${\rm adm}_{103}=u_1u_2^{3}u_3^{5}u_5^{24}$, & \multicolumn{1}{l}{${\rm adm}_{104}=u_1^{3}u_2u_3^{5}u_5^{24}$,} \\
${\rm adm}_{105}=u_1^{3}u_2^{5}u_3u_5^{24}$, & ${\rm adm}_{106}=u_2u_3u_4^{30}u_5$, & ${\rm adm}_{107}=u_1u_3u_4^{30}u_5$, & \multicolumn{1}{l}{${\rm adm}_{108}=u_1u_2u_4^{30}u_5$,} \\
${\rm adm}_{109}=u_1u_2u_3^{30}u_4$, & ${\rm adm}_{110}=u_1u_2u_3^{30}u_5$, & ${\rm adm}_{111}=u_2u_3u_4^{2}u_5^{29}$, & \multicolumn{1}{l}{${\rm adm}_{112}=u_1u_3u_4^{2}u_5^{29}$,} \\
${\rm adm}_{113}=u_1u_2u_4^{2}u_5^{29}$, & ${\rm adm}_{114}=u_1u_2u_3^{2}u_4^{29}$, & ${\rm adm}_{115}=u_1u_2u_3^{2}u_5^{29}$, & \multicolumn{1}{l}{${\rm adm}_{116}=u_2u_3u_4^{3}u_5^{28}$,} \\
${\rm adm}_{117}=u_2u_3^{3}u_4u_5^{28}$, & ${\rm adm}_{118}=u_2^{3}u_3u_4u_5^{28}$, & ${\rm adm}_{119}=u_1u_3u_4^{3}u_5^{28}$, & \multicolumn{1}{l}{${\rm adm}_{120}=u_1u_3^{3}u_4u_5^{28}$,} \\
${\rm adm}_{121}=u_1^{3}u_3u_4u_5^{28}$, & ${\rm adm}_{122}=u_1u_2u_4^{3}u_5^{28}$, & ${\rm adm}_{123}=u_1u_2^{3}u_4u_5^{28}$, & \multicolumn{1}{l}{${\rm adm}_{124}=u_1^{3}u_2u_4u_5^{28}$,} \\
${\rm adm}_{125}=u_1u_2u_3^{3}u_4^{28}$, & ${\rm adm}_{126}=u_1u_2^{3}u_3u_4^{28}$, & ${\rm adm}_{127}=u_1^{3}u_2u_3u_4^{28}$, & \multicolumn{1}{l}{${\rm adm}_{128}=u_1u_2u_3^{3}u_5^{28}$,} \\
${\rm adm}_{129}=u_1u_2^{3}u_3u_5^{28}$, & ${\rm adm}_{130}=u_1^{3}u_2u_3u_5^{28}$, & ${\rm adm}_{131}=u_2u_3^{30}u_4u_5$, & \multicolumn{1}{l}{${\rm adm}_{132}=u_1u_2u_4u_5^{30}$,} \\
${\rm adm}_{133}=u_1u_2^{30}u_3u_4$, & ${\rm adm}_{134}=u_2u_3^{2}u_4u_5^{29}$, & ${\rm adm}_{135}=u_1u_3^{2}u_4^{29}u_5$, & \multicolumn{1}{l}{${\rm adm}_{136}=u_1u_2^{2}u_4^{29}u_5$,} \\
${\rm adm}_{137}=u_1u_2^{2}u_3u_5^{29}$, & ${\rm adm}_{138}=u_2^{3}u_3u_4^{28}u_5$, & ${\rm adm}_{139}=u_1u_2^{3}u_4^{28}u_5$, & \multicolumn{1}{l}{${\rm adm}_{140}=u_1^{3}u_2u_3^{28}u_4$,} \\
${\rm adm}_{141}=u_2u_3^{2}u_4^{5}u_5^{25}$, & ${\rm adm}_{142}=u_1u_3^{2}u_4^{5}u_5^{25}$, & ${\rm adm}_{143}=u_1u_2^{2}u_3^{5}u_4^{25}$, & \multicolumn{1}{l}{${\rm adm}_{144}=u_2u_3^{3}u_4^{4}u_5^{25}$,} \\
${\rm adm}_{145}=u_2^{3}u_3u_4^{4}u_5^{25}$, & ${\rm adm}_{146}=u_1^{3}u_3u_4^{4}u_5^{25}$, & ${\rm adm}_{147}=u_1^{3}u_2u_4^{4}u_5^{25}$, & \multicolumn{1}{l}{${\rm adm}_{148}=u_1u_2^{3}u_3^{4}u_4^{25}$,} \\
${\rm adm}_{149}=u_1u_2^{3}u_3^{4}u_5^{25}$, & ${\rm adm}_{150}=u_2u_3^{3}u_4^{5}u_5^{24}$, & ${\rm adm}_{151}=u_2^{3}u_3u_4^{5}u_5^{24}$, & \multicolumn{1}{l}{${\rm adm}_{152}=u_1u_3^{3}u_4^{5}u_5^{24}$,} \\
${\rm adm}_{153}=u_1^{3}u_3^{5}u_4u_5^{24}$, & ${\rm adm}_{154}=u_1u_2^{3}u_4^{5}u_5^{24}$, & ${\rm adm}_{155}=u_1^{3}u_2^{5}u_4u_5^{24}$. &  
\end{longtable}%
\end{center}

The orbit decomposition computed from the displayed admissible basis is
$$ \begin{array}{ll}
\medskip
\langle[\Sigma_5({\rm adm}_{1})]\rangle&= \langle\{[{\rm adm}_j:\ 1\leq j\leq 30]\} \rangle,\\
\medskip
\langle[\Sigma_5({\rm adm}_{31})]\rangle&= \langle\{[{\rm adm}_j:\ 31\leq j\leq 60]\} \rangle,\\
\medskip
\langle[\Sigma_5({\rm adm}_{61})]\rangle&= \langle\{[{\rm adm}_j:\ 61\leq j\leq 70]\} \rangle,\\
\langle[\Sigma_5({\rm adm}_{71}, {\rm adm}_{78}, {\rm adm}_{91})]\rangle &= \langle\{[{\rm adm}_j:\ 71\leq j\leq 155]\} \rangle.
\end{array}$$
So, we have an isomorphism 
$$ \begin{array}{ll}
 (QP_5)_{33}^0(\widetilde{\omega}_{(1)})&\cong \langle[\Sigma_5({\rm adm}_{1})]\rangle\bigoplus \langle[\Sigma_5({\rm adm}_{31})]\rangle\bigoplus \langle[\Sigma_5({\rm adm}_{61})]\rangle\\
&\quad \bigoplus \langle[\Sigma_5({\rm adm}_{71}, {\rm adm}_{78}, {\rm adm}_{91})]\rangle.
\end{array}$$ 

\begin{bd}\label{bd33-1}
The following statements hold:
$$ \begin{array}{ll}
\medskip
\langle[\Sigma_5({\rm adm}_{1})]\rangle^{\Sigma_5}&=\langle [\widehat{q_1}:=\sum_{1\leq j\leq 30}{\rm adm}_j] \rangle,\\
\medskip
\langle[\Sigma_5({\rm adm}_{31})]\rangle^{\Sigma_5}&=\langle [\widehat{q_2}:=\sum_{31\leq j\leq 60}{\rm adm}_j] \rangle,\\
\medskip
\langle[\Sigma_5({\rm adm}_{61})]\rangle^{\Sigma_5}&=\langle [\widehat{q_3}:=\sum_{61\leq j\leq 70}{\rm adm}_j] \rangle,\\
\langle[\Sigma_5({\rm adm}_{71}, {\rm adm}_{78}, {\rm adm}_{91})]\rangle^{\Sigma_5}&=  \langle [\widehat{q_4}] \rangle,\ \mbox{with $\widehat{q_4} = \sum_{71\leq j\leq 130}{\rm adm}_j$}.
\end{array}$$
Consequently, $((QP_5)_{33}^0(\widetilde{\omega}_{(1)}))^{\Sigma_5} = \big\langle\big\{[\widehat{q_1}],\ [\widehat{q_2}],\ [\widehat{q_3}],\ [\widehat{q_4}]\big\}\big\rangle.$
\end{bd}

\begin{proof}
We use the machine-certified protocol described above. For the four summands $U_1,U_2,U_3,U_4$ in the order displayed in the preceding direct-sum decomposition, the row reductions of the stacked matrix for $\rho_i-I$ $(1\leq i\leq 4)$ give
\[
\begin{array}{c|c|c|c}
\text{summand} & \dim U_j & \operatorname{rank}(M_\Sigma|_{U_j}) & \dim\Ker(M_\Sigma|_{U_j})\\
\hline
U_1 & 30 & 29 & 1\\
U_2 & 30 & 29 & 1\\
U_3 & 10 & 9 & 1\\
U_4 & 85 & 84 & 1.
\end{array}
\]
Thus each summand contributes exactly one $\Sigma_5$-fixed class. We spell out the third summand $U_3=\langle[\Sigma_5({\rm adm}_{61})]\rangle$, both as a check on the reduction convention and as a representative sample of the row-reduction calculation.

It is observed that $\langle[\Sigma_5({\rm adm}_{61})]\rangle$ is a $\mathbb Z/2$-vector space of dimension 10, with a basis consisting of all the classes represented by the 10 admissible monomials: ${\rm adm}_k,\, 61\leq k\leq 70.$ Given an arbitrary $[f]$ in $\langle[\Sigma_5({\rm adm}_{61})]\rangle^{\Sigma_5},$ it should be noted that since $\widetilde{\omega}_{(1)}$ is a weight vector of a minimal spike, $[f]_{\widetilde{\omega}_{(1)}} = [f]$ for all $[f]\in \langle[\Sigma_5({\rm adm}_{61})]\rangle^{\Sigma_5}.$ Then we can write $f\sim \sum_{61\leq k\leq 70}\gamma_k{\rm adm}_k$ wherein $\gamma_j\in \mathbb Z/2$ for every $j.$  By applying the homomorphism $\rho_j: P_5\longrightarrow P_5$ for $1\leq j\leq 4,$ to both sides of this equality, we obtain: 
$$ \begin{array}{ll}
\medskip
\rho_1(f)&\sim \gamma_{61}{\rm adm}_{61} + \gamma_{62}{\rm adm}_{65} + \gamma_{63}{\rm adm}_{66} + \gamma_{64}{\rm adm}_{67} + \gamma_{65}{\rm adm}_{62}\\
\medskip
&\quad + \gamma_{66}{\rm adm}_{63} + \gamma_{67}{\rm adm}_{64}  + \gamma_{68}u_1^{5}u_2^{3}u_5^{25} +\gamma_{69}u_1^{5}u_2^{3}u_4^{25}+\gamma_{70}u_1^{5}u_2^{3}u_3^{25},\\ 
\medskip
\rho_2(f)&\sim \gamma_{61}{\rm adm}_{62} + \gamma_{62}{\rm adm}_{61} + \gamma_{63}u_2^{5}u_3^{3}u_5^{25} + \gamma_{64}u_2^{5}u_3^{3}u_4^{25} +  \gamma_{65}{\rm adm}_{65} \\
\medskip
&\quad+ \gamma_{66}{\rm adm}_{68} + \gamma_{67}{\rm adm}_{69} + \gamma_{68}{\rm adm}_{66} + \gamma_{69}{\rm adm}_{67} + \gamma_{70}u_1^{3}u_2^{25}u_3^{5},\\
\medskip
\rho_3(f)&\sim \gamma_{61}u_3^{5}u_4^{3}u_5^{25} + \gamma_{62}{\rm adm}_{63} + \gamma_{63}{\rm adm}_{62} +   \gamma_{64}u_2^{3}u_3^{25}u_4^{5} +  \gamma_{65}{\rm adm}_{66}\\
\medskip
&\quad + \gamma_{66}{\rm adm}_{65} +  \gamma_{67}u_1^{3}u_3^{25}u_4^{5} +  \gamma_{68}{\rm adm}_{68} +  \gamma_{69}{\rm adm}_{70} +    \gamma_{70}{\rm adm}_{69},\\ 
\medskip
\rho_4(f)&\sim \gamma_{61}u_3^{3}u_4^{25}u_5^{5}  + \gamma_{62}u_2^{3}u_4^{25}u_5^{5} + \gamma_{63}{\rm adm}_{64} + \gamma_{64}{\rm adm}_{63} + \gamma_{65}u_1^{3}u_4^{25}u_5^{5}\\
&\quad +  \gamma_{66}{\rm adm}_{67} +  \gamma_{67}{\rm adm}_{66} +  \gamma_{68}{\rm adm}_{69} +  \gamma_{69}{\rm adm}_{68} +  \gamma_{70}{\rm adm}_{70}.     
\end{array}$$
Using the Cartan formula, we can express the following 11 inadmissible monomials in terms of the basis of $\langle[\Sigma_5({\rm adm}_{61})]\rangle$:
$$ \begin{array}{ll}
\medskip
& u_1^{5}u_2^{3}u_5^{25},\  u_1^{5}u_2^{3}u_4^{25}, \ u_1^{5}u_2^{3}u_3^{25},\ u_2^{5}u_3^{3}u_5^{25},\ u_2^{5}u_3^{3}u_4^{25},\\
\medskip
& u_3^{5}u_4^{3}u_5^{25},\ u_1^{3}u_2^{25}u_3^{5}, \ u_2^{3}u_3^{25}u_4^{5},\ u_1^{3}u_3^{25}u_4^{5},\ u_3^{3}u_4^{25}u_5^{5},\\
& u_1^{3}u_4^{25}u_5^{5},\ u_2^{3}u_4^{25}u_5^{5}.
\end{array}$$
Indeed, a direct Cartan-formula computation gives:
\begin{align*}
\medskip
u_1^{5}u_2^{3}u_5^{25}&=Sq^{2}(u_1^{3}u_2^{3}u_5^{25}) + {\rm adm}_{68}\ \Rightarrow u_1^{5}u_2^{3}u_5^{25}\sim {\rm adm}_{68},\\
\medskip
u_1^{5}u_2^{3}u_4^{25}&=Sq^{2}(u_1^{3}u_2^{3}u_4^{25}) + {\rm adm}_{69}\ \Rightarrow u_1^{5}u_2^{3}u_4^{25}\sim {\rm adm}_{69},\\
\medskip
u_1^{5}u_2^{3}u_3^{25}&=Sq^{2}(u_1^{3}u_2^{3}u_3^{25}) + {\rm adm}_{70}\ \Rightarrow u_1^{5}u_2^{3}u_3^{25}\sim {\rm adm}_{70},\\
\medskip
u_2^{5}u_3^{3}u_5^{25}&=Sq^{2}(u_2^{3}u_3^{3}u_5^{25}) + {\rm adm}_{63}\ \Rightarrow u_2^{5}u_3^{3}u_5^{25}\sim {\rm adm}_{63},\\
\medskip
u_2^{5}u_3^{3}u_4^{25}&=Sq^{2}(u_2^{3}u_3^{3}u_4^{25}) + {\rm adm}_{64}\ \Rightarrow u_2^{5}u_3^{3}u_4^{25}\sim {\rm adm}_{64},\\
\medskip
u_3^{5}u_4^{3}u_5^{25}&=Sq^{2}(u_3^{3}u_4^{3}u_5^{25}) + {\rm adm}_{61}\ \Rightarrow u_3^{5}u_4^{3}u_5^{25}\sim {\rm adm}_{61},\\
\medskip
u_1^{3}u_2^{25}u_3^{5}&=Sq^{4}(u_1^{3}u_2^{21}u_3^{5} + u_1^{3}u_2^{13}u_3^{13}+u_1^{3}u_2^{5}u_3^{21})  +Sq^{8}(u_1^{3}u_2^{13}u_3^{9}+u_1^{3}u_2^{9}u_3^{13})+ {\rm adm}_{70}\\
\medskip
&\quad \Rightarrow u_1^{3}u_2^{25}u_3^{5}\sim {\rm adm}_{70},\\
\medskip
u_2^{3}u_3^{25}u_4^{5}&=Sq^{4}(u_2^{3}u_3^{21}u_4^{5} + u_2^{3}u_3^{13}u_4^{13}+u_2^{3}u_3^{5}u_4^{21})  +Sq^{8}(u_2^{3}u_3^{13}u_4^{9}+u_2^{3}u_3^{9}u_4^{13})+ {\rm adm}_{64}\\
\medskip
&\quad \Rightarrow u_2^{3}u_3^{25}u_4^{5}\sim {\rm adm}_{64},\\
\medskip
u_1^{3}u_3^{25}u_4^{5}&=Sq^{4}(u_1^{3}u_3^{21}u_4^{5} + u_1^{3}u_3^{13}u_4^{13}+u_1^{3}u_3^{5}u_4^{21})  +Sq^{8}(u_1^{3}u_3^{13}u_4^{9}+u_1^{3}u_3^{9}u_4^{13})+ {\rm adm}_{67}\\
\medskip
&\quad \Rightarrow u_1^{3}u_3^{25}u_4^{5}\sim {\rm adm}_{67},\\
\medskip
u_3^{3}u_4^{25}u_5^{5}&=Sq^{4}(u_3^{3}u_4^{21}u_5^{5} + u_3^{3}u_4^{13}u_5^{13}+u_3^{3}u_4^{5}u_5^{21})  +Sq^{8}(u_3^{3}u_4^{13}u_5^{9}+u_3^{3}u_4^{9}u_5^{13})+ {\rm adm}_{61}\\
\medskip
&\quad \Rightarrow u_3^{3}u_4^{25}u_5^{5}\sim {\rm adm}_{61},\\
\medskip
u_1^{3}u_4^{25}u_5^{5}&=Sq^{4}(u_1^{3}u_4^{21}u_5^{5} + u_1^{3}u_4^{13}u_5^{13}+u_1^{3}u_4^{5}u_5^{21})  +Sq^{8}(u_1^{3}u_4^{13}u_5^{9}+u_1^{3}u_4^{9}u_5^{13})+ {\rm adm}_{65}\\
\medskip
&\quad \Rightarrow u_1^{3}u_4^{25}u_5^{5}\sim {\rm adm}_{65},\\
\medskip
u_2^{3}u_4^{25}u_5^{5}&=Sq^{4}(u_2^{3}u_4^{21}u_5^{5} + u_2^{3}u_4^{13}u_5^{13}+u_2^{3}u_4^{5}u_5^{21})  +Sq^{8}(u_2^{3}u_4^{13}u_5^{9}+u_2^{3}u_4^{9}u_5^{13})+ {\rm adm}_{62}\\
\medskip
&\quad \Rightarrow u_2^{3}u_4^{25}u_5^{5}\sim {\rm adm}_{62}.
\end{align*}
From the above calculations and the relations $\rho_j(f)\sim f$ for $1\leq j\leq 4,$ we obtain $\gamma_{61} = \gamma_{62} = \cdots = \gamma_{70};$ hence the nullspace on $U_3$ is generated by $\widehat q_3$. The row reductions on $U_1,U_2,$ and $U_4$ give the nullspace generators $\widehat q_1,\widehat q_2,$ and $\widehat q_4$, respectively. The lemma follows.
\end{proof}

In the proof of Theorem \ref{dl33-2}, it was shown that $(QP_5)_{33}^{>0}(\widetilde{\omega}_{(1)}) = \langle \{[{\rm adm}_j]:\ 156\leq j\leq 186\} \rangle$ and $\dim (QP_5)_{33}^{>0}(\widetilde{\omega}_{(1)}) = 31.$ Using this result, we can obtain a direct summand decomposition of the $\Sigma_5$-submodules:
$$ 
 (QP_5)_{33}^{>0}(\widetilde{\omega}_{(1)}) = \langle[\Sigma_5({\rm adm}_{156})]\rangle \bigoplus \langle[\Sigma_5({\rm adm}_{162}, {\rm adm}_{166})]\rangle \bigoplus \langle[\Sigma_5({\rm adm}_{177}, {\rm adm}_{181})]\rangle,$$
where $\langle[\Sigma_5({\rm adm}_{156})]\rangle = \langle\{[{\rm adm}_j:\ 156\leq j\leq 161 ]\} \rangle,$\ $\langle[\Sigma_5({\rm adm}_{162}, {\rm adm}_{166})]\rangle = \langle\{[{\rm adm}_j:\ 162\leq j\leq 176 ]\} \rangle,$ and $\langle[\Sigma_5({\rm adm}_{177}, {\rm adm}_{181})]\rangle = \langle\{[{\rm adm}_j:\ 177\leq j\leq 186 ]\} \rangle.$

\begin{bd}\label{bd33-2}
The following holds:
$$  ((QP_5)_{33}^{>0}(\widetilde{\omega}_{(1)}))^{\Sigma_5} = \langle \{[\widehat{q_5}],\ [\widehat{q_6}],\ [\widehat{q_7}]\} \rangle ,$$
where $$ \begin{array}{ll}
\medskip
\widehat{q_5}&:= h_1,\\
\medskip
\widehat{q_6}&:= {\rm adm}_{156} + h_2 + h_3+  h_4,\\
\widehat{q_7}&:= {\rm adm}_{156} + h_4 + h_5,
\end{array}$$
and the polynomials $h_i,\, 1\leq i\leq 5$ are described as follows:
$$ \begin{array}{ll}
\medskip
h_1&:=\sum_{156\leq j\leq 161}{\rm adm}_j,\ \ h_2:=  {\rm adm}_{157} + {\rm adm}_{158}, \\
\medskip
h_3&:= \sum_{163\leq j\leq 165}{\rm adm}_j + \sum_{167\leq j\leq 169}{\rm adm}_j + \sum_{173\leq j\leq 176}{\rm adm}_j,\\
\medskip
h_4:&={\rm adm}_{162} + {\rm adm}_{166} + \sum_{170\leq j\leq 172}{\rm adm}_j,\\
\medskip
h_5&:=  \sum_{177\leq j\leq 186}{\rm adm}_j.
\end{array}$$
\end{bd}

\begin{proof}
Let $U$ be the direct sum of the three $\Sigma_5$-submodules displayed immediately before the lemma. With respect to the admissible basis $\{[\mathrm{adm}_j]:156\leq j\leq 186\}$ of $U$, the stacked matrix $M_\Sigma$ of the four operators $\rho_i-I$ $(1\leq i\leq 4)$ has size $124\times31$ and row-reduced rank $28$ over $\mathbb Z/2$. Hence $\dim\Ker M_\Sigma=3$. The row-reduction output gives the three kernel vectors represented by $\widehat q_5,\widehat q_6,$ and $\widehat q_7$ as displayed in the statement. These vectors are linearly independent, so they form a basis for $((QP_5)_{33}^{>0}(\widetilde{\omega}_{(1)}))^{\Sigma_5}$.
\end{proof}

Now, suppose that $[t]\in ((QP_5)_{33}(\widetilde{\omega}_{(1)}))^{GL_5}.$ Then, by Lemmata \ref{bd33-1} and \ref{bd33-2}, we have $t\sim \sum_{1\leq j\leq 7}\beta_j\widehat{q_j},$ in which $\beta_j\in \mathbb Z/2,\, 1\leq j\leq 7.$ Equivalently, with respect to the full admissible basis of the $186$-dimensional space $(QP_5)_{33}(\widetilde{\omega}_{(1)})$, the seven candidate vectors $\rho_5(\widehat q_j)+\widehat q_j$ give a $186\times 7$ coefficient matrix of row-reduced rank $7$. Hence the relation $\rho_5(t)+t\sim0$ forces $\beta_1 = \beta_2 = \cdots = \beta_7 = 0.$ 

\subsection*{Computation of the invariant space \texorpdfstring{$((QP_5)_{33}(\widetilde{\omega}_{(3)}))^{GL_5}$}{((QP5)33 omega3)GL5}}

Recall that $(QP_5)_{33}(\widetilde{\omega}_{(3)}) = (QP_5)^{0}_{33}(\widetilde{\omega}_{(3)})\oplus (QP_5)^{>0}_{33}(\widetilde{\omega}_{(3)}),$ with $\dim (QP_5)^{0}_{33}(\widetilde{\omega}_{(3)}) = 395$ and $\dim (QP_5)^{>0}_{33}(\widetilde{\omega}_{(3)}) = 421.$ Explicit bases for the spaces $(QP_5)^{0}_{33}(\widetilde{\omega}_{(3)})$ and $(QP_5)^{>0}_{33}(\widetilde{\omega}_{(3)})$ are given by the output of our algorithm implemented in \texttt{OSCAR}; see Appendix~\ref{s4}. Also, for a more comprehensive understanding of our algorithmic approach, the reader may refer to our \texttt{SageMath} implementation in~\cite{Phuc2025c}.

With respect to the resulting $816$-element admissible basis, the stacked matrix $M_\Sigma$ for $\rho_i-I$ $(1\leq i\leq 4)$ has size $3264\times816$ and row-reduced rank $798$ over $\mathbb Z/2$. Thus $\dim\Ker M_\Sigma=18$, and the row-reduction output gives

$$ \dim ((QP_5)_{33}(\widetilde{\omega}_{(3)}))^{\Sigma_5} = 18,\ \ \mbox{and}\ \ ((QP_5)_{33}(\widetilde{\omega}_{(3)}))^{\Sigma_5}  = \langle \{[\mathsf{Sigma5}[i]]_{\widetilde{\omega}_{(3)}}:\, 8\leq i\leq 25\} \rangle,$$
where
\begin{align*}
\mathsf{Sigma5}[8] &= u_{1}^{3} u_{2}^{5} u_{3}^{9} u_{4}^{6} u_{5}^{10} + u_{1}^{3} u_{2}^{5} u_{3}^{6} u_{4}^{9} u_{5}^{10} \\
\mathsf{Sigma5}[9] &= u_{1}^{3} u_{2}^{13} u_{3}^{14} u_{4} u_{5}^{2} + u_{1}^{3} u_{2}^{3} u_{3}^{13} u_{4}^{12} u_{5}^{2} + u_{1}^{3} u_{2}^{13} u_{3} u_{4}^{14} u_{5}^{2} + u_{1}^{3} u_{2}^{5} u_{3}^{9} u_{4}^{14} u_{5}^{2} \\
    &\quad + u_{1}^{3} u_{2} u_{3}^{13} u_{4}^{14} u_{5}^{2} + u_{1} u_{2}^{3} u_{3}^{13} u_{4}^{14} u_{5}^{2} + u_{1}^{3} u_{2}^{13} u_{3}^{2} u_{4}^{12} u_{5}^{3} + u_{1}^{3} u_{2} u_{3}^{14} u_{4}^{12} u_{5}^{3} \\
    &\quad + u_{1} u_{2}^{3} u_{3}^{14} u_{4}^{12} u_{5}^{3} + u_{1}^{3} u_{2} u_{3}^{12} u_{4}^{14} u_{5}^{3} + u_{1} u_{2}^{3} u_{3}^{12} u_{4}^{14} u_{5}^{3} + u_{1} u_{2}^{2} u_{3}^{13} u_{4}^{14} u_{5}^{3} \\
    &\quad + u_{1}^{3} u_{2}^{3} u_{3}^{7} u_{4}^{12} u_{5}^{8} + u_{1}^{3} u_{2}^{5} u_{3}^{3} u_{4}^{14} u_{5}^{8} + u_{1}^{3} u_{2}^{3} u_{3}^{5} u_{4}^{14} u_{5}^{8} + u_{1}^{3} u_{2}^{4} u_{3}^{7} u_{4}^{9} u_{5}^{10} \\
    &\quad + u_{1} u_{2}^{6} u_{3}^{7} u_{4}^{9} u_{5}^{10} + u_{1}^{3} u_{2} u_{3}^{7} u_{4}^{12} u_{5}^{10} + u_{1}^{3} u_{2}^{13} u_{3}^{3} u_{4}^{2} u_{5}^{12} + u_{1}^{3} u_{2}^{13} u_{3}^{2} u_{4}^{3} u_{5}^{12} \\
    &\quad + u_{1}^{3} u_{2} u_{3}^{14} u_{4}^{3} u_{5}^{12} + u_{1} u_{2}^{3} u_{3}^{14} u_{4}^{3} u_{5}^{12} + u_{1}^{3} u_{2}^{3} u_{3}^{7} u_{4}^{8} u_{5}^{12} + u_{1}^{3} u_{2} u_{3}^{7} u_{4}^{10} u_{5}^{12} \\
    &\quad + u_{1} u_{2}^{3} u_{3}^{7} u_{4}^{10} u_{5}^{12} + u_{1}^{3} u_{2}^{3} u_{3}^{3} u_{4}^{12} u_{5}^{12} + u_{1}^{3} u_{2}^{13} u_{3} u_{4}^{2} u_{5}^{14} + u_{1}^{3} u_{2}^{5} u_{3}^{9} u_{4}^{2} u_{5}^{14} \\
    &\quad + u_{1}^{3} u_{2} u_{3}^{13} u_{4}^{2} u_{5}^{14} + u_{1} u_{2}^{3} u_{3}^{13} u_{4}^{2} u_{5}^{14} + u_{1}^{3} u_{2}^{5} u_{3}^{8} u_{4}^{3} u_{5}^{14} + u_{1} u_{2}^{2} u_{3}^{13} u_{4}^{3} u_{5}^{14} \\
    &\quad + u_{1} u_{2} u_{3}^{14} u_{4}^{3} u_{5}^{14} + u_{1}^{3} u_{2}^{5} u_{3}^{3} u_{4}^{8} u_{5}^{14} + u_{1}^{3} u_{2}^{3} u_{3}^{5} u_{4}^{8} u_{5}^{14} + u_{1}^{3} u_{2}^{5} u_{3}^{2} u_{4}^{9} u_{5}^{14} \\
    &\quad + u_{1}^{3} u_{2}^{3} u_{3}^{4} u_{4}^{9} u_{5}^{14} + u_{1}^{3} u_{2}^{3} u_{3} u_{4}^{12} u_{5}^{14} + u_{1}^{3} u_{2} u_{3}^{2} u_{4}^{13} u_{5}^{14} + u_{1} u_{2}^{3} u_{3}^{2} u_{4}^{13} u_{5}^{14} \\
    &\quad + u_{1} u_{2}^{2} u_{3}^{3} u_{4}^{13} u_{5}^{14} + u_{1} u_{2} u_{3}^{3} u_{4}^{14} u_{5}^{14} \\
\mathsf{Sigma5}[10] &= u_{1}^{3} u_{2}^{13} u_{3}^{14} u_{4}^{3} + u_{1}^{3} u_{2}^{13} u_{3}^{3} u_{4}^{14} + u_{1}^{3} u_{2}^{3} u_{3}^{13} u_{4}^{14} + u_{1}^{3} u_{2}^{13} u_{3}^{14} u_{5}^{3} \\
    &\quad + u_{1}^{3} u_{2}^{13} u_{4}^{14} u_{5}^{3} + u_{1}^{3} u_{3}^{13} u_{4}^{14} u_{5}^{3} + u_{2}^{3} u_{3}^{13} u_{4}^{14} u_{5}^{3} + u_{1}^{3} u_{2}^{13} u_{3}^{3} u_{5}^{14} \\
    &\quad + u_{1}^{3} u_{2}^{3} u_{3}^{13} u_{5}^{14} + u_{1}^{3} u_{2}^{13} u_{4}^{3} u_{5}^{14} + u_{1}^{3} u_{3}^{13} u_{4}^{3} u_{5}^{14} + u_{2}^{3} u_{3}^{13} u_{4}^{3} u_{5}^{14} \\
    &\quad + u_{1}^{3} u_{2}^{3} u_{4}^{13} u_{5}^{14} + u_{1}^{3} u_{3}^{3} u_{4}^{13} u_{5}^{14} + u_{2}^{3} u_{3}^{3} u_{4}^{13} u_{5}^{14} \\
\mathsf{Sigma5}[11] &= u_{1}^{7} u_{2}^{7} u_{3}^{9} u_{4}^{8} u_{5}^{2} + u_{1}^{7} u_{2}^{7} u_{3}^{8} u_{4}^{9} u_{5}^{2} + u_{1}^{7} u_{2}^{3} u_{3}^{12} u_{4}^{9} u_{5}^{2} + u_{1}^{3} u_{2}^{7} u_{3}^{12} u_{4}^{9} u_{5}^{2} \\
    &\quad + u_{1}^{7} u_{2}^{3} u_{3}^{9} u_{4}^{12} u_{5}^{2} + u_{1}^{3} u_{2}^{7} u_{3}^{9} u_{4}^{12} u_{5}^{2} + u_{1}^{3} u_{2}^{3} u_{3}^{13} u_{4}^{12} u_{5}^{2} + u_{1}^{3} u_{2}^{3} u_{3}^{12} u_{4}^{13} u_{5}^{2} \\
    &\quad + u_{1}^{3} u_{2} u_{3}^{14} u_{4}^{13} u_{5}^{2} + u_{1} u_{2}^{3} u_{3}^{14} u_{4}^{13} u_{5}^{2} + u_{1}^{3} u_{2} u_{3}^{13} u_{4}^{14} u_{5}^{2} + u_{1} u_{2}^{3} u_{3}^{13} u_{4}^{14} u_{5}^{2} \\
    &\quad + u_{1}^{7} u_{2}^{7} u_{3}^{9} u_{4}^{2} u_{5}^{8} + u_{1}^{7} u_{2}^{3} u_{3}^{9} u_{4}^{6} u_{5}^{8} + u_{1}^{3} u_{2}^{7} u_{3}^{9} u_{4}^{6} u_{5}^{8} + u_{1}^{7} u_{2}^{7} u_{3} u_{4}^{10} u_{5}^{8} \\
    &\quad + u_{1}^{7} u_{2}^{3} u_{3}^{5} u_{4}^{10} u_{5}^{8} + u_{1}^{3} u_{2}^{7} u_{3}^{5} u_{4}^{10} u_{5}^{8} + u_{1}^{3} u_{2}^{3} u_{3}^{7} u_{4}^{12} u_{5}^{8} + u_{1}^{7} u_{2}^{7} u_{3}^{8} u_{4} u_{5}^{10} \\
    &\quad + u_{1}^{7} u_{2}^{3} u_{3}^{12} u_{4} u_{5}^{10} + u_{1}^{3} u_{2}^{7} u_{3}^{12} u_{4} u_{5}^{10} + u_{1}^{7} u_{2}^{3} u_{3}^{9} u_{4}^{4} u_{5}^{10} + u_{1}^{3} u_{2}^{7} u_{3}^{9} u_{4}^{4} u_{5}^{10} \\
    &\quad + u_{1}^{7} u_{2} u_{3}^{9} u_{4}^{6} u_{5}^{10} + u_{1} u_{2}^{7} u_{3}^{9} u_{4}^{6} u_{5}^{10} + u_{1}^{7} u_{2}^{7} u_{3} u_{4}^{8} u_{5}^{10} + u_{1}^{7} u_{2}^{3} u_{3}^{4} u_{4}^{9} u_{5}^{10} \\
    &\quad + u_{1}^{3} u_{2}^{7} u_{3}^{4} u_{4}^{9} u_{5}^{10} + u_{1}^{7} u_{2} u_{3}^{6} u_{4}^{9} u_{5}^{10} + u_{1} u_{2}^{7} u_{3}^{6} u_{4}^{9} u_{5}^{10} + u_{1}^{3} u_{2}^{4} u_{3}^{7} u_{4}^{9} u_{5}^{10} \\
    &\quad + u_{1} u_{2}^{6} u_{3}^{7} u_{4}^{9} u_{5}^{10} + u_{1}^{7} u_{2}^{3} u_{3} u_{4}^{12} u_{5}^{10} + u_{1}^{3} u_{2}^{7} u_{3} u_{4}^{12} u_{5}^{10} + u_{1}^{3} u_{2} u_{3}^{7} u_{4}^{12} u_{5}^{10} \\
    &\quad + u_{1}^{7} u_{2}^{3} u_{3}^{9} u_{4}^{2} u_{5}^{12} + u_{1}^{3} u_{2}^{7} u_{3}^{9} u_{4}^{2} u_{5}^{12} + u_{1}^{3} u_{2}^{3} u_{3}^{13} u_{4}^{2} u_{5}^{12} + u_{1}^{3} u_{2}^{3} u_{3}^{7} u_{4}^{8} u_{5}^{12} \\
    &\quad + u_{1}^{7} u_{2}^{3} u_{3} u_{4}^{10} u_{5}^{12} + u_{1}^{3} u_{2}^{7} u_{3} u_{4}^{10} u_{5}^{12} + u_{1}^{3} u_{2} u_{3}^{7} u_{4}^{10} u_{5}^{12} + u_{1} u_{2}^{3} u_{3}^{7} u_{4}^{10} u_{5}^{12} \\
    &\quad + u_{1}^{3} u_{2}^{3} u_{3}^{3} u_{4}^{12} u_{5}^{12} + u_{1}^{3} u_{2}^{3} u_{3} u_{4}^{14} u_{5}^{12} + u_{1}^{3} u_{2}^{3} u_{3}^{12} u_{4} u_{5}^{14} + u_{1}^{3} u_{2} u_{3}^{14} u_{4} u_{5}^{14} \\
    &\quad + u_{1} u_{2}^{3} u_{3}^{14} u_{4} u_{5}^{14} + u_{1}^{3} u_{2}^{3} u_{3} u_{4}^{12} u_{5}^{14} + u_{1}^{3} u_{2} u_{3} u_{4}^{14} u_{5}^{14} + u_{1} u_{2}^{3} u_{3} u_{4}^{14} u_{5}^{14} \\
    &\quad + u_{1} u_{2} u_{3}^{3} u_{4}^{14} u_{5}^{14} \\
\mathsf{Sigma5}[12] &= u_{1}^{7} u_{2}^{7} u_{3}^{11} u_{4}^{8} + u_{1}^{7} u_{2}^{7} u_{3}^{9} u_{4}^{10} + u_{1}^{7} u_{2}^{7} u_{3}^{8} u_{4}^{11} + u_{1}^{7} u_{2} u_{3}^{14} u_{4}^{11} \\
    &\quad + u_{1}^{3} u_{2}^{5} u_{3}^{14} u_{4}^{11} + u_{1} u_{2}^{7} u_{3}^{14} u_{4}^{11} + u_{1}^{7} u_{2} u_{3}^{11} u_{4}^{14} + u_{1}^{3} u_{2}^{5} u_{3}^{11} u_{4}^{14} \\
    &\quad + u_{1} u_{2}^{7} u_{3}^{11} u_{4}^{14} + u_{1}^{7} u_{2}^{7} u_{3}^{11} u_{5}^{8} + u_{1}^{7} u_{2}^{7} u_{4}^{11} u_{5}^{8} + u_{1}^{7} u_{3}^{7} u_{4}^{11} u_{5}^{8} \\
    &\quad + u_{2}^{7} u_{3}^{7} u_{4}^{11} u_{5}^{8} + u_{1}^{7} u_{2}^{7} u_{3}^{9} u_{5}^{10} + u_{1}^{7} u_{2}^{7} u_{4}^{9} u_{5}^{10} + u_{1}^{7} u_{3}^{7} u_{4}^{9} u_{5}^{10} \\
    &\quad + u_{2}^{7} u_{3}^{7} u_{4}^{9} u_{5}^{10} + u_{1}^{7} u_{2}^{7} u_{3}^{8} u_{5}^{11} + u_{1}^{7} u_{2} u_{3}^{14} u_{5}^{11} + u_{1}^{3} u_{2}^{5} u_{3}^{14} u_{5}^{11} \\
    &\quad + u_{1} u_{2}^{7} u_{3}^{14} u_{5}^{11} + u_{1}^{7} u_{2}^{7} u_{4}^{8} u_{5}^{11} + u_{1}^{7} u_{3}^{7} u_{4}^{8} u_{5}^{11} + u_{2}^{7} u_{3}^{7} u_{4}^{8} u_{5}^{11} \\
    &\quad + u_{1}^{7} u_{2} u_{4}^{14} u_{5}^{11} + u_{1}^{3} u_{2}^{5} u_{4}^{14} u_{5}^{11} + u_{1} u_{2}^{7} u_{4}^{14} u_{5}^{11} + u_{1}^{7} u_{3} u_{4}^{14} u_{5}^{11} \\
    &\quad + u_{2}^{7} u_{3} u_{4}^{14} u_{5}^{11} + u_{1}^{3} u_{3}^{5} u_{4}^{14} u_{5}^{11} + u_{2}^{3} u_{3}^{5} u_{4}^{14} u_{5}^{11} + u_{1} u_{3}^{7} u_{4}^{14} u_{5}^{11} \\
    &\quad + u_{2} u_{3}^{7} u_{4}^{14} u_{5}^{11} + u_{1}^{7} u_{2} u_{3}^{11} u_{5}^{14} + u_{1}^{3} u_{2}^{5} u_{3}^{11} u_{5}^{14} + u_{1} u_{2}^{7} u_{3}^{11} u_{5}^{14} \\
    &\quad + u_{1}^{7} u_{2} u_{4}^{11} u_{5}^{14} + u_{1}^{3} u_{2}^{5} u_{4}^{11} u_{5}^{14} + u_{1} u_{2}^{7} u_{4}^{11} u_{5}^{14} + u_{1}^{7} u_{3} u_{4}^{11} u_{5}^{14} \\
    &\quad + u_{2}^{7} u_{3} u_{4}^{11} u_{5}^{14} + u_{1}^{3} u_{3}^{5} u_{4}^{11} u_{5}^{14} + u_{2}^{3} u_{3}^{5} u_{4}^{11} u_{5}^{14} + u_{1} u_{3}^{7} u_{4}^{11} u_{5}^{14} \\
    &\quad + u_{2} u_{3}^{7} u_{4}^{11} u_{5}^{14} \\
\mathsf{Sigma5}[13] &= u_{1}^{7} u_{2}^{7} u_{3}^{8} u_{4}^{3} u_{5}^{8} + u_{1}^{7} u_{2}^{3} u_{3}^{11} u_{4}^{4} u_{5}^{8} + u_{1}^{3} u_{2}^{7} u_{3}^{11} u_{4}^{4} u_{5}^{8} + u_{1}^{7} u_{2}^{9} u_{3}^{3} u_{4}^{6} u_{5}^{8} \\
    &\quad + u_{1}^{3} u_{2}^{13} u_{3}^{3} u_{4}^{6} u_{5}^{8} + u_{1}^{7} u_{2}^{3} u_{3}^{9} u_{4}^{6} u_{5}^{8} + u_{1}^{3} u_{2}^{3} u_{3}^{13} u_{4}^{6} u_{5}^{8} + u_{1}^{3} u_{2}^{5} u_{3}^{10} u_{4}^{7} u_{5}^{8} \\
    &\quad + u_{1}^{7} u_{2}^{7} u_{3}^{3} u_{4}^{8} u_{5}^{8} + u_{1}^{7} u_{2}^{3} u_{3}^{5} u_{4}^{10} u_{5}^{8} + u_{1}^{3} u_{2}^{7} u_{3}^{5} u_{4}^{10} u_{5}^{8} + u_{1}^{3} u_{2}^{5} u_{3}^{7} u_{4}^{10} u_{5}^{8} \\
    &\quad + u_{1}^{7} u_{2} u_{3}^{6} u_{4}^{11} u_{5}^{8} + u_{1} u_{2}^{7} u_{3}^{6} u_{4}^{11} u_{5}^{8} + u_{1}^{7} u_{2}^{3} u_{3}^{3} u_{4}^{12} u_{5}^{8} + u_{1}^{3} u_{2}^{3} u_{3}^{7} u_{4}^{12} u_{5}^{8} \\
    &\quad + u_{1}^{7} u_{2} u_{3}^{3} u_{4}^{14} u_{5}^{8} + u_{1} u_{2}^{7} u_{3}^{3} u_{4}^{14} u_{5}^{8} + u_{1}^{3} u_{2}^{3} u_{3}^{5} u_{4}^{14} u_{5}^{8} + u_{1}^{7} u_{2} u_{3}^{14} u_{4} u_{5}^{10} \\
    &\quad + u_{1}^{3} u_{2}^{5} u_{3}^{14} u_{4} u_{5}^{10} + u_{1} u_{2}^{7} u_{3}^{14} u_{4} u_{5}^{10} + u_{1}^{7} u_{2}^{9} u_{3}^{3} u_{4}^{4} u_{5}^{10} + u_{1}^{3} u_{2}^{13} u_{3}^{3} u_{4}^{4} u_{5}^{10} \\
    &\quad + u_{1}^{7} u_{2}^{3} u_{3}^{9} u_{4}^{4} u_{5}^{10} + u_{1}^{7} u_{2}^{9} u_{3}^{2} u_{4}^{5} u_{5}^{10} + u_{1}^{3} u_{2}^{13} u_{3}^{2} u_{4}^{5} u_{5}^{10} + u_{1}^{7} u_{2}^{8} u_{3}^{3} u_{4}^{5} u_{5}^{10} \\
    &\quad + u_{1} u_{2}^{14} u_{3}^{3} u_{4}^{5} u_{5}^{10} + u_{1}^{7} u_{2}^{3} u_{3}^{8} u_{4}^{5} u_{5}^{10} + u_{1}^{7} u_{2} u_{3}^{10} u_{4}^{5} u_{5}^{10} + u_{1}^{3} u_{2} u_{3}^{14} u_{4}^{5} u_{5}^{10} \\
    &\quad + u_{1} u_{2}^{3} u_{3}^{14} u_{4}^{5} u_{5}^{10} + u_{1}^{7} u_{2}^{9} u_{3} u_{4}^{6} u_{5}^{10} + u_{1}^{3} u_{2}^{13} u_{3} u_{4}^{6} u_{5}^{10} + u_{1}^{3} u_{2} u_{3}^{12} u_{4}^{7} u_{5}^{10} \\
    &\quad + u_{1} u_{2}^{3} u_{3}^{12} u_{4}^{7} u_{5}^{10} + u_{1} u_{2}^{2} u_{3}^{13} u_{4}^{7} u_{5}^{10} + u_{1} u_{2} u_{3}^{14} u_{4}^{7} u_{5}^{10} + u_{1}^{7} u_{2}^{3} u_{3}^{5} u_{4}^{8} u_{5}^{10} \\
    &\quad + u_{1}^{3} u_{2}^{7} u_{3}^{5} u_{4}^{8} u_{5}^{10} + u_{1}^{7} u_{2} u_{3}^{7} u_{4}^{8} u_{5}^{10} + u_{1} u_{2}^{7} u_{3}^{7} u_{4}^{8} u_{5}^{10} + u_{1}^{7} u_{2}^{3} u_{3}^{4} u_{4}^{9} u_{5}^{10} \\
    &\quad + u_{1}^{3} u_{2}^{7} u_{3}^{4} u_{4}^{9} u_{5}^{10} + u_{1}^{7} u_{2} u_{3}^{6} u_{4}^{9} u_{5}^{10} + u_{1} u_{2}^{7} u_{3}^{6} u_{4}^{9} u_{5}^{10} + u_{1}^{3} u_{2}^{4} u_{3}^{7} u_{4}^{9} u_{5}^{10} \\
    &\quad + u_{1} u_{2}^{6} u_{3}^{7} u_{4}^{9} u_{5}^{10} + u_{1}^{3} u_{2}^{7} u_{3} u_{4}^{12} u_{5}^{10} + u_{1}^{7} u_{2} u_{3}^{3} u_{4}^{12} u_{5}^{10} + u_{1}^{3} u_{2} u_{3}^{7} u_{4}^{12} u_{5}^{10} \\
    &\quad + u_{1}^{3} u_{2} u_{3}^{6} u_{4}^{13} u_{5}^{10} + u_{1} u_{2}^{3} u_{3}^{6} u_{4}^{13} u_{5}^{10} + u_{1} u_{2}^{2} u_{3}^{7} u_{4}^{13} u_{5}^{10} + u_{1}^{7} u_{2} u_{3} u_{4}^{14} u_{5}^{10} \\
    &\quad + u_{1}^{3} u_{2}^{5} u_{3} u_{4}^{14} u_{5}^{10} + u_{1} u_{2}^{7} u_{3} u_{4}^{14} u_{5}^{10} + u_{1}^{3} u_{2} u_{3}^{5} u_{4}^{14} u_{5}^{10} + u_{1} u_{2}^{3} u_{3}^{5} u_{4}^{14} u_{5}^{10} \\
    &\quad + u_{1} u_{2} u_{3}^{7} u_{4}^{14} u_{5}^{10} + u_{1}^{3} u_{2}^{5} u_{3}^{10} u_{4}^{4} u_{5}^{11} + u_{1}^{3} u_{2} u_{3}^{12} u_{4}^{6} u_{5}^{11} + u_{1} u_{2}^{3} u_{3}^{12} u_{4}^{6} u_{5}^{11} \\
    &\quad + u_{1} u_{2}^{2} u_{3}^{13} u_{4}^{6} u_{5}^{11} + u_{1} u_{2} u_{3}^{14} u_{4}^{6} u_{5}^{11} + u_{1}^{7} u_{2} u_{3}^{6} u_{4}^{8} u_{5}^{11} + u_{1}^{3} u_{2}^{5} u_{3}^{6} u_{4}^{8} u_{5}^{11} \\
    &\quad + u_{1} u_{2}^{7} u_{3}^{6} u_{4}^{8} u_{5}^{11} + u_{1}^{7} u_{2} u_{3}^{2} u_{4}^{12} u_{5}^{11} + u_{1}^{3} u_{2}^{5} u_{3}^{2} u_{4}^{12} u_{5}^{11} + u_{1} u_{2}^{7} u_{3}^{2} u_{4}^{12} u_{5}^{11} \\
    &\quad + u_{1}^{3} u_{2} u_{3}^{4} u_{4}^{14} u_{5}^{11} + u_{1} u_{2}^{3} u_{3}^{4} u_{4}^{14} u_{5}^{11} + u_{1} u_{2}^{2} u_{3}^{5} u_{4}^{14} u_{5}^{11} + u_{1}^{7} u_{2} u_{3}^{11} u_{4}^{2} u_{5}^{12} \\
    &\quad + u_{1}^{3} u_{2}^{5} u_{3}^{11} u_{4}^{2} u_{5}^{12} + u_{1} u_{2}^{7} u_{3}^{11} u_{4}^{2} u_{5}^{12} + u_{1}^{3} u_{2}^{3} u_{3}^{12} u_{4}^{3} u_{5}^{12} + u_{1}^{3} u_{2}^{3} u_{3}^{7} u_{4}^{8} u_{5}^{12} \\
    &\quad + u_{1}^{7} u_{2} u_{3}^{2} u_{4}^{11} u_{5}^{12} + u_{1}^{3} u_{2}^{5} u_{3}^{2} u_{4}^{11} u_{5}^{12} + u_{1} u_{2}^{7} u_{3}^{2} u_{4}^{11} u_{5}^{12} + u_{1}^{3} u_{2} u_{3}^{3} u_{4}^{14} u_{5}^{12} \\
    &\quad + u_{1} u_{2}^{3} u_{3}^{3} u_{4}^{14} u_{5}^{12} + u_{1}^{3} u_{2}^{3} u_{3}^{12} u_{4} u_{5}^{14} + u_{1}^{7} u_{2}^{9} u_{3} u_{4}^{2} u_{5}^{14} + u_{1}^{3} u_{2}^{13} u_{3} u_{4}^{2} u_{5}^{14} \\
    &\quad + u_{1}^{3} u_{2}^{5} u_{3}^{9} u_{4}^{2} u_{5}^{14} + u_{1}^{3} u_{2} u_{3}^{13} u_{4}^{2} u_{5}^{14} + u_{1} u_{2}^{3} u_{3}^{13} u_{4}^{2} u_{5}^{14} + u_{1}^{7} u_{2}^{3} u_{3} u_{4}^{8} u_{5}^{14} \\
    &\quad + u_{1}^{7} u_{2} u_{3}^{3} u_{4}^{8} u_{5}^{14} + u_{1}^{3} u_{2}^{5} u_{3}^{3} u_{4}^{8} u_{5}^{14} + u_{1} u_{2}^{7} u_{3}^{3} u_{4}^{8} u_{5}^{14} + u_{1}^{3} u_{2} u_{3}^{6} u_{4}^{9} u_{5}^{14} \\
    &\quad + u_{1} u_{2}^{3} u_{3}^{6} u_{4}^{9} u_{5}^{14} + u_{1} u_{2}^{2} u_{3}^{7} u_{4}^{9} u_{5}^{14} + u_{1}^{7} u_{2} u_{3} u_{4}^{10} u_{5}^{14} + u_{1}^{3} u_{2}^{5} u_{3} u_{4}^{10} u_{5}^{14} \\
    &\quad + u_{1} u_{2}^{7} u_{3} u_{4}^{10} u_{5}^{14} + u_{1}^{3} u_{2} u_{3}^{4} u_{4}^{11} u_{5}^{14} + u_{1} u_{2}^{3} u_{3}^{4} u_{4}^{11} u_{5}^{14} + u_{1} u_{2}^{2} u_{3}^{5} u_{4}^{11} u_{5}^{14} \\
    &\quad + u_{1}^{3} u_{2} u_{3}^{3} u_{4}^{12} u_{5}^{14} + u_{1} u_{2}^{3} u_{3}^{3} u_{4}^{12} u_{5}^{14} + u_{1} u_{2} u_{3}^{3} u_{4}^{14} u_{5}^{14} \\
\mathsf{Sigma5}[14] &= u_{1}^{7} u_{2} u_{3}^{9} u_{4}^{14} u_{5}^{2} + u_{1} u_{2}^{7} u_{3}^{9} u_{4}^{14} u_{5}^{2} + u_{1}^{3} u_{2} u_{3}^{13} u_{4}^{14} u_{5}^{2} + u_{1} u_{2}^{3} u_{3}^{13} u_{4}^{14} u_{5}^{2} \\
    &\quad + u_{1}^{3} u_{2}^{13} u_{3}^{2} u_{4}^{12} u_{5}^{3} + u_{1}^{3} u_{2} u_{3}^{14} u_{4}^{12} u_{5}^{3} + u_{1} u_{2}^{3} u_{3}^{14} u_{4}^{12} u_{5}^{3} + u_{1}^{3} u_{2} u_{3}^{12} u_{4}^{14} u_{5}^{3} \\
    &\quad + u_{1} u_{2}^{3} u_{3}^{12} u_{4}^{14} u_{5}^{3} + u_{1} u_{2}^{2} u_{3}^{13} u_{4}^{14} u_{5}^{3} + u_{1}^{7} u_{2}^{7} u_{3}^{8} u_{4}^{3} u_{5}^{8} + u_{1}^{7} u_{2}^{11} u_{3}^{3} u_{4}^{4} u_{5}^{8} \\
    &\quad + u_{1}^{3} u_{2}^{5} u_{3}^{6} u_{4}^{11} u_{5}^{8} + u_{1}^{7} u_{2} u_{3}^{3} u_{4}^{14} u_{5}^{8} + u_{1}^{3} u_{2}^{5} u_{3}^{3} u_{4}^{14} u_{5}^{8} + u_{1} u_{2}^{7} u_{3}^{3} u_{4}^{14} u_{5}^{8} \\
    &\quad + u_{1}^{3} u_{2}^{3} u_{3}^{5} u_{4}^{14} u_{5}^{8} + u_{1}^{7} u_{2} u_{3}^{14} u_{4} u_{5}^{10} + u_{1}^{3} u_{2}^{5} u_{3}^{14} u_{4} u_{5}^{10} + u_{1} u_{2}^{7} u_{3}^{14} u_{4} u_{5}^{10} \\
    &\quad + u_{1}^{3} u_{2}^{3} u_{3}^{13} u_{4}^{4} u_{5}^{10} + u_{1}^{3} u_{2}^{13} u_{3} u_{4}^{6} u_{5}^{10} + u_{1}^{7} u_{2} u_{3}^{9} u_{4}^{6} u_{5}^{10} + u_{1} u_{2}^{7} u_{3}^{9} u_{4}^{6} u_{5}^{10} \\
    &\quad + u_{1}^{7} u_{2} u_{3}^{7} u_{4}^{8} u_{5}^{10} + u_{1}^{3} u_{2}^{5} u_{3}^{7} u_{4}^{8} u_{5}^{10} + u_{1} u_{2}^{7} u_{3}^{7} u_{4}^{8} u_{5}^{10} + u_{1}^{3} u_{2}^{7} u_{3} u_{4}^{12} u_{5}^{10} \\
    &\quad + u_{1}^{7} u_{2} u_{3} u_{4}^{14} u_{5}^{10} + u_{1} u_{2}^{7} u_{3} u_{4}^{14} u_{5}^{10} + u_{1} u_{2} u_{3}^{7} u_{4}^{14} u_{5}^{10} + u_{1} u_{2}^{2} u_{3}^{12} u_{4}^{7} u_{5}^{11} \\
    &\quad + u_{1}^{3} u_{2}^{5} u_{3}^{6} u_{4}^{8} u_{5}^{11} + u_{1}^{7} u_{2} u_{3}^{2} u_{4}^{12} u_{5}^{11} + u_{1} u_{2}^{7} u_{3}^{2} u_{4}^{12} u_{5}^{11} + u_{1} u_{2}^{2} u_{3}^{7} u_{4}^{12} u_{5}^{11} \\
    &\quad + u_{1}^{3} u_{2} u_{3}^{4} u_{4}^{14} u_{5}^{11} + u_{1} u_{2}^{3} u_{3}^{4} u_{4}^{14} u_{5}^{11} + u_{1} u_{2} u_{3}^{6} u_{4}^{14} u_{5}^{11} + u_{1}^{7} u_{2} u_{3}^{11} u_{4}^{2} u_{5}^{12} \\
    &\quad + u_{1}^{3} u_{2}^{5} u_{3}^{11} u_{4}^{2} u_{5}^{12} + u_{1} u_{2}^{7} u_{3}^{11} u_{4}^{2} u_{5}^{12} + u_{1}^{3} u_{2}^{13} u_{3}^{2} u_{4}^{3} u_{5}^{12} + u_{1}^{3} u_{2}^{3} u_{3}^{12} u_{4}^{3} u_{5}^{12} \\
    &\quad + u_{1}^{3} u_{2} u_{3}^{14} u_{4}^{3} u_{5}^{12} + u_{1} u_{2}^{3} u_{3}^{14} u_{4}^{3} u_{5}^{12} + u_{1}^{3} u_{2}^{7} u_{3} u_{4}^{10} u_{5}^{12} + u_{1}^{7} u_{2} u_{3}^{3} u_{4}^{10} u_{5}^{12} \\
    &\quad + u_{1} u_{2}^{7} u_{3}^{3} u_{4}^{10} u_{5}^{12} + u_{1}^{7} u_{2} u_{3}^{2} u_{4}^{11} u_{5}^{12} + u_{1} u_{2}^{7} u_{3}^{2} u_{4}^{11} u_{5}^{12} + u_{1}^{3} u_{2}^{3} u_{3}^{4} u_{4}^{11} u_{5}^{12} \\
    &\quad + u_{1} u_{2}^{2} u_{3}^{7} u_{4}^{11} u_{5}^{12} + u_{1}^{3} u_{2}^{3} u_{3}^{12} u_{4} u_{5}^{14} + u_{1}^{3} u_{2}^{5} u_{3}^{9} u_{4}^{2} u_{5}^{14} + u_{1}^{3} u_{2} u_{3}^{13} u_{4}^{2} u_{5}^{14} \\
    &\quad + u_{1} u_{2}^{3} u_{3}^{13} u_{4}^{2} u_{5}^{14} + u_{1}^{3} u_{2} u_{3}^{12} u_{4}^{3} u_{5}^{14} + u_{1} u_{2}^{3} u_{3}^{12} u_{4}^{3} u_{5}^{14} + u_{1} u_{2}^{2} u_{3}^{13} u_{4}^{3} u_{5}^{14} \\
    &\quad + u_{1}^{7} u_{2} u_{3} u_{4}^{10} u_{5}^{14} + u_{1} u_{2}^{7} u_{3} u_{4}^{10} u_{5}^{14} + u_{1} u_{2} u_{3}^{7} u_{4}^{10} u_{5}^{14} + u_{1}^{3} u_{2} u_{3}^{4} u_{4}^{11} u_{5}^{14} \\
    &\quad + u_{1} u_{2}^{3} u_{3}^{4} u_{4}^{11} u_{5}^{14} + u_{1} u_{2}^{2} u_{3}^{3} u_{4}^{13} u_{5}^{14} + u_{1}^{3} u_{2} u_{3} u_{4}^{14} u_{5}^{14} + u_{1} u_{2}^{3} u_{3} u_{4}^{14} u_{5}^{14} \\
    &\quad + u_{1} u_{2} u_{3}^{3} u_{4}^{14} u_{5}^{14} \\
\mathsf{Sigma5}[15] &= u_{1}^{7} u_{2}^{11} u_{3}^{5} u_{4}^{8} u_{5}^{2} + u_{1}^{7} u_{2}^{3} u_{3}^{13} u_{4}^{8} u_{5}^{2} + u_{1}^{3} u_{2}^{7} u_{3}^{13} u_{4}^{8} u_{5}^{2} + u_{1}^{7} u_{2}^{7} u_{3}^{8} u_{4}^{9} u_{5}^{2} \\
    &\quad + u_{1}^{7} u_{2}^{11} u_{3} u_{4}^{12} u_{5}^{2} + u_{1}^{3} u_{2}^{13} u_{3}^{3} u_{4}^{12} u_{5}^{2} + u_{1}^{7} u_{2} u_{3}^{11} u_{4}^{12} u_{5}^{2} + u_{1} u_{2}^{7} u_{3}^{11} u_{4}^{12} u_{5}^{2} \\
    &\quad + u_{1}^{7} u_{2}^{3} u_{3}^{8} u_{4}^{13} u_{5}^{2} + u_{1}^{3} u_{2}^{7} u_{3}^{8} u_{4}^{13} u_{5}^{2} + u_{1}^{7} u_{2} u_{3}^{10} u_{4}^{13} u_{5}^{2} + u_{1} u_{2}^{7} u_{3}^{10} u_{4}^{13} u_{5}^{2} \\
    &\quad + u_{1}^{3} u_{2}^{4} u_{3}^{11} u_{4}^{13} u_{5}^{2} + u_{1} u_{2}^{6} u_{3}^{11} u_{4}^{13} u_{5}^{2} + u_{1}^{3} u_{2} u_{3}^{14} u_{4}^{13} u_{5}^{2} + u_{1}^{3} u_{2}^{13} u_{3} u_{4}^{14} u_{5}^{2} \\
    &\quad + u_{1}^{3} u_{2}^{5} u_{3}^{9} u_{4}^{14} u_{5}^{2} + u_{1}^{3} u_{2} u_{3}^{13} u_{4}^{14} u_{5}^{2} + u_{1} u_{2}^{3} u_{3}^{13} u_{4}^{14} u_{5}^{2} + u_{1}^{7} u_{2}^{11} u_{3}^{5} u_{4}^{2} u_{5}^{8} \\
    &\quad + u_{1}^{7} u_{2}^{3} u_{3}^{13} u_{4}^{2} u_{5}^{8} + u_{1}^{3} u_{2}^{7} u_{3}^{13} u_{4}^{2} u_{5}^{8} + u_{1}^{7} u_{2}^{7} u_{3}^{3} u_{4}^{8} u_{5}^{8} + u_{1}^{7} u_{2} u_{3}^{7} u_{4}^{10} u_{5}^{8} \\
    &\quad + u_{1} u_{2}^{7} u_{3}^{7} u_{4}^{10} u_{5}^{8} + u_{1}^{3} u_{2}^{3} u_{3}^{7} u_{4}^{12} u_{5}^{8} + u_{1}^{7} u_{2}^{3} u_{3} u_{4}^{14} u_{5}^{8} + u_{1}^{3} u_{2}^{7} u_{3} u_{4}^{14} u_{5}^{8} \\
    &\quad + u_{1}^{7} u_{2} u_{3}^{3} u_{4}^{14} u_{5}^{8} + u_{1}^{3} u_{2}^{5} u_{3}^{3} u_{4}^{14} u_{5}^{8} + u_{1} u_{2}^{7} u_{3}^{3} u_{4}^{14} u_{5}^{8} + u_{1}^{3} u_{2}^{3} u_{3}^{5} u_{4}^{14} u_{5}^{8} \\
    &\quad + u_{1}^{3} u_{2} u_{3}^{7} u_{4}^{14} u_{5}^{8} + u_{1} u_{2}^{3} u_{3}^{7} u_{4}^{14} u_{5}^{8} + u_{1}^{7} u_{2}^{7} u_{3}^{8} u_{4} u_{5}^{10} + u_{1}^{7} u_{2}^{9} u_{3} u_{4}^{6} u_{5}^{10} \\
    &\quad + u_{1}^{3} u_{2} u_{3}^{13} u_{4}^{6} u_{5}^{10} + u_{1} u_{2}^{3} u_{3}^{13} u_{4}^{6} u_{5}^{10} + u_{1}^{7} u_{2} u_{3}^{7} u_{4}^{8} u_{5}^{10} + u_{1} u_{2}^{7} u_{3}^{7} u_{4}^{8} u_{5}^{10} \\
    &\quad + u_{1}^{3} u_{2}^{4} u_{3}^{7} u_{4}^{9} u_{5}^{10} + u_{1} u_{2}^{6} u_{3}^{7} u_{4}^{9} u_{5}^{10} + u_{1} u_{2}^{3} u_{3}^{7} u_{4}^{12} u_{5}^{10} + u_{1}^{7} u_{2} u_{3} u_{4}^{14} u_{5}^{10} \\
    &\quad + u_{1} u_{2}^{7} u_{3} u_{4}^{14} u_{5}^{10} + u_{1} u_{2} u_{3}^{7} u_{4}^{14} u_{5}^{10} + u_{1}^{7} u_{2}^{11} u_{3} u_{4}^{2} u_{5}^{12} + u_{1}^{3} u_{2}^{13} u_{3}^{3} u_{4}^{2} u_{5}^{12} \\
    &\quad + u_{1}^{7} u_{2} u_{3}^{11} u_{4}^{2} u_{5}^{12} + u_{1} u_{2}^{7} u_{3}^{11} u_{4}^{2} u_{5}^{12} + u_{1}^{3} u_{2}^{3} u_{3}^{7} u_{4}^{8} u_{5}^{12} + u_{1}^{7} u_{2}^{3} u_{3} u_{4}^{10} u_{5}^{12} \\
    &\quad + u_{1}^{3} u_{2}^{3} u_{3}^{3} u_{4}^{12} u_{5}^{12} + u_{1}^{7} u_{2}^{3} u_{3}^{8} u_{4} u_{5}^{14} + u_{1}^{3} u_{2}^{7} u_{3}^{8} u_{4} u_{5}^{14} + u_{1}^{7} u_{2} u_{3}^{10} u_{4} u_{5}^{14} \\
    &\quad + u_{1} u_{2}^{7} u_{3}^{10} u_{4} u_{5}^{14} + u_{1}^{3} u_{2}^{4} u_{3}^{11} u_{4} u_{5}^{14} + u_{1} u_{2}^{6} u_{3}^{11} u_{4} u_{5}^{14} + u_{1}^{3} u_{2} u_{3}^{14} u_{4} u_{5}^{14} \\
    &\quad + u_{1}^{3} u_{2}^{13} u_{3} u_{4}^{2} u_{5}^{14} + u_{1}^{7} u_{2} u_{3}^{9} u_{4}^{2} u_{5}^{14} + u_{1}^{3} u_{2}^{5} u_{3}^{9} u_{4}^{2} u_{5}^{14} + u_{1} u_{2}^{7} u_{3}^{9} u_{4}^{2} u_{5}^{14} \\
    &\quad + u_{1}^{3} u_{2}^{5} u_{3}^{8} u_{4}^{3} u_{5}^{14} + u_{1}^{3} u_{2} u_{3}^{12} u_{4}^{3} u_{5}^{14} + u_{1} u_{2}^{3} u_{3}^{12} u_{4}^{3} u_{5}^{14} + u_{1} u_{2} u_{3}^{14} u_{4}^{3} u_{5}^{14} \\
    &\quad + u_{1}^{7} u_{2}^{3} u_{3} u_{4}^{8} u_{5}^{14} + u_{1}^{3} u_{2}^{7} u_{3} u_{4}^{8} u_{5}^{14} + u_{1}^{7} u_{2} u_{3}^{3} u_{4}^{8} u_{5}^{14} + u_{1}^{3} u_{2}^{5} u_{3}^{3} u_{4}^{8} u_{5}^{14} \\
    &\quad + u_{1} u_{2}^{7} u_{3}^{3} u_{4}^{8} u_{5}^{14} + u_{1}^{3} u_{2}^{3} u_{3}^{5} u_{4}^{8} u_{5}^{14} + u_{1}^{3} u_{2} u_{3}^{7} u_{4}^{8} u_{5}^{14} + u_{1}^{3} u_{2}^{5} u_{3}^{2} u_{4}^{9} u_{5}^{14} \\
    &\quad + u_{1}^{3} u_{2}^{3} u_{3}^{4} u_{4}^{9} u_{5}^{14} + u_{1}^{7} u_{2} u_{3} u_{4}^{10} u_{5}^{14} + u_{1} u_{2}^{7} u_{3} u_{4}^{10} u_{5}^{14} + u_{1} u_{2} u_{3}^{7} u_{4}^{10} u_{5}^{14} \\
    &\quad + u_{1}^{3} u_{2}^{4} u_{3} u_{4}^{11} u_{5}^{14} + u_{1} u_{2}^{6} u_{3} u_{4}^{11} u_{5}^{14} + u_{1}^{3} u_{2} u_{3}^{4} u_{4}^{11} u_{5}^{14} + u_{1} u_{2} u_{3}^{6} u_{4}^{11} u_{5}^{14} \\
    &\quad + u_{1}^{3} u_{2}^{3} u_{3} u_{4}^{12} u_{5}^{14} + u_{1} u_{2}^{3} u_{3} u_{4}^{14} u_{5}^{14} + u_{1} u_{2} u_{3}^{3} u_{4}^{14} u_{5}^{14} \\
\mathsf{Sigma5}[16] &= u_{1}^{7} u_{2}^{11} u_{3}^{13} u_{4}^{2} + u_{1}^{7} u_{2}^{7} u_{3}^{9} u_{4}^{10} + u_{1}^{7} u_{2}^{3} u_{3}^{11} u_{4}^{12} + u_{1}^{3} u_{2}^{7} u_{3}^{11} u_{4}^{12} \\
    &\quad + u_{1}^{7} u_{2}^{11} u_{3} u_{4}^{14} + u_{1}^{7} u_{2} u_{3}^{11} u_{4}^{14} + u_{1} u_{2}^{7} u_{3}^{11} u_{4}^{14} + u_{1}^{3} u_{2}^{3} u_{3}^{13} u_{4}^{14} \\
    &\quad + u_{1}^{7} u_{2}^{11} u_{3}^{13} u_{5}^{2} + u_{1}^{7} u_{2}^{11} u_{4}^{13} u_{5}^{2} + u_{1}^{7} u_{3}^{11} u_{4}^{13} u_{5}^{2} + u_{2}^{7} u_{3}^{11} u_{4}^{13} u_{5}^{2} \\
    &\quad + u_{1}^{7} u_{2}^{7} u_{3}^{9} u_{5}^{10} + u_{1}^{7} u_{2}^{7} u_{4}^{9} u_{5}^{10} + u_{1}^{7} u_{3}^{7} u_{4}^{9} u_{5}^{10} + u_{2}^{7} u_{3}^{7} u_{4}^{9} u_{5}^{10} \\
    &\quad + u_{1}^{7} u_{2}^{3} u_{3}^{11} u_{5}^{12} + u_{1}^{3} u_{2}^{7} u_{3}^{11} u_{5}^{12} + u_{1}^{7} u_{2}^{3} u_{4}^{11} u_{5}^{12} + u_{1}^{3} u_{2}^{7} u_{4}^{11} u_{5}^{12} \\
    &\quad + u_{1}^{7} u_{3}^{3} u_{4}^{11} u_{5}^{12} + u_{2}^{7} u_{3}^{3} u_{4}^{11} u_{5}^{12} + u_{1}^{3} u_{3}^{7} u_{4}^{11} u_{5}^{12} + u_{2}^{3} u_{3}^{7} u_{4}^{11} u_{5}^{12} \\
    &\quad + u_{1}^{7} u_{2}^{11} u_{3} u_{5}^{14} + u_{1}^{7} u_{2} u_{3}^{11} u_{5}^{14} + u_{1} u_{2}^{7} u_{3}^{11} u_{5}^{14} + u_{1}^{3} u_{2}^{3} u_{3}^{13} u_{5}^{14} \\
    &\quad + u_{1}^{7} u_{2}^{11} u_{4} u_{5}^{14} + u_{1}^{7} u_{3}^{11} u_{4} u_{5}^{14} + u_{2}^{7} u_{3}^{11} u_{4} u_{5}^{14} + u_{1}^{7} u_{2} u_{4}^{11} u_{5}^{14} \\
    &\quad + u_{1} u_{2}^{7} u_{4}^{11} u_{5}^{14} + u_{1}^{7} u_{3} u_{4}^{11} u_{5}^{14} + u_{2}^{7} u_{3} u_{4}^{11} u_{5}^{14} + u_{1} u_{3}^{7} u_{4}^{11} u_{5}^{14} \\
    &\quad + u_{2} u_{3}^{7} u_{4}^{11} u_{5}^{14} + u_{1}^{3} u_{2}^{3} u_{4}^{13} u_{5}^{14} + u_{1}^{3} u_{3}^{3} u_{4}^{13} u_{5}^{14} + u_{2}^{3} u_{3}^{3} u_{4}^{13} u_{5}^{14} \\
\mathsf{Sigma5}[17] &= u_{1}^{15} u_{2} u_{3}^{2} u_{4}^{12} u_{5}^{3} + u_{1} u_{2}^{15} u_{3}^{2} u_{4}^{12} u_{5}^{3} + u_{1} u_{2}^{2} u_{3}^{15} u_{4}^{12} u_{5}^{3} + u_{1} u_{2}^{2} u_{3}^{12} u_{4}^{15} u_{5}^{3} \\
    &\quad + u_{1}^{15} u_{2}^{3} u_{3}^{3} u_{4}^{4} u_{5}^{8} + u_{1}^{3} u_{2}^{15} u_{3}^{3} u_{4}^{4} u_{5}^{8} + u_{1}^{3} u_{2}^{3} u_{3}^{15} u_{4}^{4} u_{5}^{8} + u_{1}^{3} u_{2}^{3} u_{3}^{4} u_{4}^{15} u_{5}^{8} \\
    &\quad + u_{1}^{15} u_{2} u_{3} u_{4}^{6} u_{5}^{10} + u_{1} u_{2}^{15} u_{3} u_{4}^{6} u_{5}^{10} + u_{1} u_{2} u_{3}^{15} u_{4}^{6} u_{5}^{10} + u_{1} u_{2} u_{3}^{6} u_{4}^{15} u_{5}^{10} \\
    &\quad + u_{1}^{15} u_{2} u_{3}^{2} u_{4}^{3} u_{5}^{12} + u_{1} u_{2}^{15} u_{3}^{2} u_{4}^{3} u_{5}^{12} + u_{1} u_{2}^{2} u_{3}^{15} u_{4}^{3} u_{5}^{12} + u_{1} u_{2}^{2} u_{3}^{3} u_{4}^{15} u_{5}^{12} \\
    &\quad + u_{1} u_{2}^{2} u_{3}^{12} u_{4}^{3} u_{5}^{15} + u_{1}^{3} u_{2}^{3} u_{3}^{4} u_{4}^{8} u_{5}^{15} + u_{1} u_{2} u_{3}^{6} u_{4}^{10} u_{5}^{15} + u_{1} u_{2}^{2} u_{3}^{3} u_{4}^{12} u_{5}^{15} \\
\mathsf{Sigma5}[18] &= u_{1}^{15} u_{2}^{3} u_{3}^{5} u_{4}^{8} u_{5}^{2} + u_{1}^{3} u_{2}^{15} u_{3}^{5} u_{4}^{8} u_{5}^{2} + u_{1}^{3} u_{2}^{5} u_{3}^{15} u_{4}^{8} u_{5}^{2} + u_{1}^{15} u_{2}^{3} u_{3} u_{4}^{12} u_{5}^{2} \\
    &\quad + u_{1}^{3} u_{2}^{15} u_{3} u_{4}^{12} u_{5}^{2} + u_{1}^{15} u_{2} u_{3}^{3} u_{4}^{12} u_{5}^{2} + u_{1} u_{2}^{15} u_{3}^{3} u_{4}^{12} u_{5}^{2} + u_{1}^{3} u_{2} u_{3}^{15} u_{4}^{12} u_{5}^{2} \\
    &\quad + u_{1} u_{2}^{3} u_{3}^{15} u_{4}^{12} u_{5}^{2} + u_{1}^{15} u_{2} u_{3} u_{4}^{14} u_{5}^{2} + u_{1} u_{2}^{15} u_{3} u_{4}^{14} u_{5}^{2} + u_{1} u_{2} u_{3}^{15} u_{4}^{14} u_{5}^{2} \\
    &\quad + u_{1}^{3} u_{2}^{5} u_{3}^{8} u_{4}^{15} u_{5}^{2} + u_{1}^{3} u_{2} u_{3}^{12} u_{4}^{15} u_{5}^{2} + u_{1} u_{2}^{3} u_{3}^{12} u_{4}^{15} u_{5}^{2} + u_{1} u_{2} u_{3}^{14} u_{4}^{15} u_{5}^{2} \\
    &\quad + u_{1}^{15} u_{2}^{3} u_{3}^{5} u_{4}^{2} u_{5}^{8} + u_{1}^{3} u_{2}^{15} u_{3}^{5} u_{4}^{2} u_{5}^{8} + u_{1}^{3} u_{2}^{5} u_{3}^{15} u_{4}^{2} u_{5}^{8} + u_{1}^{3} u_{2}^{5} u_{3}^{2} u_{4}^{15} u_{5}^{8} \\
    &\quad + u_{1}^{15} u_{2} u_{3} u_{4}^{6} u_{5}^{10} + u_{1} u_{2}^{15} u_{3} u_{4}^{6} u_{5}^{10} + u_{1} u_{2} u_{3}^{15} u_{4}^{6} u_{5}^{10} + u_{1} u_{2} u_{3}^{6} u_{4}^{15} u_{5}^{10} \\
    &\quad + u_{1}^{15} u_{2}^{3} u_{3} u_{4}^{2} u_{5}^{12} + u_{1}^{3} u_{2}^{15} u_{3} u_{4}^{2} u_{5}^{12} + u_{1}^{15} u_{2} u_{3}^{3} u_{4}^{2} u_{5}^{12} + u_{1} u_{2}^{15} u_{3}^{3} u_{4}^{2} u_{5}^{12} \\
    &\quad + u_{1}^{3} u_{2} u_{3}^{15} u_{4}^{2} u_{5}^{12} + u_{1} u_{2}^{3} u_{3}^{15} u_{4}^{2} u_{5}^{12} + u_{1}^{3} u_{2} u_{3}^{2} u_{4}^{15} u_{5}^{12} + u_{1} u_{2}^{3} u_{3}^{2} u_{4}^{15} u_{5}^{12} \\
    &\quad + u_{1}^{15} u_{2} u_{3} u_{4}^{2} u_{5}^{14} + u_{1} u_{2}^{15} u_{3} u_{4}^{2} u_{5}^{14} + u_{1} u_{2} u_{3}^{15} u_{4}^{2} u_{5}^{14} + u_{1} u_{2} u_{3}^{2} u_{4}^{15} u_{5}^{14} \\
    &\quad + u_{1}^{3} u_{2}^{5} u_{3}^{8} u_{4}^{2} u_{5}^{15} + u_{1}^{3} u_{2} u_{3}^{12} u_{4}^{2} u_{5}^{15} + u_{1} u_{2}^{3} u_{3}^{12} u_{4}^{2} u_{5}^{15} + u_{1} u_{2} u_{3}^{14} u_{4}^{2} u_{5}^{15} \\
    &\quad + u_{1}^{3} u_{2}^{5} u_{3}^{2} u_{4}^{8} u_{5}^{15} + u_{1} u_{2} u_{3}^{6} u_{4}^{10} u_{5}^{15} + u_{1}^{3} u_{2} u_{3}^{2} u_{4}^{12} u_{5}^{15} + u_{1} u_{2}^{3} u_{3}^{2} u_{4}^{12} u_{5}^{15} + u_{1} u_{2} u_{3}^{2} u_{4}^{14} u_{5}^{15} \\
\mathsf{Sigma5}[19] &= u_{1}^{15} u_{2} u_{3}^{14} u_{4}^{3} + u_{1} u_{2}^{15} u_{3}^{14} u_{4}^{3} + u_{1} u_{2}^{14} u_{3}^{15} u_{4}^{3} + u_{1}^{15} u_{2}^{3} u_{3}^{5} u_{4}^{10} \\
    &\quad + u_{1}^{3} u_{2}^{15} u_{3}^{5} u_{4}^{10} + u_{1}^{3} u_{2}^{5} u_{3}^{15} u_{4}^{10} + u_{1}^{15} u_{2} u_{3}^{3} u_{4}^{14} + u_{1} u_{2}^{15} u_{3}^{3} u_{4}^{14} \\
    &\quad + u_{1} u_{2}^{3} u_{3}^{15} u_{4}^{14} + u_{1} u_{2}^{14} u_{3}^{3} u_{4}^{15} + u_{1}^{3} u_{2}^{5} u_{3}^{10} u_{4}^{15} + u_{1} u_{2}^{3} u_{3}^{14} u_{4}^{15} \\
    &\quad + u_{1}^{15} u_{2} u_{3}^{14} u_{5}^{3} + u_{1} u_{2}^{15} u_{3}^{14} u_{5}^{3} + u_{1} u_{2}^{14} u_{3}^{15} u_{5}^{3} + u_{1}^{15} u_{2} u_{4}^{14} u_{5}^{3} \\
    &\quad + u_{1} u_{2}^{15} u_{4}^{14} u_{5}^{3} + u_{1}^{15} u_{3} u_{4}^{14} u_{5}^{3} + u_{2}^{15} u_{3} u_{4}^{14} u_{5}^{3} + u_{1} u_{3}^{15} u_{4}^{14} u_{5}^{3} \\
    &\quad + u_{2} u_{3}^{15} u_{4}^{14} u_{5}^{3} + u_{1} u_{2}^{14} u_{4}^{15} u_{5}^{3} + u_{1} u_{3}^{14} u_{4}^{15} u_{5}^{3} + u_{2} u_{3}^{14} u_{4}^{15} u_{5}^{3} \\
    &\quad + u_{1}^{15} u_{2}^{3} u_{3}^{5} u_{5}^{10} + u_{1}^{3} u_{2}^{15} u_{3}^{5} u_{5}^{10} + u_{1}^{3} u_{2}^{5} u_{3}^{15} u_{5}^{10} + u_{1}^{15} u_{2}^{3} u_{4}^{5} u_{5}^{10} \\
    &\quad + u_{1}^{3} u_{2}^{15} u_{4}^{5} u_{5}^{10} + u_{1}^{15} u_{3}^{3} u_{4}^{5} u_{5}^{10} + u_{2}^{15} u_{3}^{3} u_{4}^{5} u_{5}^{10} + u_{1}^{3} u_{3}^{15} u_{4}^{5} u_{5}^{10} \\
    &\quad + u_{2}^{3} u_{3}^{15} u_{4}^{5} u_{5}^{10} + u_{1}^{3} u_{2}^{5} u_{4}^{15} u_{5}^{10} + u_{1}^{3} u_{3}^{5} u_{4}^{15} u_{5}^{10} + u_{2}^{3} u_{3}^{5} u_{4}^{15} u_{5}^{10} \\
    &\quad + u_{1}^{15} u_{2} u_{3}^{3} u_{5}^{14} + u_{1} u_{2}^{15} u_{3}^{3} u_{5}^{14} + u_{1} u_{2}^{3} u_{3}^{15} u_{5}^{14} + u_{1}^{15} u_{2} u_{4}^{3} u_{5}^{14} \\
    &\quad + u_{1} u_{2}^{15} u_{4}^{3} u_{5}^{14} + u_{1}^{15} u_{3} u_{4}^{3} u_{5}^{14} + u_{2}^{15} u_{3} u_{4}^{3} u_{5}^{14} + u_{1} u_{3}^{15} u_{4}^{3} u_{5}^{14} \\
    &\quad + u_{2} u_{3}^{15} u_{4}^{3} u_{5}^{14} + u_{1} u_{2}^{3} u_{4}^{15} u_{5}^{14} + u_{1} u_{3}^{3} u_{4}^{15} u_{5}^{14} + u_{2} u_{3}^{3} u_{4}^{15} u_{5}^{14} \\
    &\quad + u_{1} u_{2}^{14} u_{3}^{3} u_{5}^{15} + u_{1}^{3} u_{2}^{5} u_{3}^{10} u_{5}^{15} + u_{1} u_{2}^{3} u_{3}^{14} u_{5}^{15} + u_{1} u_{2}^{14} u_{4}^{3} u_{5}^{15} \\
    &\quad + u_{1} u_{3}^{14} u_{4}^{3} u_{5}^{15} + u_{2} u_{3}^{14} u_{4}^{3} u_{5}^{15} + u_{1}^{3} u_{2}^{5} u_{4}^{10} u_{5}^{15} + u_{1}^{3} u_{3}^{5} u_{4}^{10} u_{5}^{15} \\
    &\quad + u_{2}^{3} u_{3}^{5} u_{4}^{10} u_{5}^{15} + u_{1} u_{2}^{3} u_{4}^{14} u_{5}^{15} + u_{1} u_{3}^{3} u_{4}^{14} u_{5}^{15} + u_{2} u_{3}^{3} u_{4}^{14} u_{5}^{15} \\
\mathsf{Sigma5}[20] &= u_{1}^{15} u_{2}^{3} u_{3}^{13} u_{4}^{2} + u_{1}^{3} u_{2}^{15} u_{3}^{13} u_{4}^{2} + u_{1}^{3} u_{2}^{13} u_{3}^{15} u_{4}^{2} + u_{1}^{15} u_{2} u_{3}^{14} u_{4}^{3} \\
    &\quad + u_{1} u_{2}^{15} u_{3}^{14} u_{4}^{3} + u_{1} u_{2}^{14} u_{3}^{15} u_{4}^{3} + u_{1}^{15} u_{2}^{3} u_{3}^{3} u_{4}^{12} + u_{1}^{3} u_{2}^{15} u_{3}^{3} u_{4}^{12} \\
    &\quad + u_{1}^{3} u_{2}^{3} u_{3}^{15} u_{4}^{12} + u_{1}^{15} u_{2}^{3} u_{3} u_{4}^{14} + u_{1}^{3} u_{2}^{15} u_{3} u_{4}^{14} + u_{1}^{3} u_{2} u_{3}^{15} u_{4}^{14} \\
    &\quad + u_{1}^{3} u_{2}^{13} u_{3}^{2} u_{4}^{15} + u_{1} u_{2}^{14} u_{3}^{3} u_{4}^{15} + u_{1}^{3} u_{2}^{3} u_{3}^{12} u_{4}^{15} + u_{1}^{3} u_{2} u_{3}^{14} u_{4}^{15} \\
    &\quad + u_{1}^{15} u_{2}^{3} u_{3}^{13} u_{5}^{2} + u_{1}^{3} u_{2}^{15} u_{3}^{13} u_{5}^{2} + u_{1}^{3} u_{2}^{13} u_{3}^{15} u_{5}^{2} + u_{1}^{15} u_{2}^{3} u_{4}^{13} u_{5}^{2} \\
    &\quad + u_{1}^{3} u_{2}^{15} u_{4}^{13} u_{5}^{2} + u_{1}^{15} u_{3}^{3} u_{4}^{13} u_{5}^{2} + u_{2}^{15} u_{3}^{3} u_{4}^{13} u_{5}^{2} + u_{1}^{3} u_{3}^{15} u_{4}^{13} u_{5}^{2} \\
    &\quad + u_{2}^{3} u_{3}^{15} u_{4}^{13} u_{5}^{2} + u_{1}^{3} u_{2}^{13} u_{4}^{15} u_{5}^{2} + u_{1}^{3} u_{3}^{13} u_{4}^{15} u_{5}^{2} + u_{2}^{3} u_{3}^{13} u_{4}^{15} u_{5}^{2} \\
    &\quad + u_{1}^{15} u_{2} u_{3}^{14} u_{5}^{3} + u_{1} u_{2}^{15} u_{3}^{14} u_{5}^{3} + u_{1} u_{2}^{14} u_{3}^{15} u_{5}^{3} + u_{1}^{15} u_{2} u_{4}^{14} u_{5}^{3} \\
    &\quad + u_{1} u_{2}^{15} u_{4}^{14} u_{5}^{3} + u_{1}^{15} u_{3} u_{4}^{14} u_{5}^{3} + u_{2}^{15} u_{3} u_{4}^{14} u_{5}^{3} + u_{1} u_{3}^{15} u_{4}^{14} u_{5}^{3} \\
    &\quad + u_{2} u_{3}^{15} u_{4}^{14} u_{5}^{3} + u_{1} u_{2}^{14} u_{4}^{15} u_{5}^{3} + u_{1} u_{3}^{14} u_{4}^{15} u_{5}^{3} + u_{2} u_{3}^{14} u_{4}^{15} u_{5}^{3} \\
    &\quad + u_{1}^{15} u_{2}^{3} u_{3}^{3} u_{5}^{12} + u_{1}^{3} u_{2}^{15} u_{3}^{3} u_{5}^{12} + u_{1}^{3} u_{2}^{3} u_{3}^{15} u_{5}^{12} + u_{1}^{15} u_{2}^{3} u_{4}^{3} u_{5}^{12} \\
    &\quad + u_{1}^{3} u_{2}^{15} u_{4}^{3} u_{5}^{12} + u_{1}^{15} u_{3}^{3} u_{4}^{3} u_{5}^{12} + u_{2}^{15} u_{3}^{3} u_{4}^{3} u_{5}^{12} + u_{1}^{3} u_{3}^{15} u_{4}^{3} u_{5}^{12} \\
    &\quad + u_{2}^{3} u_{3}^{15} u_{4}^{3} u_{5}^{12} + u_{1}^{3} u_{2}^{3} u_{4}^{15} u_{5}^{12} + u_{1}^{3} u_{3}^{3} u_{4}^{15} u_{5}^{12} + u_{2}^{3} u_{3}^{3} u_{4}^{15} u_{5}^{12} \\
    &\quad + u_{1}^{15} u_{2}^{3} u_{3} u_{5}^{14} + u_{1}^{3} u_{2}^{15} u_{3} u_{5}^{14} + u_{1}^{3} u_{2} u_{3}^{15} u_{5}^{14} + u_{1}^{15} u_{2}^{3} u_{4} u_{5}^{14} \\
    &\quad + u_{1}^{3} u_{2}^{15} u_{4} u_{5}^{14} + u_{1}^{15} u_{3}^{3} u_{4} u_{5}^{14} + u_{2}^{15} u_{3}^{3} u_{4} u_{5}^{14} + u_{1}^{3} u_{3}^{15} u_{4} u_{5}^{14} \\
    &\quad + u_{2}^{3} u_{3}^{15} u_{4} u_{5}^{14} + u_{1}^{3} u_{2} u_{4}^{15} u_{5}^{14} + u_{1}^{3} u_{3} u_{4}^{15} u_{5}^{14} + u_{2}^{3} u_{3} u_{4}^{15} u_{5}^{14} \\
    &\quad + u_{1}^{3} u_{2}^{13} u_{3}^{2} u_{5}^{15} + u_{1} u_{2}^{14} u_{3}^{3} u_{5}^{15} + u_{1}^{3} u_{2}^{3} u_{3}^{12} u_{5}^{15} + u_{1}^{3} u_{2} u_{3}^{14} u_{5}^{15} \\
    &\quad + u_{1}^{3} u_{2}^{13} u_{4}^{2} u_{5}^{15} + u_{1}^{3} u_{3}^{13} u_{4}^{2} u_{5}^{15} + u_{2}^{3} u_{3}^{13} u_{4}^{2} u_{5}^{15} + u_{1} u_{2}^{14} u_{4}^{3} u_{5}^{15} \\
    &\quad + u_{1} u_{3}^{14} u_{4}^{3} u_{5}^{15} + u_{2} u_{3}^{14} u_{4}^{3} u_{5}^{15} + u_{1}^{3} u_{2}^{3} u_{4}^{12} u_{5}^{15} + u_{1}^{3} u_{3}^{3} u_{4}^{12} u_{5}^{15} \\
    &\quad + u_{2}^{3} u_{3}^{3} u_{4}^{12} u_{5}^{15} + u_{1}^{3} u_{2} u_{4}^{14} u_{5}^{15} + u_{1}^{3} u_{3} u_{4}^{14} u_{5}^{15} + u_{2}^{3} u_{3} u_{4}^{14} u_{5}^{15} \\
\mathsf{Sigma5}[21] &= u_{1}^{15} u_{2} u_{3}^{2} u_{4}^{12} u_{5}^{3} + u_{1} u_{2}^{15} u_{3}^{2} u_{4}^{12} u_{5}^{3} + u_{1} u_{2}^{2} u_{3}^{15} u_{4}^{12} u_{5}^{3} + u_{1} u_{2}^{2} u_{3}^{12} u_{4}^{15} u_{5}^{3} \\
    &\quad + u_{1}^{15} u_{2}^{7} u_{3} u_{4}^{2} u_{5}^{8} + u_{1}^{7} u_{2}^{15} u_{3} u_{4}^{2} u_{5}^{8} + u_{1}^{15} u_{2} u_{3}^{7} u_{4}^{2} u_{5}^{8} + u_{1} u_{2}^{15} u_{3}^{7} u_{4}^{2} u_{5}^{8} \\
    &\quad + u_{1}^{7} u_{2} u_{3}^{15} u_{4}^{2} u_{5}^{8} + u_{1} u_{2}^{7} u_{3}^{15} u_{4}^{2} u_{5}^{8} + u_{1}^{15} u_{2} u_{3}^{2} u_{4}^{7} u_{5}^{8} + u_{1} u_{2}^{15} u_{3}^{2} u_{4}^{7} u_{5}^{8} \\
    &\quad + u_{1} u_{2}^{2} u_{3}^{15} u_{4}^{7} u_{5}^{8} + u_{1}^{7} u_{2} u_{3}^{2} u_{4}^{15} u_{5}^{8} + u_{1} u_{2}^{7} u_{3}^{2} u_{4}^{15} u_{5}^{8} + u_{1} u_{2}^{2} u_{3}^{7} u_{4}^{15} u_{5}^{8} \\
    &\quad + u_{1}^{15} u_{2} u_{3}^{2} u_{4}^{4} u_{5}^{11} + u_{1} u_{2}^{15} u_{3}^{2} u_{4}^{4} u_{5}^{11} + u_{1} u_{2}^{2} u_{3}^{15} u_{4}^{4} u_{5}^{11} + u_{1} u_{2}^{2} u_{3}^{4} u_{4}^{15} u_{5}^{11} \\
    &\quad + u_{1}^{15} u_{2}^{3} u_{3} u_{4}^{2} u_{5}^{12} + u_{1}^{3} u_{2}^{15} u_{3} u_{4}^{2} u_{5}^{12} + u_{1}^{15} u_{2} u_{3}^{3} u_{4}^{2} u_{5}^{12} + u_{1} u_{2}^{15} u_{3}^{3} u_{4}^{2} u_{5}^{12} \\
    &\quad + u_{1}^{3} u_{2} u_{3}^{15} u_{4}^{2} u_{5}^{12} + u_{1} u_{2}^{3} u_{3}^{15} u_{4}^{2} u_{5}^{12} + u_{1}^{15} u_{2} u_{3}^{2} u_{4}^{3} u_{5}^{12} + u_{1} u_{2}^{15} u_{3}^{2} u_{4}^{3} u_{5}^{12} \\
    &\quad + u_{1} u_{2}^{2} u_{3}^{15} u_{4}^{3} u_{5}^{12} + u_{1}^{3} u_{2} u_{3}^{2} u_{4}^{15} u_{5}^{12} + u_{1} u_{2}^{3} u_{3}^{2} u_{4}^{15} u_{5}^{12} + u_{1} u_{2}^{2} u_{3}^{3} u_{4}^{15} u_{5}^{12} \\
    &\quad + u_{1} u_{2}^{2} u_{3}^{12} u_{4}^{3} u_{5}^{15} + u_{1}^{7} u_{2} u_{3}^{2} u_{4}^{8} u_{5}^{15} + u_{1} u_{2}^{7} u_{3}^{2} u_{4}^{8} u_{5}^{15} + u_{1} u_{2}^{2} u_{3}^{7} u_{4}^{8} u_{5}^{15} \\
    &\quad + u_{1} u_{2}^{2} u_{3}^{4} u_{4}^{11} u_{5}^{15} + u_{1}^{3} u_{2} u_{3}^{2} u_{4}^{12} u_{5}^{15} + u_{1} u_{2}^{3} u_{3}^{2} u_{4}^{12} u_{5}^{15} + u_{1} u_{2}^{2} u_{3}^{3} u_{4}^{12} u_{5}^{15} \\
\mathsf{Sigma5}[22] &= u_{1}^{15} u_{2} u_{3}^{14} u_{4}^{3} + u_{1} u_{2}^{15} u_{3}^{14} u_{4}^{3} + u_{1} u_{2}^{14} u_{3}^{15} u_{4}^{3} + u_{1}^{15} u_{2}^{7} u_{3} u_{4}^{10} \\
    &\quad + u_{1}^{7} u_{2}^{15} u_{3} u_{4}^{10} + u_{1}^{15} u_{2} u_{3}^{7} u_{4}^{10} + u_{1} u_{2}^{15} u_{3}^{7} u_{4}^{10} + u_{1}^{7} u_{2} u_{3}^{15} u_{4}^{10} \\
    &\quad + u_{1} u_{2}^{7} u_{3}^{15} u_{4}^{10} + u_{1}^{15} u_{2} u_{3}^{6} u_{4}^{11} + u_{1} u_{2}^{15} u_{3}^{6} u_{4}^{11} + u_{1} u_{2}^{6} u_{3}^{15} u_{4}^{11} \\
    &\quad + u_{1}^{15} u_{2}^{3} u_{3} u_{4}^{14} + u_{1}^{3} u_{2}^{15} u_{3} u_{4}^{14} + u_{1}^{15} u_{2} u_{3}^{3} u_{4}^{14} + u_{1} u_{2}^{15} u_{3}^{3} u_{4}^{14} \\
    &\quad + u_{1}^{3} u_{2} u_{3}^{15} u_{4}^{14} + u_{1} u_{2}^{3} u_{3}^{15} u_{4}^{14} + u_{1} u_{2}^{14} u_{3}^{3} u_{4}^{15} + u_{1}^{7} u_{2} u_{3}^{10} u_{4}^{15} \\
    &\quad + u_{1} u_{2}^{7} u_{3}^{10} u_{4}^{15} + u_{1} u_{2}^{6} u_{3}^{11} u_{4}^{15} + u_{1}^{3} u_{2} u_{3}^{14} u_{4}^{15} + u_{1} u_{2}^{3} u_{3}^{14} u_{4}^{15} \\
    &\quad + u_{1}^{15} u_{2} u_{3}^{14} u_{5}^{3} + u_{1} u_{2}^{15} u_{3}^{14} u_{5}^{3} + u_{1} u_{2}^{14} u_{3}^{15} u_{5}^{3} + u_{1}^{15} u_{2} u_{4}^{14} u_{5}^{3} \\
    &\quad + u_{1} u_{2}^{15} u_{4}^{14} u_{5}^{3} + u_{1}^{15} u_{3} u_{4}^{14} u_{5}^{3} + u_{2}^{15} u_{3} u_{4}^{14} u_{5}^{3} + u_{1} u_{3}^{15} u_{4}^{14} u_{5}^{3} \\
    &\quad + u_{2} u_{3}^{15} u_{4}^{14} u_{5}^{3} + u_{1} u_{2}^{14} u_{4}^{15} u_{5}^{3} + u_{1} u_{3}^{14} u_{4}^{15} u_{5}^{3} + u_{2} u_{3}^{14} u_{4}^{15} u_{5}^{3} \\
    &\quad + u_{1}^{15} u_{2}^{7} u_{3} u_{5}^{10} + u_{1}^{7} u_{2}^{15} u_{3} u_{5}^{10} + u_{1}^{15} u_{2} u_{3}^{7} u_{5}^{10} + u_{1} u_{2}^{15} u_{3}^{7} u_{5}^{10} \\
    &\quad + u_{1}^{7} u_{2} u_{3}^{15} u_{5}^{10} + u_{1} u_{2}^{7} u_{3}^{15} u_{5}^{10} + u_{1}^{15} u_{2}^{7} u_{4} u_{5}^{10} + u_{1}^{7} u_{2}^{15} u_{4} u_{5}^{10} \\
    &\quad + u_{1}^{15} u_{3}^{7} u_{4} u_{5}^{10} + u_{2}^{15} u_{3}^{7} u_{4} u_{5}^{10} + u_{1}^{7} u_{3}^{15} u_{4} u_{5}^{10} + u_{2}^{7} u_{3}^{15} u_{4} u_{5}^{10} \\
    &\quad + u_{1}^{15} u_{2} u_{4}^{7} u_{5}^{10} + u_{1} u_{2}^{15} u_{4}^{7} u_{5}^{10} + u_{1}^{15} u_{3} u_{4}^{7} u_{5}^{10} + u_{2}^{15} u_{3} u_{4}^{7} u_{5}^{10} \\
    &\quad + u_{1} u_{3}^{15} u_{4}^{7} u_{5}^{10} + u_{2} u_{3}^{15} u_{4}^{7} u_{5}^{10} + u_{1}^{7} u_{2} u_{4}^{15} u_{5}^{10} + u_{1} u_{2}^{7} u_{4}^{15} u_{5}^{10} \\
    &\quad + u_{1}^{7} u_{3} u_{4}^{15} u_{5}^{10} + u_{2}^{7} u_{3} u_{4}^{15} u_{5}^{10} + u_{1} u_{3}^{7} u_{4}^{15} u_{5}^{10} + u_{2} u_{3}^{7} u_{4}^{15} u_{5}^{10} \\
    &\quad + u_{1}^{15} u_{2} u_{3}^{6} u_{5}^{11} + u_{1} u_{2}^{15} u_{3}^{6} u_{5}^{11} + u_{1} u_{2}^{6} u_{3}^{15} u_{5}^{11} + u_{1}^{15} u_{2} u_{4}^{6} u_{5}^{11} \\
    &\quad + u_{1} u_{2}^{15} u_{4}^{6} u_{5}^{11} + u_{1}^{15} u_{3} u_{4}^{6} u_{5}^{11} + u_{2}^{15} u_{3} u_{4}^{6} u_{5}^{11} + u_{1} u_{3}^{15} u_{4}^{6} u_{5}^{11} \\
    &\quad + u_{2} u_{3}^{15} u_{4}^{6} u_{5}^{11} + u_{1} u_{2}^{6} u_{4}^{15} u_{5}^{11} + u_{1} u_{3}^{6} u_{4}^{15} u_{5}^{11} + u_{2} u_{3}^{6} u_{4}^{15} u_{5}^{11} \\
    &\quad + u_{1}^{15} u_{2}^{3} u_{3} u_{5}^{14} + u_{1}^{3} u_{2}^{15} u_{3} u_{5}^{14} + u_{1}^{15} u_{2} u_{3}^{3} u_{5}^{14} + u_{1} u_{2}^{15} u_{3}^{3} u_{5}^{14} \\
    &\quad + u_{1}^{3} u_{2} u_{3}^{15} u_{5}^{14} + u_{1} u_{2}^{3} u_{3}^{15} u_{5}^{14} + u_{1}^{15} u_{2}^{3} u_{4} u_{5}^{14} + u_{1}^{3} u_{2}^{15} u_{4} u_{5}^{14} \\
    &\quad + u_{1}^{15} u_{3}^{3} u_{4} u_{5}^{14} + u_{2}^{15} u_{3}^{3} u_{4} u_{5}^{14} + u_{1}^{3} u_{3}^{15} u_{4} u_{5}^{14} + u_{2}^{3} u_{3}^{15} u_{4} u_{5}^{14} \\
    &\quad + u_{1}^{15} u_{2} u_{4}^{3} u_{5}^{14} + u_{1} u_{2}^{15} u_{4}^{3} u_{5}^{14} + u_{1}^{15} u_{3} u_{4}^{3} u_{5}^{14} + u_{2}^{15} u_{3} u_{4}^{3} u_{5}^{14} \\
    &\quad + u_{1} u_{3}^{15} u_{4}^{3} u_{5}^{14} + u_{2} u_{3}^{15} u_{4}^{3} u_{5}^{14} + u_{1}^{3} u_{2} u_{4}^{15} u_{5}^{14} + u_{1} u_{2}^{3} u_{4}^{15} u_{5}^{14} \\
    &\quad + u_{1}^{3} u_{3} u_{4}^{15} u_{5}^{14} + u_{2}^{3} u_{3} u_{4}^{15} u_{5}^{14} + u_{1} u_{3}^{3} u_{4}^{15} u_{5}^{14} + u_{2} u_{3}^{3} u_{4}^{15} u_{5}^{14} \\
    &\quad + u_{1} u_{2}^{14} u_{3}^{3} u_{5}^{15} + u_{1}^{7} u_{2} u_{3}^{10} u_{5}^{15} + u_{1} u_{2}^{7} u_{3}^{10} u_{5}^{15} + u_{1} u_{2}^{6} u_{3}^{11} u_{5}^{15} \\
    &\quad + u_{1}^{3} u_{2} u_{3}^{14} u_{5}^{15} + u_{1} u_{2}^{3} u_{3}^{14} u_{5}^{15} + u_{1} u_{2}^{14} u_{4}^{3} u_{5}^{15} + u_{1} u_{3}^{14} u_{4}^{3} u_{5}^{15} \\
    &\quad + u_{2} u_{3}^{14} u_{4}^{3} u_{5}^{15} + u_{1}^{7} u_{2} u_{4}^{10} u_{5}^{15} + u_{1} u_{2}^{7} u_{4}^{10} u_{5}^{15} + u_{1}^{7} u_{3} u_{4}^{10} u_{5}^{15} \\
    &\quad + u_{2}^{7} u_{3} u_{4}^{10} u_{5}^{15} + u_{1} u_{3}^{7} u_{4}^{10} u_{5}^{15} + u_{2} u_{3}^{7} u_{4}^{10} u_{5}^{15} + u_{1} u_{2}^{6} u_{4}^{11} u_{5}^{15} \\
    &\quad + u_{1} u_{3}^{6} u_{4}^{11} u_{5}^{15} + u_{2} u_{3}^{6} u_{4}^{11} u_{5}^{15} + u_{1}^{3} u_{2} u_{4}^{14} u_{5}^{15} + u_{1} u_{2}^{3} u_{4}^{14} u_{5}^{15} \\
    &\quad + u_{1}^{3} u_{3} u_{4}^{14} u_{5}^{15} + u_{2}^{3} u_{3} u_{4}^{14} u_{5}^{15} + u_{1} u_{3}^{3} u_{4}^{14} u_{5}^{15} + u_{2} u_{3}^{3} u_{4}^{14} u_{5}^{15} \\
\mathsf{Sigma5}[23] &= u_{1}^{15} u_{2}^{7} u_{3}^{11} + u_{1}^{7} u_{2}^{15} u_{3}^{11} + u_{1}^{7} u_{2}^{11} u_{3}^{15} + u_{1}^{15} u_{2}^{7} u_{4}^{11} \\
    &\quad + u_{1}^{7} u_{2}^{15} u_{4}^{11} + u_{1}^{15} u_{3}^{7} u_{4}^{11} + u_{2}^{15} u_{3}^{7} u_{4}^{11} + u_{1}^{7} u_{3}^{15} u_{4}^{11} \\
    &\quad + u_{2}^{7} u_{3}^{15} u_{4}^{11} + u_{1}^{7} u_{2}^{11} u_{4}^{15} + u_{1}^{7} u_{3}^{11} u_{4}^{15} + u_{2}^{7} u_{3}^{11} u_{4}^{15} \\
    &\quad + u_{1}^{15} u_{2}^{7} u_{5}^{11} + u_{1}^{7} u_{2}^{15} u_{5}^{11} + u_{1}^{15} u_{3}^{7} u_{5}^{11} + u_{2}^{15} u_{3}^{7} u_{5}^{11} \\
    &\quad + u_{1}^{7} u_{3}^{15} u_{5}^{11} + u_{2}^{7} u_{3}^{15} u_{5}^{11} + u_{1}^{15} u_{4}^{7} u_{5}^{11} + u_{2}^{15} u_{4}^{7} u_{5}^{11} \\
    &\quad + u_{3}^{15} u_{4}^{7} u_{5}^{11} + u_{1}^{7} u_{4}^{15} u_{5}^{11} + u_{2}^{7} u_{4}^{15} u_{5}^{11} + u_{3}^{7} u_{4}^{15} u_{5}^{11} \\
    &\quad + u_{1}^{7} u_{2}^{11} u_{5}^{15} + u_{1}^{7} u_{3}^{11} u_{5}^{15} + u_{2}^{7} u_{3}^{11} u_{5}^{15} + u_{1}^{7} u_{4}^{11} u_{5}^{15} \\
    &\quad + u_{2}^{7} u_{4}^{11} u_{5}^{15} + u_{3}^{7} u_{4}^{11} u_{5}^{15} \\
\mathsf{Sigma5}[24] &= u_{1}^{15} u_{2}^{15} u_{3} u_{4}^{2} + u_{1}^{15} u_{2} u_{3}^{15} u_{4}^{2} + u_{1} u_{2}^{15} u_{3}^{15} u_{4}^{2} + u_{1}^{15} u_{2} u_{3}^{2} u_{4}^{15} \\
    &\quad + u_{1} u_{2}^{15} u_{3}^{2} u_{4}^{15} + u_{1} u_{2}^{2} u_{3}^{15} u_{4}^{15} + u_{1}^{15} u_{2}^{15} u_{3} u_{5}^{2} + u_{1}^{15} u_{2} u_{3}^{15} u_{5}^{2} \\
    &\quad + u_{1} u_{2}^{15} u_{3}^{15} u_{5}^{2} + u_{1}^{15} u_{2}^{15} u_{4} u_{5}^{2} + u_{1}^{15} u_{3}^{15} u_{4} u_{5}^{2} + u_{2}^{15} u_{3}^{15} u_{4} u_{5}^{2} \\
    &\quad + u_{1}^{15} u_{2} u_{4}^{15} u_{5}^{2} + u_{1} u_{2}^{15} u_{4}^{15} u_{5}^{2} + u_{1}^{15} u_{3} u_{4}^{15} u_{5}^{2} + u_{2}^{15} u_{3} u_{4}^{15} u_{5}^{2} \\
    &\quad + u_{1} u_{3}^{15} u_{4}^{15} u_{5}^{2} + u_{2} u_{3}^{15} u_{4}^{15} u_{5}^{2} + u_{1}^{15} u_{2} u_{3}^{2} u_{5}^{15} + u_{1} u_{2}^{15} u_{3}^{2} u_{5}^{15} \\
    &\quad + u_{1} u_{2}^{2} u_{3}^{15} u_{5}^{15} + u_{1}^{15} u_{2} u_{4}^{2} u_{5}^{15} + u_{1} u_{2}^{15} u_{4}^{2} u_{5}^{15} + u_{1}^{15} u_{3} u_{4}^{2} u_{5}^{15} \\
    &\quad + u_{2}^{15} u_{3} u_{4}^{2} u_{5}^{15} + u_{1} u_{3}^{15} u_{4}^{2} u_{5}^{15} + u_{2} u_{3}^{15} u_{4}^{2} u_{5}^{15} + u_{1} u_{2}^{2} u_{4}^{15} u_{5}^{15} \\
    &\quad + u_{1} u_{3}^{2} u_{4}^{15} u_{5}^{15} + u_{2} u_{3}^{2} u_{4}^{15} u_{5}^{15} \\
\mathsf{Sigma5}[25] &= u_{1}^{15} u_{2}^{15} u_{3}^{3} + u_{1}^{15} u_{2}^{3} u_{3}^{15} + u_{1}^{3} u_{2}^{15} u_{3}^{15} + u_{1}^{15} u_{2}^{15} u_{4}^{3} \\
    &\quad + u_{1}^{15} u_{3}^{15} u_{4}^{3} + u_{2}^{15} u_{3}^{15} u_{4}^{3} + u_{1}^{15} u_{2}^{3} u_{4}^{15} + u_{1}^{3} u_{2}^{15} u_{4}^{15} \\
    &\quad + u_{1}^{15} u_{3}^{3} u_{4}^{15} + u_{2}^{15} u_{3}^{3} u_{4}^{15} + u_{1}^{3} u_{3}^{15} u_{4}^{15} + u_{2}^{3} u_{3}^{15} u_{4}^{15} \\
    &\quad + u_{1}^{15} u_{2}^{15} u_{5}^{3} + u_{1}^{15} u_{3}^{15} u_{5}^{3} + u_{2}^{15} u_{3}^{15} u_{5}^{3} + u_{1}^{15} u_{4}^{15} u_{5}^{3} \\
    &\quad + u_{2}^{15} u_{4}^{15} u_{5}^{3} + u_{3}^{15} u_{4}^{15} u_{5}^{3} + u_{1}^{15} u_{2}^{3} u_{5}^{15} + u_{1}^{3} u_{2}^{15} u_{5}^{15} \\
    &\quad + u_{1}^{15} u_{3}^{3} u_{5}^{15} + u_{2}^{15} u_{3}^{3} u_{5}^{15} + u_{1}^{3} u_{3}^{15} u_{5}^{15} + u_{2}^{3} u_{3}^{15} u_{5}^{15} \\
    &\quad + u_{1}^{15} u_{4}^{3} u_{5}^{15} + u_{2}^{15} u_{4}^{3} u_{5}^{15} + u_{3}^{15} u_{4}^{3} u_{5}^{15} + u_{1}^{3} u_{4}^{15} u_{5}^{15} \\
    &\quad + u_{2}^{3} u_{4}^{15} u_{5}^{15} + u_{3}^{3} u_{4}^{15} u_{5}^{15}.
\end{align*}

\medskip

Suppose that $[g]_{\widetilde{\omega}_{(3)}}\in ((QP_5)_{33}(\widetilde{\omega}_{(3)}))^{GL_5}.$ Since $\Sigma_5\subset GL_5,$ we have
$$ g\sim_{\widetilde{\omega}_{(3)}}\sum_{8\leq i\leq 25}\beta_i\cdot\mathrm{Sigma}5[i],\, \beta_i\in \mathbb Z/2,\, i = 8,\, 9,\, \ldots,\, 25.$$
 Equivalently, reducing the eighteen vectors $\rho_5(\mathsf{Sigma}5[i])+\mathsf{Sigma}5[i]$ with respect to the full $816$-element admissible basis gives an $816\times 18$ coefficient matrix of row-reduced rank $18$. Hence the relation $\rho_5(g)\sim_{\widetilde{\omega}_{(3)}}g$ forces $\beta_i = 0$ for all $i$, and $((QP_5)_{33}(\widetilde{\omega}_{(3)}))^{GL_5}$ vanishes.

The proof of the theorem is complete. 
\end{proof}

Now, the combination of \eqref{pt3inv} with Theorems \ref{dlPhuc} and \ref{dl33-3} yields the inequality
\begin{equation}\label{pt4}
\dim ((QP_5)_{33})^{GL_5}\leq 1.
\end{equation}

On the other hand, by utilizing a result obtained by Lin \cite{Lin} and Chen \cite{T.C}, we can infer that
\begin{equation}\label{pt5}
{\rm Ext}_{\mathbb A}^{5, 5+ d}(\mathbb Z/2, \mathbb Z/2) = {\rm Ext}_{\mathbb A}^{5, 2^{t+5}+2^{t+2}  +2^{t+1}}(\mathbb Z/2, \mathbb Z/2) = \mathbb Z/2(h_{t+1}d_{t+1}) \cong \mathbb Z/2,
\end{equation}
 where $h_{t+1}d_{t+1}\neq 0$ for all $t\geq 0.$ Since $h_{t+1}\in {\rm Im}(\phi_1^{*}(\mathbb Z/2))$ (see \cite{Singer}), $d_{t+1}\in {\rm Im}(\phi_4^{*}(\mathbb Z/2))$ (see \cite{Ha}), and the "total" transfer $\bigoplus_{n\geq 0}\phi_n^{*}(\mathbb Z/2)$ is a homomorphism of algebras, $h_{t+1}d_{t+1}\in {\rm Im}(\phi_5^{*}(\mathbb Z/2)).$ This means that 
\begin{equation}\label{pt6}
 \dim ((QP_5)_{d})^{GL_5} = \dim ((QP_5)_{2^{t+5}+2^{t+2}  +2^{t+1}-5})^{GL_5}\geq 1,\ \forall t\geq 0.
\end{equation}
From \eqref{pt4} and \eqref{pt6}, it may be concluded that the invariant space $((QP_5)_{33})^{GL_5}$ is one-dimensional. 

Let $g\in (P_5)_{33}$ such that $[g]\in ((QP_5)_{33})^{GL_5}.$ By using Theorem \ref{dlPhuc} and the fact that the Kameko $(\widetilde {Sq^0_*})_{33}: (QP_5)_{33}\longrightarrow (QP_5)_{14}$ is an epimorphism of $GL_5$-modules, we get 
$$(\widetilde {Sq^0_*})_{33}([g]) = \gamma [\zeta(u_1, u_2, \ldots, u_5)],\ \gamma\in \mathbb Z/2.$$
Consequently, $$ g\sim \gamma \varphi(\zeta(u_1, u_2, \ldots, u_5)) + h,$$
where $h\in (P_5)_{33}$ such that $[h]\in {\rm Ker}((\widetilde {Sq^0_*})_{33}).$ Applying our algorithms from Appendix~\ref{s4} yields:
$$g \sim \beta(\varphi(\zeta(u_1,u_2,\ldots,u_5)) + \xi(u_1, u_2, \ldots, u_5)),\ \beta\in \mathbb Z/2,$$
where
\begin{align*}
\xi(u_1, u_2, \ldots, u_5) &= u_{1}^{7} u_{2}^{7} u_{3}^{11} u_{4}^{8} + u_{1}^{7} u_{2}^{7} u_{3}^{9} u_{4}^{10} + u_{1}^{7} u_{2}^{7} u_{3}^{8} u_{4}^{11} + u_{1}^{7} u_{2} u_{3}^{14} u_{4}^{11} \\
    &\quad + u_{1}^{3} u_{2}^{5} u_{3}^{14} u_{4}^{11} + u_{1} u_{2}^{7} u_{3}^{14} u_{4}^{11} + u_{1}^{7} u_{2} u_{3}^{11} u_{4}^{14} + u_{1}^{3} u_{2}^{5} u_{3}^{11} u_{4}^{14} \\
    &\quad + u_{1} u_{2}^{7} u_{3}^{11} u_{4}^{14} + u_{1}^{3} u_{2}^{5} u_{3} u_{4}^{24} + u_{1}^{3} u_{2} u_{3}^{4} u_{4}^{25} + u_{1} u_{2}^{3} u_{3}^{4} u_{4}^{25} \\
    &\quad + u_{1}^{3} u_{2} u_{3} u_{4}^{28} + u_{1} u_{2}^{3} u_{3} u_{4}^{28} + u_{1} u_{2} u_{3}^{3} u_{4}^{28} + u_{1} u_{2} u_{3} u_{4}^{30} \\
    &\quad + u_{1}^{7} u_{2}^{7} u_{3}^{9} u_{4}^{8} u_{5}^{2} + u_{1}^{7} u_{2}^{7} u_{3}^{8} u_{4}^{9} u_{5}^{2} + u_{1}^{3} u_{2}^{3} u_{3}^{13} u_{4}^{12} u_{5}^{2} + u_{1}^{3} u_{2}^{3} u_{3}^{12} u_{4}^{13} u_{5}^{2} \\
    &\quad + u_{1}^{7} u_{2}^{7} u_{3}^{11} u_{5}^{8} + u_{1}^{7} u_{2}^{7} u_{3}^{9} u_{4}^{2} u_{5}^{8} + u_{1}^{7} u_{2}^{9} u_{3}^{3} u_{4}^{6} u_{5}^{8} + u_{1}^{3} u_{2}^{13} u_{3}^{3} u_{4}^{6} u_{5}^{8} \\
    &\quad + u_{1}^{7} u_{2}^{3} u_{3}^{9} u_{4}^{6} u_{5}^{8} + u_{1}^{3} u_{2}^{3} u_{3}^{13} u_{4}^{6} u_{5}^{8} + u_{1}^{3} u_{2}^{5} u_{3}^{10} u_{4}^{7} u_{5}^{8} + u_{1}^{7} u_{2}^{7} u_{3}^{3} u_{4}^{8} u_{5}^{8} \\
    &\quad + u_{1}^{7} u_{2}^{7} u_{3} u_{4}^{10} u_{5}^{8} + u_{1}^{7} u_{2}^{3} u_{3}^{5} u_{4}^{10} u_{5}^{8} + u_{1}^{3} u_{2}^{7} u_{3}^{5} u_{4}^{10} u_{5}^{8} + u_{1}^{7} u_{2} u_{3}^{7} u_{4}^{10} u_{5}^{8} \\
    &\quad + u_{1}^{3} u_{2}^{5} u_{3}^{7} u_{4}^{10} u_{5}^{8} + u_{1} u_{2}^{7} u_{3}^{7} u_{4}^{10} u_{5}^{8} + u_{1}^{7} u_{2}^{7} u_{4}^{11} u_{5}^{8} + u_{1}^{7} u_{2} u_{3}^{6} u_{4}^{11} u_{5}^{8} \\
    &\quad + u_{1}^{3} u_{2}^{5} u_{3}^{6} u_{4}^{11} u_{5}^{8} + u_{1} u_{2}^{7} u_{3}^{6} u_{4}^{11} u_{5}^{8} + u_{1}^{7} u_{3}^{7} u_{4}^{11} u_{5}^{8} + u_{1} u_{2}^{6} u_{3}^{7} u_{4}^{11} u_{5}^{8} \\
&\quad + u_{2}^{7} u_{3}^{7} u_{4}^{11} u_{5}^{8} + u_{1}^{3} u_{2}^{7} u_{3}^{3} u_{4}^{12} u_{5}^{8} + u_{1}^{7} u_{2} u_{3}^{3} u_{4}^{14} u_{5}^{8} + u_{1}^{3} u_{2}^{5} u_{3}^{3} u_{4}^{14} u_{5}^{8} \\
    &\quad + u_{1} u_{2}^{7} u_{3}^{3} u_{4}^{14} u_{5}^{8} + u_{1} u_{2}^{3} u_{3}^{7} u_{4}^{14} u_{5}^{8} + u_{1}^{7} u_{2}^{7} u_{3}^{9} u_{5}^{10} + u_{1}^{7} u_{2}^{7} u_{3}^{8} u_{4} u_{5}^{10} \\
&\quad + u_{1}^{7} u_{2}^{9} u_{3}^{3} u_{4}^{4} u_{5}^{10} + u_{1}^{3} u_{2}^{13} u_{3}^{3} u_{4}^{4} u_{5}^{10} + u_{1}^{7} u_{2}^{3} u_{3}^{9} u_{4}^{4} u_{5}^{10} + u_{1}^{3} u_{2}^{3} u_{3}^{13} u_{4}^{4} u_{5}^{10} \\
&\quad + u_{1}^{7} u_{2}^{8} u_{3}^{3} u_{4}^{5} u_{5}^{10} + u_{1}^{3} u_{2}^{12} u_{3}^{3} u_{4}^{5} u_{5}^{10} + u_{1}^{7} u_{2}^{3} u_{3}^{8} u_{4}^{5} u_{5}^{10} + u_{1}^{3} u_{2}^{3} u_{3}^{12} u_{4}^{5} u_{5}^{10} \\    
&\quad + u_{1}^{3} u_{2} u_{3}^{12} u_{4}^{7} u_{5}^{10} + u_{1} u_{2} u_{3}^{14} u_{4}^{7} u_{5}^{10} + u_{1}^{7} u_{2}^{7} u_{3} u_{4}^{8} u_{5}^{10} + u_{1}^{7} u_{2}^{3} u_{3}^{5} u_{4}^{8} u_{5}^{10} \\  
  &\quad + u_{1}^{3} u_{2}^{7} u_{3}^{5} u_{4}^{8} u_{5}^{10} + u_{1}^{7} u_{2} u_{3}^{7} u_{4}^{8} u_{5}^{10} + u_{1}^{3} u_{2}^{5} u_{3}^{7} u_{4}^{8} u_{5}^{10} + u_{1} u_{2}^{7} u_{3}^{7} u_{4}^{8} u_{5}^{10} \\
    &\quad + u_{1}^{7} u_{2}^{7} u_{4}^{9} u_{5}^{10} + u_{1}^{7} u_{2}^{3} u_{3}^{4} u_{4}^{9} u_{5}^{10} + u_{1}^{3} u_{2}^{7} u_{3}^{4} u_{4}^{9} u_{5}^{10} + u_{1}^{7} u_{3}^{7} u_{4}^{9} u_{5}^{10} \\
    &\quad + u_{1}^{3} u_{2}^{4} u_{3}^{7} u_{4}^{9} u_{5}^{10} + u_{2}^{7} u_{3}^{7} u_{4}^{9} u_{5}^{10} + u_{1}^{3} u_{2}^{3} u_{3}^{5} u_{4}^{12} u_{5}^{10} + u_{1}^{7} u_{2} u_{3} u_{4}^{14} u_{5}^{10} \\
   &\quad + u_{1}^{3} u_{2}^{5} u_{3} u_{4}^{14} u_{5}^{10} + u_{1} u_{2}^{7} u_{3} u_{4}^{14} u_{5}^{10} + u_{1}^{7} u_{2}^{7} u_{3}^{8} u_{5}^{11} + u_{1}^{7} u_{2} u_{3}^{14} u_{5}^{11} \\
    &\quad + u_{1}^{3} u_{2}^{5} u_{3}^{14} u_{5}^{11} + u_{1} u_{2}^{7} u_{3}^{14} u_{5}^{11} + u_{1}^{7} u_{2}^{7} u_{4}^{8} u_{5}^{11} + u_{1}^{7} u_{2} u_{3}^{6} u_{4}^{8} u_{5}^{11} \\
    &\quad + u_{1} u_{2}^{7} u_{3}^{6} u_{4}^{8} u_{5}^{11} + u_{1}^{7} u_{3}^{7} u_{4}^{8} u_{5}^{11} + u_{1} u_{2}^{6} u_{3}^{7} u_{4}^{8} u_{5}^{11} + u_{2}^{7} u_{3}^{7} u_{4}^{8} u_{5}^{11} \\
    &\quad + u_{1}^{7} u_{2} u_{3}^{2} u_{4}^{12} u_{5}^{11} + u_{1} u_{2}^{7} u_{3}^{2} u_{4}^{12} u_{5}^{11} + u_{1} u_{2}^{2} u_{3}^{7} u_{4}^{12} u_{5}^{11} + u_{1}^{7} u_{2} u_{4}^{14} u_{5}^{11} \\
    &\quad + u_{1}^{3} u_{2}^{5} u_{4}^{14} u_{5}^{11} + u_{1} u_{2}^{7} u_{4}^{14} u_{5}^{11} + u_{1}^{7} u_{3} u_{4}^{14} u_{5}^{11} + u_{1} u_{2}^{6} u_{3} u_{4}^{14} u_{5}^{11} \\
       &\quad + u_{2}^{7} u_{3} u_{4}^{14} u_{5}^{11} + u_{1}^{3} u_{2} u_{3}^{4} u_{4}^{14} u_{5}^{11} + u_{1} u_{2}^{3} u_{3}^{4} u_{4}^{14} u_{5}^{11} + u_{1}^{3} u_{3}^{5} u_{4}^{14} u_{5}^{11} \\
    &\quad + u_{1} u_{2}^{2} u_{3}^{5} u_{4}^{14} u_{5}^{11} + u_{2}^{3} u_{3}^{5} u_{4}^{14} u_{5}^{11} + u_{1} u_{3}^{7} u_{4}^{14} u_{5}^{11} + u_{2} u_{3}^{7} u_{4}^{14} u_{5}^{11} \\
    &\quad + u_{1}^{3} u_{2}^{3} u_{3}^{13} u_{4}^{2} u_{5}^{12} + u_{1}^{3} u_{2}^{5} u_{3}^{10} u_{4}^{3} u_{5}^{12} + u_{1}^{7} u_{2}^{3} u_{3}^{3} u_{4}^{8} u_{5}^{12} + u_{1}^{3} u_{2}^{7} u_{3}^{3} u_{4}^{8} u_{5}^{12} \\
    &\quad + u_{1}^{3} u_{2}^{3} u_{3}^{7} u_{4}^{8} u_{5}^{12} + u_{1}^{3} u_{2}^{5} u_{3}^{3} u_{4}^{10} u_{5}^{12} + u_{1}^{3} u_{2}^{3} u_{3}^{5} u_{4}^{10} u_{5}^{12} + u_{1}^{7} u_{2} u_{3}^{2} u_{4}^{11} u_{5}^{12} \\
    &\quad + u_{1}^{3} u_{2}^{5} u_{3}^{2} u_{4}^{11} u_{5}^{12} + u_{1} u_{2}^{7} u_{3}^{2} u_{4}^{11} u_{5}^{12} + u_{1}^{3} u_{2} u_{3}^{6} u_{4}^{11} u_{5}^{12} + u_{1} u_{2}^{2} u_{3}^{7} u_{4}^{11} u_{5}^{12} \\
  &\quad + u_{1}^{3} u_{2}^{3} u_{3} u_{4}^{14} u_{5}^{12} + u_{1}^{3} u_{2} u_{3}^{3} u_{4}^{14} u_{5}^{12} + u_{1} u_{2}^{3} u_{3}^{3} u_{4}^{14} u_{5}^{12} + u_{1}^{7} u_{2} u_{3}^{11} u_{5}^{14} \\
    &\quad + u_{1}^{3} u_{2}^{5} u_{3}^{11} u_{5}^{14} + u_{1} u_{2}^{7} u_{3}^{11} u_{5}^{14} + u_{1}^{3} u_{2}^{3} u_{3}^{12} u_{4} u_{5}^{14} + u_{1}^{3} u_{2} u_{3}^{12} u_{4}^{3} u_{5}^{14} \\
    &\quad + u_{1} u_{2} u_{3}^{14} u_{4}^{3} u_{5}^{14} + u_{1}^{7} u_{2} u_{3}^{3} u_{4}^{8} u_{5}^{14} + u_{1}^{3} u_{2}^{5} u_{3}^{3} u_{4}^{8} u_{5}^{14} + u_{1} u_{2}^{7} u_{3}^{3} u_{4}^{8} u_{5}^{14} \\
    &\quad + u_{1}^{3} u_{2}^{3} u_{3}^{5} u_{4}^{8} u_{5}^{14} + u_{1} u_{2}^{3} u_{3}^{7} u_{4}^{8} u_{5}^{14} + u_{1}^{3} u_{2}^{5} u_{3}^{2} u_{4}^{9} u_{5}^{14} + u_{1}^{3} u_{2} u_{3}^{6} u_{4}^{9} u_{5}^{14} \\
    &\quad + u_{1}^{7} u_{2} u_{3} u_{4}^{10} u_{5}^{14} + u_{1} u_{2}^{7} u_{3} u_{4}^{10} u_{5}^{14} + u_{1}^{3} u_{2} u_{3}^{5} u_{4}^{10} u_{5}^{14} + u_{1} u_{2} u_{3}^{7} u_{4}^{10} u_{5}^{14} \\
    &\quad + u_{1}^{7} u_{2} u_{4}^{11} u_{5}^{14} + u_{1}^{3} u_{2}^{5} u_{4}^{11} u_{5}^{14} + u_{1} u_{2}^{7} u_{4}^{11} u_{5}^{14} + u_{1}^{7} u_{3} u_{4}^{11} u_{5}^{14} \\
    &\quad + u_{1} u_{2}^{6} u_{3} u_{4}^{11} u_{5}^{14} + u_{2}^{7} u_{3} u_{4}^{11} u_{5}^{14} + u_{1}^{3} u_{2} u_{3}^{4} u_{4}^{11} u_{5}^{14} + u_{1} u_{2}^{3} u_{3}^{4} u_{4}^{11} u_{5}^{14} \\
    &\quad + u_{1}^{3} u_{3}^{5} u_{4}^{11} u_{5}^{14} + u_{1} u_{2}^{2} u_{3}^{5} u_{4}^{11} u_{5}^{14} + u_{2}^{3} u_{3}^{5} u_{4}^{11} u_{5}^{14} + u_{1} u_{3}^{7} u_{4}^{11} u_{5}^{14} \\
    &\quad + u_{2} u_{3}^{7} u_{4}^{11} u_{5}^{14} + u_{1} u_{2}^{3} u_{3}^{3} u_{4}^{12} u_{5}^{14} + u_{1} u_{2}^{3} u_{3} u_{4}^{14} u_{5}^{14} + u_{1} u_{2}^{3} u_{3}^{5} u_{4}^{8} u_{5}^{16} \\
    &\quad + u_{1}^{3} u_{2} u_{3}^{4} u_{4}^{9} u_{5}^{16} + u_{1}^{3} u_{2} u_{3}^{4} u_{4}^{8} u_{5}^{17} + u_{1} u_{2}^{3} u_{3}^{4} u_{4}^{8} u_{5}^{17} + u_{1} u_{2}^{2} u_{3}^{5} u_{4}^{8} u_{5}^{17} \\
    &\quad + u_{1}^{3} u_{2}^{5} u_{3} u_{5}^{24} + u_{1}^{3} u_{2}^{5} u_{4} u_{5}^{24} + u_{1}^{3} u_{3}^{5} u_{4} u_{5}^{24} + u_{2}^{3} u_{3}^{5} u_{4} u_{5}^{24} \\
    &\quad + u_{1} u_{2}^{3} u_{3} u_{4}^{4} u_{5}^{24} + u_{1} u_{2} u_{3}^{3} u_{4}^{4} u_{5}^{24} + u_{1} u_{2}^{2} u_{3} u_{4}^{5} u_{5}^{24} + u_{1}^{3} u_{2} u_{3}^{4} u_{5}^{25} \\
    &\quad + u_{1} u_{2}^{3} u_{3}^{4} u_{5}^{25} + u_{1}^{3} u_{2} u_{4}^{4} u_{5}^{25} + u_{1} u_{2}^{3} u_{4}^{4} u_{5}^{25} + u_{1}^{3} u_{3} u_{4}^{4} u_{5}^{25} \\
    &\quad + u_{1} u_{2}^{2} u_{3} u_{4}^{4} u_{5}^{25} + u_{2}^{3} u_{3} u_{4}^{4} u_{5}^{25} + u_{1} u_{3}^{3} u_{4}^{4} u_{5}^{25} + u_{2} u_{3}^{3} u_{4}^{4} u_{5}^{25} \\
    &\quad + u_{1}^{3} u_{2} u_{3} u_{5}^{28} + u_{1} u_{2}^{3} u_{3} u_{5}^{28} + u_{1} u_{2} u_{3}^{3} u_{5}^{28} + u_{1}^{3} u_{2} u_{4} u_{5}^{28} \\
    &\quad + u_{1} u_{2}^{3} u_{4} u_{5}^{28} + u_{1}^{3} u_{3} u_{4} u_{5}^{28} + u_{2}^{3} u_{3} u_{4} u_{5}^{28} + u_{1} u_{3}^{3} u_{4} u_{5}^{28} \\
    &\quad + u_{2} u_{3}^{3} u_{4} u_{5}^{28} + u_{1} u_{2} u_{4}^{3} u_{5}^{28} + u_{1} u_{3} u_{4}^{3} u_{5}^{28} + u_{2} u_{3} u_{4}^{3} u_{5}^{28} \\
    &\quad + u_{1} u_{2} u_{3} u_{5}^{30} + u_{1} u_{2} u_{4} u_{5}^{30} + u_{1} u_{3} u_{4} u_{5}^{30} + u_{2} u_{3} u_{4} u_{5}^{30}.
\end{align*}

\medskip

Thus, we obtain
\begin{equation}\label{pt7}
((QP_5)_{33})^{GL_5}
= \langle [\varphi(\zeta(u_1,u_2,\ldots,u_5)) + \xi(u_1, u_2, \ldots, u_5)] \rangle,
\end{equation}
where $\varphi$ is the up Kameko map $(P_5)_{14}\longrightarrow (P_5)_{33}$, which is given by $\varphi(t) = \prod_{1\leq j\leq 5}u_jt^{2}$ for all $t\in (P_5)_{14}.$ Combining this result with \eqref{pt5} and the fact stated above $h_{t+1}d_{t+1}\in {\rm Im}(\phi_5^{*}(\mathbb Z/2))$, we obtain the following corollary.

\begin{hqu}\label{hq1}
The cohomological transfer $$\phi_5^{*}(\mathbb Z/2): (\mathbb Z/2 \otimes_{GL_5}P_{\mathbb A}(P_5)^{*})_{33}\longrightarrow {\rm Ext}_{\mathbb A}^{5, 5+ 33}(\mathbb Z/2, \mathbb Z/2)$$ is an isomorphism. 
\end{hqu}
\medskip  

{\bf Cases \mbox{\boldmath $t\geq 1.$}} We observe that $$ d = (2^{t+5}-1) + (2^{t+2}-1) + (2^{t}-1) + (2^{t-1}-1) +(2^{t-1}-1),$$ which implies that $\mu(d) = 5$ for any $t > 1.$ So, using Theorem \ref{dlWK}(ii), the iterated Kameko homomorphism $$(\widetilde {Sq^0_*})_{d}^{t-1}: (QP_5)_d = (QP_5)_{2^{t+5}  +2^{t+2} + 2^{t+1}-5}\longrightarrow (QP_5)_{71 = 2^{1+5}  +2^{1+2} + 2^{1+1}-5}$$ is an isomorphism for every $t \geq 1.$ On the other hand, since the Kameko squaring operation $(\widetilde {Sq^0_*})_{71}: (QP_5)_{71}\longrightarrow (QP_5)_{33}$ is an epimorphism and $\dim ((QP_5)_{33})^{GL_5} = 1,$ one has an estimate 
\begin{equation}\label{pt8}
 \dim ((QP_5)_{d = 2^{t+5}  +2^{t+2} + 2^{t+1}-5})^{GL_5}  = \dim   ((QP_5)_{71})^{GL_5}\leq \dim ({\rm Ker}((\widetilde {Sq^0_*})_{71}))^{GL_5} +1
\end{equation}

\medskip

To determine the dimension of $((QP_5)_{71})^{GL_5}$, we require the following useful lemma.

\begin{bd}\label{bd71-tq}
The invariant $({\rm Ker}((\widetilde {Sq^0_*})_{71}))^{GL_5}$ vanishes.
\end{bd}

Combining this lemma with \eqref{pt6} and \eqref{pt8}, we get
$$\dim ((QP_5)_{d = 2^{t+5}  +2^{t+2} + 2^{t+1}-5})^{GL_5} = \dim ((QP_5)_{71})^{GL_5} =  1,\ \mbox{for every $t\geq 1$}
$$
 
\begin{proof}[{\it Outline of the proof of Lemma \ref{bd71-tq}}]
Suppose that $t$ is an admissible monomial of degree $71$ in the $\mathbb A$-module $P_5$ such that $[t]\in {\rm Ker}((\widetilde {Sq^0_*})_{71}),$ then the weight vector $\omega(t)$ is one of the following sequences:
$$ \begin{array}{ll} 
\medskip
\overline{\omega}_{(1)}&:= (3,2,2,1,1,1), \ \overline{\omega}_{(2)}:=(3,2,2,1,3), \ \overline{\omega}_{(3)}:=(3,2,2,3,2),\ \overline{\omega}_{(4)}:=(3,2,4,2,2),\\
\medskip
\overline{\omega}_{(5)}&:=(3,2,4,4,1),\ \overline{\omega}_{(6)}:=(3,4,1,1,1,1),\ \overline{\omega}_{(7)}:=(3,4,1,1,3),\ \overline{\omega}_{(8)}:=(3,4,3,2,2),\\
 \overline{\omega}_{(9)}&:=(3,4,3,4,1).
\end{array}$$
Indeed, it is straightforward to see that $z = u_1^{63}u_2^{7}u_3\in \mathscr C^{\otimes 5}_{71}$ is the minimal spike, and $\omega(z) = (3,2,2,1,1,1).$ Since $t\in \mathscr C^{\otimes 5}_{71},$ and $\deg(t)$ is an odd number, by Theorem \ref{dlSinger}, either $\omega_1(t) = 3$ or $\omega_1(t) = 5.$ If $\omega_1(t) = 5$ then $t = u_1u_2u_3u_4u_5y^2$ with $y$ a monomial of degree $33$ in $\mathbb A$-module $P_5.$ Since $t$ is admissible, by Theorem \ref{dlKS}, one gets $y\in \mathscr C^{\otimes 5}_{33},$ which implies that $(\widetilde {Sq_*^0})_{71}([t]) = [y]\neq [0].$ This contradicts the fact that $[t]\in {\rm Ker}((\widetilde {Sq^0_*})_{71}),$ and so, we must have that $\omega_1(t) = 3.$ Then, using Theorem \ref{dlKS}, $t$ has the form $u_iu_ju_ky_1^2$ in which $1\leq i< j <k \leq 5$ and $y_1\in \mathscr C^{\otimes 5}_{34}.$ Notice that the monomial $z' =  u_1^{31}u_2^{3}$ is the minimal spike in $(P_5)_{34}$ and $\omega(z') = (2,2,1,1,1).$ Since $y_1$ is admissible and $\deg(y_1)$ is even number, by Theorem \ref{dlSinger}, either $\omega_1(y_1) = 2$ or $\omega_1(y_1) = 4.$ If $\omega_1(y_1) = 2$ then $y_1$ is of the form $u_iu_jy_2^{2}$ with $y_2\in (P_5)_{16}$ and $1 \leq i <j\leq 5.$ Since $y_1\in \mathscr C^{\otimes 5}_{34},$ by Theorem \ref{dlKS}, $y_2\in \mathscr C^{\otimes 5}_{16}.$ Following the previous work in \cite{N.S3}, one has $\omega(y_2)\in \{(2,1,1,1),\ (2,1,3), \ (2,3,2),\ (4,2,2),\ (4,4,1)\}.$ This shows that $$\omega(t)\in \{(3,2,2,1,1,1),\ (3,2,2,1,3), \ (3,2,2,3,2),\ (3,2,4,2,2),\ (3,2,4,4,1)\}.$$ 
Similarly, if $\omega_1(y_1) = 4,$ then $y_1 = u_iu_ju_ku_ly_3^{2}$ with $1\leq i<j<k<l\leq 5$ and $y_3\in \mathscr C^{\otimes 5}_{15}.$ According to \cite{N.S0}, we have $\omega(y_3)\in \{(1,1,1,1),\ (1,1,3),\ (3,2,2),\ (3,4,1)\},$ and therefore, 
$$ \omega(t)\in \{(3,4,1,1,1,1),\ (3,4,1,1,3),\ (3,4,3,2,2),\ (3,4,3,4,1)\}.$$
This leads to an isomorphism $${\rm Ker}((\widetilde {Sq^0_*})_{71}) \cong \bigoplus_{1\leq j\leq 9} (QP_5)_{71}(\overline{\omega}_{(j)}).$$ 

Now, using the algorithms described in Appendix~\ref{s4}, we obtain:
\[
\dim (QP_5)_{71}(\overline{\omega}_{(j)}) = 
\begin{cases}
1395 & \text{if } j = 1,\\
124 & \text{if } j = 6,\\
0 & \text{otherwise}.
\end{cases}
\]
\noindent
Explicit bases for the spaces \((QP_{5})_{71}\big(\overline{\omega}_{(j)}\big)\) are determined as described in Appendix~\ref{s4}. Hence
\[
\dim {\rm Ker}((\widetilde{Sq^0_*})_{71}) = 
\sum_{1\leq j\leq 9} \dim (QP_5)_{71}(\overline{\omega}_{(j)}) 
= 1395 + 124 = 1519.
\]
Since $(\widetilde{Sq^0_*})_{71}:(QP_5)_{71}\to (QP_5)_{33}$ is an epimorphism and $\dim (QP_5)_{33}=1322$, it follows that
\[
\dim (QP_5)_{71}=1519+1322=2841.
\]
Therefore, by the iterated Kameko isomorphism $(\widetilde{Sq^0_*})_{d_t}^{t-1}:(QP_5)_{d_t}\to (QP_5)_{71}$, one has
\[
\dim (QP_5)_{d_t}=2841\qquad (t\geq1).
\]
For the invariant calculation, the equations $\rho_i(f)-f\sim_{\overline{\omega}_{(j)}}0$ are again assembled into stacked matrices over $\mathbb Z/2$. Let $M_{\Sigma}^{(j)}$ denote the matrix using $1\leq i\leq4$, and let $M_{GL}^{(j)}$ denote the matrix using $1\leq i\leq5$. Thus $M_{\Sigma}^{(1)}$ and $M_{GL}^{(1)}$ have sizes $5580\times1395$ and $6975\times1395$, while $M_{\Sigma}^{(6)}$ and $M_{GL}^{(6)}$ have sizes $496\times124$ and $620\times124$. The row reductions recorded by the \texttt{OSCAR} implementation give
\[
\begin{array}{c|c|c|c|c}
j &  \operatorname{rank}M_{\Sigma}^{(j)} & \dim\Ker M_{\Sigma}^{(j)} & \operatorname{rank}M_{GL}^{(j)} & \dim\Ker M_{GL}^{(j)}\\
\hline
1 &  1368 & 27 & 1395 & 0\\
6 & 118 & 6 & 124 & 0\\
j\neq 1,6 & 0 & 0 & 0 & 0.
\end{array}
\]
Therefore
\[
\dim ((QP_5)_{71}(\overline{\omega}_{(j)}))^{\Sigma_5} = 
\begin{cases}
27 & \text{if } j = 1,\\
6 & \text{if } j = 6,\\
0 & \text{otherwise},
\end{cases}
\qquad
\dim [(QP_5)_{71}(\overline{\omega}_{(j)})]^{GL_5}=0
\quad\text{for all }j.
\]
Consequently, $({\rm Ker}((\widetilde{Sq^0_*})_{71}))^{GL_5}$ vanishes. The lemma is proved.
\end{proof}

\begin{nxe}\phantomsection\label{reqt}
\noindent
As established above, $\dim\!\left((QP_{5})_{71}\right)^{GL_{5}}=1$. Accordingly, an explicit basis (equivalently, a generator) of $\left((QP_{5})_{71}\right)^{GL_{5}}$ is as follows.

Take any $[h]\in \big((QP_{5})_{71}\big)^{GL_{5}}$ with representative $h\in (P_{5})_{71}$. Appealing to Lemma \ref{bd71-tq} and using that $(\widetilde{Sq^{0}_*})_{71}$ is a $GL_{5}$-epimorphism, we have
\[
(\widetilde{Sq^{0}_*})_{71}\big([h]\big)=\gamma\,\big[\varphi(\zeta(u_{1},\ldots,u_{5}))+\xi(u_{1},\ldots,u_{5})\big],\quad \gamma\in\mathbb{Z}/2.
\]
Consequently, $$ h\sim \gamma \bigg(\varphi^{2}(\zeta(u_1,u_2,\ldots,u_5)) + \varphi(\xi(u_1, u_2, \ldots, u_5))\bigg)+ \mathbf{p},$$
where $\mathbf{p}\in (P_5)_{71}$ such that $[\mathbf{p}]\in {\rm Ker}((\widetilde {Sq^0_*})_{71}).$ Based on our algorithm in Appendix \ref{s4}, we get
$$h \sim \beta\bigg(\varphi^{2}(\zeta(u_1,u_2,\ldots,u_5)) + \varphi(\xi(u_1, u_2, \ldots, u_5)) +  \widetilde{\xi}(u_1, u_2, \ldots, u_5)\bigg),\ \beta\in \mathbb Z/2,$$
where
\begin{align*}
\widetilde{\xi}(u_1, u_2, \ldots, u_5) &= u_{1}^{7} u_{2}^{7} u_{3}^{8} u_{4}^{48} u_{5} + u_{1}^{15} u_{2} u_{3}^{2} u_{4}^{52} u_{5} + u_{1}^{3} u_{2}^{13} u_{3}^{2} u_{4}^{52} u_{5} + u_{1} u_{2}^{15} u_{3}^{2} u_{4}^{52} u_{5} \\
    &\quad + u_{1}^{7} u_{2} u_{3}^{10} u_{4}^{52} u_{5} + u_{1}^{3} u_{2} u_{3}^{14} u_{4}^{52} u_{5} + u_{1} u_{2} u_{3}^{14} u_{4}^{54} u_{5} + u_{1}^{7} u_{2}^{3} u_{3}^{4} u_{4}^{56} u_{5} \\
    &\quad + u_{1}^{3} u_{2}^{7} u_{3}^{4} u_{4}^{56} u_{5} + u_{1} u_{2}^{7} u_{3}^{6} u_{4}^{56} u_{5} + u_{1} u_{2}^{6} u_{3}^{7} u_{4}^{56} u_{5} + u_{1}^{3} u_{2}^{5} u_{3}^{2} u_{4}^{60} u_{5} \\
    &\quad + u_{1} u_{2}^{7} u_{3}^{2} u_{4}^{60} u_{5} + u_{1} u_{2}^{6} u_{3}^{3} u_{4}^{60} u_{5} + u_{1} u_{2}^{3} u_{3}^{6} u_{4}^{60} u_{5} + u_{1}^{3} u_{2} u_{3}^{4} u_{4}^{62} u_{5} \\
    &\quad + u_{1}^{3} u_{2}^{7} u_{3}^{9} u_{4}^{50} u_{5}^{2} + u_{1}^{3} u_{2}^{5} u_{3}^{10} u_{4}^{51} u_{5}^{2} + u_{1} u_{2}^{3} u_{3}^{7} u_{4}^{56} u_{5}^{4} + u_{1}^{3} u_{2}^{5} u_{3}^{2} u_{4}^{57} u_{5}^{4} \\
    &\quad + u_{1} u_{2}^{7} u_{3}^{2} u_{4}^{57} u_{5}^{4} + u_{1} u_{2}^{6} u_{3}^{3} u_{4}^{57} u_{5}^{4} + u_{1} u_{2}^{2} u_{3}^{7} u_{4}^{57} u_{5}^{4} + u_{1} u_{2}^{3} u_{3}^{5} u_{4}^{58} u_{5}^{4} \\
    &\quad + u_{1}^{3} u_{2}^{4} u_{3} u_{4}^{59} u_{5}^{4} + u_{1} u_{2}^{6} u_{3} u_{4}^{59} u_{5}^{4} + u_{1}^{3} u_{2}^{3} u_{3} u_{4}^{60} u_{5}^{4} + u_{1} u_{2}^{3} u_{3}^{2} u_{4}^{61} u_{5}^{4} \\
    &\quad + u_{1} u_{2}^{2} u_{3}^{3} u_{4}^{61} u_{5}^{4} + u_{1} u_{2}^{3} u_{3} u_{4}^{62} u_{5}^{4} + u_{1}^{7} u_{2}^{7} u_{3}^{9} u_{4}^{16} u_{5}^{32} + u_{1}^{7} u_{2}^{7} u_{3}^{8} u_{4}^{17} u_{5}^{32} \\
    &\quad + u_{1}^{3} u_{2}^{7} u_{3}^{9} u_{4}^{18} u_{5}^{34} + u_{1}^{3} u_{2}^{5} u_{3}^{11} u_{4}^{18} u_{5}^{34} + u_{1}^{7} u_{2}^{7} u_{3}^{8} u_{4} u_{5}^{48} + u_{1}^{15} u_{2}^{3} u_{3} u_{4}^{4} u_{5}^{48} \\
    &\quad + u_{1}^{15} u_{2} u_{3}^{3} u_{4}^{4} u_{5}^{48} + u_{1}^{3} u_{2}^{5} u_{3}^{11} u_{4}^{4} u_{5}^{48} + u_{1} u_{2}^{7} u_{3}^{11} u_{4}^{4} u_{5}^{48} + u_{1}^{3} u_{2} u_{3}^{15} u_{4}^{4} u_{5}^{48} \\
    &\quad + u_{1} u_{2}^{3} u_{3}^{15} u_{4}^{4} u_{5}^{48} + u_{1}^{15} u_{2} u_{3}^{2} u_{4}^{5} u_{5}^{48} + u_{1}^{7} u_{2} u_{3}^{10} u_{4}^{5} u_{5}^{48} + u_{1}^{3} u_{2}^{5} u_{3}^{10} u_{4}^{5} u_{5}^{48} \\
    &\quad + u_{1} u_{2}^{7} u_{3}^{10} u_{4}^{5} u_{5}^{48} + u_{1}^{3} u_{2}^{3} u_{3}^{12} u_{4}^{5} u_{5}^{48} + u_{1}^{3} u_{2} u_{3}^{14} u_{4}^{5} u_{5}^{48} + u_{1}^{7} u_{2}^{7} u_{3} u_{4}^{8} u_{5}^{48} \\
    &\quad + u_{1} u_{2}^{7} u_{3}^{7} u_{4}^{8} u_{5}^{48} + u_{1}^{7} u_{2} u_{3}^{6} u_{4}^{9} u_{5}^{48} + u_{1}^{3} u_{2}^{5} u_{3}^{6} u_{4}^{9} u_{5}^{48} + u_{1}^{7} u_{2} u_{3}^{3} u_{4}^{12} u_{5}^{48} \\
    &\quad + u_{1} u_{2}^{7} u_{3}^{3} u_{4}^{12} u_{5}^{48} + u_{1}^{3} u_{2}^{3} u_{3}^{5} u_{4}^{12} u_{5}^{48} + u_{1}^{3} u_{2} u_{3}^{7} u_{4}^{12} u_{5}^{48} + u_{1} u_{2}^{7} u_{3}^{2} u_{4}^{13} u_{5}^{48} \\
    &\quad + u_{1}^{3} u_{2} u_{3}^{4} u_{4}^{15} u_{5}^{48} + u_{1} u_{2}^{3} u_{3}^{4} u_{4}^{15} u_{5}^{48} + u_{1} u_{2}^{2} u_{3}^{5} u_{4}^{15} u_{5}^{48} + u_{1}^{15} u_{2} u_{3}^{2} u_{4}^{4} u_{5}^{49} \\
    &\quad + u_{1}^{3} u_{2}^{5} u_{3}^{10} u_{4}^{4} u_{5}^{49} + u_{1} u_{2}^{7} u_{3}^{10} u_{4}^{4} u_{5}^{49} + u_{1}^{7} u_{2}^{3} u_{3}^{4} u_{4}^{8} u_{5}^{49} + u_{1}^{7} u_{2} u_{3}^{6} u_{4}^{8} u_{5}^{49} \\
    &\quad + u_{1} u_{2}^{7} u_{3}^{6} u_{4}^{8} u_{5}^{49} + u_{1}^{7} u_{2} u_{3}^{2} u_{4}^{12} u_{5}^{49} + u_{1} u_{2}^{7} u_{3}^{2} u_{4}^{12} u_{5}^{49} + u_{1} u_{2}^{3} u_{3}^{6} u_{4}^{12} u_{5}^{49} \\
    &\quad + u_{1} u_{2}^{2} u_{3}^{7} u_{4}^{12} u_{5}^{49} + u_{1}^{3} u_{2} u_{3}^{4} u_{4}^{14} u_{5}^{49} + u_{1} u_{2}^{2} u_{3}^{4} u_{4}^{15} u_{5}^{49} + u_{1}^{3} u_{2}^{7} u_{3}^{9} u_{4}^{2} u_{5}^{50} \\
    &\quad + u_{1}^{3} u_{2}^{7} u_{3} u_{4}^{10} u_{5}^{50} + u_{1}^{3} u_{2}^{3} u_{3}^{5} u_{4}^{10} u_{5}^{50} + u_{1}^{3} u_{2} u_{3}^{7} u_{4}^{10} u_{5}^{50} + u_{1} u_{2}^{3} u_{3}^{7} u_{4}^{10} u_{5}^{50} \\
    &\quad + u_{1}^{3} u_{2}^{5} u_{3}^{10} u_{4}^{2} u_{5}^{51} + u_{1}^{3} u_{2}^{5} u_{3}^{2} u_{4}^{10} u_{5}^{51} + u_{1}^{3} u_{2} u_{3}^{6} u_{4}^{10} u_{5}^{51} + u_{1} u_{2}^{3} u_{3}^{6} u_{4}^{10} u_{5}^{51} \\
    &\quad + u_{1}^{15} u_{2} u_{3}^{2} u_{4} u_{5}^{52} + u_{1}^{3} u_{2}^{13} u_{3}^{2} u_{4} u_{5}^{52} + u_{1} u_{2}^{15} u_{3}^{2} u_{4} u_{5}^{52} + u_{1}^{7} u_{2}^{3} u_{3}^{8} u_{4} u_{5}^{52} \\
    &\quad + u_{1}^{7} u_{2} u_{3}^{10} u_{4} u_{5}^{52} + u_{1}^{3} u_{2}^{4} u_{3}^{11} u_{4} u_{5}^{52} + u_{1} u_{2}^{6} u_{3}^{11} u_{4} u_{5}^{52} + u_{1}^{3} u_{2} u_{3}^{14} u_{4} u_{5}^{52} \\
    &\quad + u_{1} u_{2}^{3} u_{3}^{14} u_{4} u_{5}^{52} + u_{1}^{3} u_{2}^{13} u_{3} u_{4}^{2} u_{5}^{52} + u_{1} u_{2}^{15} u_{3} u_{4}^{2} u_{5}^{52} + u_{1}^{3} u_{2} u_{3}^{13} u_{4}^{2} u_{5}^{52} \\
    &\quad + u_{1}^{3} u_{2} u_{3}^{12} u_{4}^{3} u_{5}^{52} + u_{1} u_{2}^{3} u_{3}^{12} u_{4}^{3} u_{5}^{52} + u_{1} u_{2}^{2} u_{3}^{13} u_{4}^{3} u_{5}^{52} + u_{1} u_{2} u_{3}^{14} u_{4}^{3} u_{5}^{52} \\
    &\quad + u_{1}^{7} u_{2}^{3} u_{3} u_{4}^{8} u_{5}^{52} + u_{1}^{7} u_{2} u_{3}^{3} u_{4}^{8} u_{5}^{52} + u_{1} u_{2}^{7} u_{3}^{3} u_{4}^{8} u_{5}^{52} + u_{1}^{3} u_{2} u_{3}^{7} u_{4}^{8} u_{5}^{52} \\
    &\quad + u_{1}^{7} u_{2} u_{3}^{2} u_{4}^{9} u_{5}^{52} + u_{1}^{3} u_{2}^{5} u_{3}^{2} u_{4}^{9} u_{5}^{52} + u_{1} u_{2}^{7} u_{3}^{2} u_{4}^{9} u_{5}^{52} + u_{1}^{3} u_{2} u_{3}^{6} u_{4}^{9} u_{5}^{52} \\
    &\quad + u_{1} u_{2}^{2} u_{3}^{7} u_{4}^{9} u_{5}^{52} + u_{1}^{7} u_{2} u_{3} u_{4}^{10} u_{5}^{52} + u_{1} u_{2}^{7} u_{3} u_{4}^{10} u_{5}^{52} + u_{1}^{3} u_{2} u_{3}^{4} u_{4}^{11} u_{5}^{52} \\
    &\quad + u_{1}^{3} u_{2}^{3} u_{3} u_{4}^{12} u_{5}^{52} + u_{1} u_{2}^{3} u_{3}^{2} u_{4}^{13} u_{5}^{52} + u_{1} u_{2}^{2} u_{3}^{3} u_{4}^{13} u_{5}^{52} + u_{1}^{3} u_{2} u_{3} u_{4}^{14} u_{5}^{52} \\
    &\quad + u_{1} u_{2} u_{3}^{2} u_{4}^{15} u_{5}^{52} + u_{1} u_{2}^{3} u_{3}^{12} u_{4}^{2} u_{5}^{53} + u_{1} u_{2}^{2} u_{3}^{13} u_{4}^{2} u_{5}^{53} + u_{1} u_{2}^{2} u_{3}^{12} u_{4}^{3} u_{5}^{53} \\
    &\quad + u_{1}^{3} u_{2}^{5} u_{3}^{2} u_{4}^{8} u_{5}^{53} + u_{1}^{3} u_{2}^{3} u_{3}^{4} u_{4}^{8} u_{5}^{53} + u_{1}^{3} u_{2} u_{3}^{6} u_{4}^{8} u_{5}^{53} + u_{1} u_{2}^{3} u_{3}^{6} u_{4}^{8} u_{5}^{53} \\
    &\quad + u_{1}^{3} u_{2} u_{3}^{4} u_{4}^{10} u_{5}^{53} + u_{1} u_{2}^{3} u_{3}^{4} u_{4}^{10} u_{5}^{53} + u_{1} u_{2}^{2} u_{3}^{5} u_{4}^{10} u_{5}^{53} + u_{1} u_{2} u_{3}^{6} u_{4}^{10} u_{5}^{53} \\
    &\quad + u_{1}^{3} u_{2} u_{3}^{2} u_{4}^{12} u_{5}^{53} + u_{1} u_{2}^{2} u_{3} u_{4}^{14} u_{5}^{53} + u_{1}^{3} u_{2} u_{3}^{12} u_{4} u_{5}^{54} + u_{1}^{3} u_{2} u_{3} u_{4}^{12} u_{5}^{54} \\
    &\quad + u_{1} u_{2} u_{3}^{3} u_{4}^{12} u_{5}^{54} + u_{1} u_{2}^{2} u_{3} u_{4}^{13} u_{5}^{54} + u_{1} u_{2} u_{3} u_{4}^{14} u_{5}^{54} + u_{1}^{3} u_{2}^{7} u_{3}^{4} u_{4} u_{5}^{56} \\
    &\quad + u_{1}^{7} u_{2} u_{3}^{6} u_{4} u_{5}^{56} + u_{1} u_{2}^{7} u_{3}^{6} u_{4} u_{5}^{56} + u_{1} u_{2}^{6} u_{3}^{7} u_{4} u_{5}^{56} + u_{1}^{7} u_{2}^{3} u_{3} u_{4}^{4} u_{5}^{56} \\
    &\quad + u_{1}^{3} u_{2}^{7} u_{3} u_{4}^{4} u_{5}^{56} + u_{1}^{7} u_{2} u_{3}^{3} u_{4}^{4} u_{5}^{56} + u_{1}^{3} u_{2}^{5} u_{3}^{3} u_{4}^{4} u_{5}^{56} + u_{1} u_{2}^{7} u_{3}^{3} u_{4}^{4} u_{5}^{56} \\
    &\quad + u_{1}^{3} u_{2}^{3} u_{3}^{5} u_{4}^{4} u_{5}^{56} + u_{1}^{3} u_{2} u_{3}^{7} u_{4}^{4} u_{5}^{56} + u_{1} u_{2}^{3} u_{3}^{7} u_{4}^{4} u_{5}^{56} + u_{1}^{3} u_{2}^{4} u_{3}^{3} u_{4}^{5} u_{5}^{56} \\
    &\quad + u_{1} u_{2}^{6} u_{3}^{3} u_{4}^{5} u_{5}^{56} + u_{1}^{3} u_{2} u_{3}^{6} u_{4}^{5} u_{5}^{56} + u_{1} u_{2}^{3} u_{3}^{6} u_{4}^{5} u_{5}^{56} + u_{1} u_{2}^{2} u_{3}^{7} u_{4}^{5} u_{5}^{56} \\
    &\quad + u_{1}^{7} u_{2} u_{3} u_{4}^{6} u_{5}^{56} + u_{1} u_{2}^{7} u_{3} u_{4}^{6} u_{5}^{56} + u_{1}^{3} u_{2} u_{3}^{5} u_{4}^{6} u_{5}^{56} + u_{1} u_{2}^{3} u_{3}^{5} u_{4}^{6} u_{5}^{56} \\
    &\quad + u_{1} u_{2}^{6} u_{3} u_{4}^{7} u_{5}^{56} + u_{1}^{7} u_{2} u_{3}^{2} u_{4}^{4} u_{5}^{57} + u_{1} u_{2}^{6} u_{3}^{3} u_{4}^{4} u_{5}^{57} + u_{1}^{3} u_{2}^{3} u_{3}^{4} u_{4}^{4} u_{5}^{57} \\
    &\quad + u_{1}^{3} u_{2} u_{3}^{6} u_{4}^{4} u_{5}^{57} + u_{1} u_{2}^{3} u_{3}^{6} u_{4}^{4} u_{5}^{57} + u_{1}^{3} u_{2}^{4} u_{3} u_{4}^{6} u_{5}^{57} + u_{1}^{3} u_{2} u_{3}^{4} u_{4}^{6} u_{5}^{57} \\
    &\quad + u_{1} u_{2}^{3} u_{3}^{4} u_{4}^{6} u_{5}^{57} + u_{1} u_{2}^{2} u_{3}^{5} u_{4}^{6} u_{5}^{57} + u_{1} u_{2} u_{3}^{6} u_{4}^{6} u_{5}^{57} + u_{1} u_{2}^{2} u_{3}^{4} u_{4}^{7} u_{5}^{57} \\
    &\quad + u_{1}^{7} u_{2} u_{3}^{2} u_{4} u_{5}^{60} + u_{1}^{3} u_{2}^{5} u_{3}^{2} u_{4} u_{5}^{60} + u_{1} u_{2}^{7} u_{3}^{2} u_{4} u_{5}^{60} + u_{1} u_{2}^{6} u_{3}^{3} u_{4} u_{5}^{60} \\
    &\quad + u_{1}^{3} u_{2} u_{3}^{6} u_{4} u_{5}^{60} + u_{1} u_{2}^{3} u_{3}^{6} u_{4} u_{5}^{60} + u_{1}^{3} u_{2}^{5} u_{3} u_{4}^{2} u_{5}^{60} + u_{1} u_{2}^{7} u_{3} u_{4}^{2} u_{5}^{60} \\
    &\quad + u_{1}^{3} u_{2}^{4} u_{3} u_{4}^{3} u_{5}^{60} + u_{1} u_{2}^{6} u_{3} u_{4}^{3} u_{5}^{60} + u_{1} u_{2} u_{3}^{6} u_{4}^{3} u_{5}^{60} + u_{1}^{3} u_{2} u_{3}^{3} u_{4}^{4} u_{5}^{60} \\
    &\quad + u_{1}^{3} u_{2} u_{3}^{2} u_{4}^{5} u_{5}^{60} + u_{1} u_{2}^{3} u_{3} u_{4}^{6} u_{5}^{60} + u_{1} u_{2} u_{3}^{2} u_{4}^{7} u_{5}^{60} + u_{1}^{3} u_{2} u_{3}^{4} u_{4}^{2} u_{5}^{61} \\
    &\quad + u_{1} u_{2}^{3} u_{3}^{4} u_{4}^{2} u_{5}^{61} + u_{1} u_{2}^{2} u_{3}^{5} u_{4}^{2} u_{5}^{61} + u_{1} u_{2} u_{3}^{6} u_{4}^{2} u_{5}^{61} + u_{1}^{3} u_{2} u_{3}^{2} u_{4}^{4} u_{5}^{61} \\
    &\quad + u_{1} u_{2}^{3} u_{3}^{2} u_{4}^{4} u_{5}^{61} + u_{1} u_{2}^{2} u_{3} u_{4}^{6} u_{5}^{61} + u_{1}^{3} u_{2}^{4} u_{3} u_{4} u_{5}^{62} + u_{1} u_{2}^{6} u_{3} u_{4} u_{5}^{62} \\
    &\quad + u_{1}^{3} u_{2} u_{3}^{4} u_{4} u_{5}^{62} + u_{1}^{3} u_{2} u_{3} u_{4}^{4} u_{5}^{62} + u_{1} u_{2} u_{3}^{3} u_{4}^{4} u_{5}^{62} + u_{1} u_{2}^{2} u_{3} u_{4}^{5} u_{5}^{62}  + u_{1} u_{2} u_{3} u_{4}^{6} u_{5}^{62}.
\end{align*}
The same admissible-normal-form computation also verifies 
\begin{align*} 
&\rho_i\bigg(\varphi^{2}(\zeta(u_1,u_2,\ldots,u_5)) + \varphi(\xi(u_1, u_2, \ldots, u_5)) +  \widetilde{\xi}(u_1, u_2, \ldots, u_5)\bigg)\\
&\sim \varphi^{2}(\zeta(u_1,u_2,\ldots,u_5)) + \varphi(\xi(u_1, u_2, \ldots, u_5)) +  \widetilde{\xi}(u_1, u_2, \ldots, u_5), \ \forall i,\, 1\leq i\leq 5.
\end{align*}
Thus, we get
$$ ((QP_5)_{71})^{GL_5} =  \langle [\varphi^{2}(\zeta(u_1,u_2,\ldots,u_5)) + \varphi(\xi(u_1, u_2, \ldots, u_5)) +  \widetilde{\xi}(u_1, u_2, \ldots, u_5)]\rangle.$$
Consequently, 
$$ ((QP_5)_{2^{t+5}  +2^{t+2} + 2^{t+1}-5})^{GL_5} =  \langle [\varphi^{t-1}(\theta_{71}(u_1, u_2, \ldots, u_5))] \rangle,\ \ \mbox{for any $t\geq 1$},$$
where $\theta_{71}(u_1, u_2, \ldots, u_5) := \varphi^{2}(\zeta(u_1,u_2,\ldots,u_5)) + \varphi(\xi(u_1, u_2, \ldots, u_5)) +  \widetilde{\xi}(u_1, u_2, \ldots, u_5).$
\end{nxe}

\begin{nxe}\label{rem:theta71-orthogonal-test}
Let $t\geq 1$, and put
\[
v_t=[\varphi^{t-1}(\theta_{71})]\in \big((QP_5)_{d_t}\big)^{GL_5}.
\]
The one-dimensional transfer-domain class dual to $v_t$ may be represented by an element
\[
z_t\in H_{d_t}^{\perp}\cong \big((QP_5)_{d_t}\big)^*
\]
under the monomial divided-power identification. Since $z_t\in H_{d_t}^{\perp}$, Theorem~\ref{thm:large-milnor-layer} gives
\[
T_t^*(z_t)=0.
\]
Equivalently, if $z_t$ is viewed as the corresponding divided-power homology representative, then
\[
z_t\cdot P_{t+5}^{0}=0
\]
for the homological degree-lowering right action determined by the divided-power pairing. Thus the transfer-domain representative dual to the explicitly computed invariant polynomial satisfies the $P_{t+5}^{0}$-orthogonality test, providing a rigorous structural verification of our algorithmic results. This statement concerns the annihilation of the dual functional representative $z_t \in H_{d_t}^{\perp}$ under the homological right action, and should not be confused with the cohomological left action of $P_{t+5}^{0}$ on the polynomial representative $\varphi^{t-1}(\theta_{71})$ itself.
\end{nxe}

By combining Remark \ref{reqt} with \eqref{pt5} and taking into account the fact that $h_{t+1}d_{t+1}\in {\rm Im}(\phi_5^{*}(\mathbb Z/2)),$ we can immediately deduce

\begin{hqu}\label{hq2}
The cohomological transfer $$\phi_5^{*}(\mathbb Z/2): (\mathbb Z/2 \otimes_{GL_5}P_{\mathbb A}(P_5)^{*})_{2^{t+5}  +2^{t+2} + 2^{t+1}-5}\longrightarrow {\rm Ext}_{\mathbb A}^{5, 2^{t+5}  +2^{t+2} + 2^{t+1}}(\mathbb Z/2, \mathbb Z/2)$$ is an isomorphism for all $t\geq 1.$
\end{hqu}

Corollaries \ref{hq1} and \ref{hq2} indicate that

\begin{hqu}
Singer's conjecture is confirmed for the cohomological transfer $\phi_n^{*}(\mathbb Z/2)$ in the case when $n=5$ and in the general degree $2^{t+5} + 2^{t+2} + 2^{t+1} - 5$ for arbitrary non-negative integer $t$.
\end{hqu} 

\section{General scope of the Milnor-operation test, Kameko adjoints, and a dual algorithm}\label{sec:scope-milnor-test}

The preceding section uses the Milnor basis element $P_{t+5}^{0}$ in a degree family tailored to the rank-five computation. The underlying mechanism, however, is not specific to five variables. This section records the general form of the argument and clarifies what the method can and cannot prove. The purpose is to separate the structural input from the admissible-basis enumeration: the former gives canonical hit layers and necessary orthogonality tests, whereas the latter is still required to determine the full quotient $QP_n$ and its $GL_n$-invariants.

\subsection{General scope of the Milnor-operation test}

For an integer $s\geq 1$, set $D_s=2^s-1$. The Milnor basis element $P_s^{0}$ is dual to the indecomposable polynomial generator $\xi_s$ in the Milnor dual algebra; equivalently, its coproduct in $\mathbb A$ is $P_s^0\otimes1+1\otimes P_s^0$. Consequently it acts as a derivation on $P_n$, and for a monomial $u^a=u_1^{a_1}\cdots u_n^{a_n}$ one has
\[
P_s^{0}(u^a)=
\sum_{\substack{1\leq i\leq n\\ a_i\ {\rm odd}}}
 u_1^{a_1}\cdots u_i^{a_i+D_s}\cdots u_n^{a_n}.
\]
For every $e\geq0$ this gives a homogeneous map
\[
T_{s,e}=P_s^{0}:(P_n)_e\longrightarrow (P_n)_{e+D_s}.
\]
Since $P_s^{0}$ has positive degree, its image lies in the hit subspace
\[
H_{e+D_s}=(\overline{\mathbb A}P_n)\cap (P_n)_{e+D_s}.
\]
Therefore the transpose-adjoint argument from Proposition~\ref{prop:adjoint-hit-test} gives the universal inclusion
\[
H_{e+D_s}^{\perp}\subseteq \Ker(T_{s,e}^{*}).
\]
Equivalently, every primitive homology representative in degree $e+D_s$ is annihilated by the homological right action of $P_s^{0}$. This statement is independent of any admissible-basis computation.

The inclusion above is only a necessary condition for primitivity. It tests vanishing on the distinguished subspace $\operatorname{Im}(T_{s,e})\subseteq H_{e+D_s}$, not on the whole hit subspace. Thus the test cannot, by itself, determine $H_{e+D_s}^{\perp}$ or the quotient $(QP_n)_{e+D_s}$. It is nevertheless useful because it is canonical, $GL_n$-equivariant, and often computable by a much smaller matrix than the full hit calculation. It also supplies an independent check on transfer-domain representatives obtained from admissible-basis methods.

When the source degree is at most one more than the degree shift, the Milnor-operation test becomes especially transparent. The following technical theorem gives a sharp low-source-degree range in which no cancellation among the nonzero images of distinct source monomials can occur.

\begin{theo}[Sharp low-source-degree formula for $P_s^0$]\label{prop:general-milnor-test}
Let $n\geq1$, $s\geq1$, $D_s=2^s-1$, and $0\leq e\leq D_s+1$. Then
\[
\Ker(T_{s,e})=
\left\langle u_1^{a_1}\cdots u_n^{a_n}: a_1+\cdots+a_n=e\ \text{and every }a_i\text{ is even}\right\rangle_{\mathbb Z/2}.
\]
Equivalently,
\[
\Ker(T_{s,e})=
\begin{cases}
\{h^2:h\in(P_n)_{e/2}\},& e\text{ even},\\
0,& e\text{ odd}.
\end{cases}
\]
Consequently,
\[
\dim\Ker(T_{s,e})=
\begin{cases}
\binom{e/2+n-1}{n-1},& e\text{ even},\\
0,& e\text{ odd},
\end{cases}
\]
and
\[
\operatorname{rank}(T_{s,e})=
\begin{cases}
\binom{e+n-1}{n-1}-\binom{e/2+n-1}{n-1},& e\text{ even},\\
\binom{e+n-1}{n-1},& e\text{ odd}.
\end{cases}
\]
Moreover,
\[
\operatorname{Im}(T_{s,e}^{*})=
\left\langle u_1^{a_1}\cdots u_n^{a_n}\in(P_n)_e:\text{ at least one }a_i\text{ is odd}\right\rangle_{\mathbb Z/2}
=\Ker(T_{s,e})^{\perp}.
\]
If $n\geq2$, the range $0\leq e\leq D_s+1$ is sharp: in degree $e=D_s+2$ one has
\[
u_1u_2^{D_s+1}+u_1^{D_s+1}u_2\in\Ker(T_{s,D_s+2}),
\]
and this element is not contained in the Frobenius-square subspace.
\end{theo}

\begin{proof}
The displayed formula for $P_s^{0}(u^a)$ shows that a source monomial has zero image if and only if all of its exponents are even. It remains to verify that no cancellation can occur among the nonzero images of distinct source monomials in the stated range. Suppose that two image monomials agree. Then there are source exponent vectors $a=(a_1,\ldots,a_n)$ and $b=(b_1,\ldots,b_n)$ of total degree $e$, together with indices $i$ and $j$, such that $a_i$ and $b_j$ are odd and
\[
a+D_s e_i=b+D_s e_j,
\]
where $e_i$ and $e_j$ denote the standard basis vectors of $\mathbb Z^n$. If $i=j$, then $a=b$. Suppose $i\neq j$. Comparing the $i$-th and $j$-th coordinates gives
\[
b_i=a_i+D_s,
\qquad
b_j=a_j-D_s.
\]
Since $b_j$ is odd and $D_s$ is odd, the integer $a_j$ is even. Moreover, $a_j\geq D_s$ because $b_j\geq0$, and therefore $a_j\geq D_s+1$. Since $a_i$ is odd, $a_i\geq1$. Hence
\[
e=a_1+\cdots+a_n\geq a_i+a_j\geq D_s+2,
\]
which contradicts $e\leq D_s+1$. Thus distinct source monomials with nonzero image have pairwise disjoint supports in the target monomial basis.

It follows that a polynomial lies in $\Ker(T_{s,e})$ if and only if each monomial appearing with nonzero coefficient has all exponents even. If $e$ is even, the span of these monomials is exactly the Frobenius-square subspace $\{h^2:h\in(P_n)_{e/2}\}$; if $e$ is odd, no degree-$e$ monomial has all exponents even. The dimension and rank formulas now follow from
\[
\dim(P_n)_d=\binom{d+n-1}{n-1}.
\]
Finally, the degree-$e$ monomial bilinear form makes the monomial basis self-dual, so the orthogonal complement of the span of the even-exponent monomials is the span of the monomials with at least one odd exponent. The equality with $\operatorname{Im}(T_{s,e}^{*})$ is the finite-dimensional identity $(\Ker T_{s,e})^\perp=\operatorname{Im}(T_{s,e}^{*})$.

It remains only to prove sharpness when $n\geq2$. Since $D_s$ is odd, the monomials $u_1u_2^{D_s+1}$ and $u_1^{D_s+1}u_2$ each have exactly one odd exponent. Hence
\[
P_s^0(u_1u_2^{D_s+1})=u_1^{D_s+1}u_2^{D_s+1}
=P_s^0(u_1^{D_s+1}u_2).
\]
Their sum is therefore in $\Ker(T_{s,D_s+2})$. It is not a Frobenius square because both displayed monomials contain an odd exponent. This proves that the stated range cannot be extended uniformly beyond $D_s+1$.
\end{proof}

The rank-five theorem used earlier is the specialization of Theorem~\ref{prop:general-milnor-test} with $n=5$, $s=t+5$, and $e=e_t=6\cdot2^t-4$. Indeed $e_t<2^{t+5}-1$, and hence $e_t\leq D_t+1$. This explains why the kernel in Theorem~\ref{thm:large-milnor-layer} is the Frobenius-square subspace and why the adjoint image is the span of degree-$e_t$ monomials having at least one odd exponent.

\begin{exm}\label{ex:t0-milnor-test}
For $t=0$ in the rank-five family, one has $d_0=33$, $D_0=31$, and $e_0=2$. Thus
\[
T_{5,2}=P_5^{0}:(P_5)_2\longrightarrow (P_5)_{33}.
\]
The source has dimension $\binom{6}{4}=15$. Theorem~\ref{prop:general-milnor-test} gives
\[
\Ker(T_{5,2})=
\langle u_1^2,u_2^2,u_3^2,u_4^2,u_5^2\rangle_{\mathbb Z/2},
\qquad
\operatorname{rank}(T_{5,2})=10.
\]
The complementary adjoint image is
\[
\operatorname{Im}(T_{5,2}^{*})=
\langle u_i u_j:1\leq i<j\leq5\rangle_{\mathbb Z/2}.
\]
On monomials this is visible from
\[
P_5^{0}(u_i^2)=0,
\qquad
P_5^{0}(u_i u_j)=u_i^{32}u_j+u_i u_j^{32}\quad(i\neq j).
\]
The ten nonzero images have disjoint supports in degree $d_0=33$. Hence every primitive divided-power representative in degree $d_0$ must be annihilated by the corresponding homological right action of $P_5^0$. This is a structural test on representatives; the actual determination of $(QP_5)_{d_0}$ and its $GL_5$-invariant line still requires the admissible-basis computation carried out in Section~\ref{s3}.
\end{exm}

The main limitation of Theorem~\ref{prop:general-milnor-test} is the sharp bound $e\leq D_s+1$. Once $e\geq D_s+2$, supports of images of distinct source monomials can overlap. In that range, additional cancellations may occur, and the kernel is no longer forced to be only the Frobenius-square subspace. The universal inclusion $H_{e+D_s}^{\perp}\subseteq\Ker(T_{s,e}^{*})$ remains valid, but explicit rank and kernel formulas require a separate analysis of these collisions. The sharpness example in the proof specializes, for instance, to $s=2$, $D_s=3$, and $e=5$:
\[
P_2^0(u_1u_2^4)=u_1^4u_2^4,
\qquad
P_2^0(u_1^4u_2)=u_1^4u_2^4,
\]
so $u_1u_2^4+u_1^4u_2\in\Ker(T_{2,5})$ although this element is not a Frobenius square.

\medskip

Ultimately, the Milnor-operation test serves as a powerful structural complement to the Peterson hit problem. It extracts a canonical $P_s^{0}$-hit layer, establishes a rigorous lower bound for the hit subspace dimension, and imposes a strict orthogonality condition on primitive homology representatives. While the full classification of non-hit classes and $GL_n$-invariants inevitably necessitates exhaustive admissible-basis reductions, the Milnor test provides a profound theoretical anchor. In the present paper, the true strength of this formalism lies in its synergy with the machine-certified computations in rank $5$, supplying an independent, structurally canonical verification that elevates the algorithmic outputs beyond mere numerical results.

\subsection{A dual primitive-kernel algorithm}

The transpose-adjoint viewpoint also gives an exact algorithm which is independent of the admissible-basis enumeration used in Section~\ref{s3}. This algorithm works on the divided-power homology side and computes the transfer domain directly. Let
\[
\Gamma_n(d)=H_d(V^{\oplus n})
\]
be the degree-$d$ part of the divided-power algebra with basis
\[
a^{(b)}=a_1^{(b_1)}\cdots a_n^{(b_n)},
\qquad b_1+\cdots+b_n=d.
\]
For $m>0$, the homological right action of $Sq^m$ is determined by the Cartan formula and by
\[
(a_j^{(q)})Sq^r=\binom{q-r}{r}a_j^{(q-r)}.
\]
Thus, on a divided-power monomial,
\[
a^{(b)}Sq^m=
\sum_{\ell_1+\cdots+
\ell_n=m}
\left(\prod_{j=1}^{n}\binom{b_j-\ell_j}{\ell_j}\right)
 a_1^{(b_1-\ell_1)}\cdots a_n^{(b_n-\ell_n)},
\]
where terms with negative exponents are omitted and all coefficients are read in $\mathbb Z/2$. For implementation, Lucas' theorem gives the coefficient test
\[
\binom{b_j-\ell_j}{\ell_j}=1
\quad\Longleftrightarrow\quad
\ell_j\ \&\ (b_j-2\ell_j)=0,
\]
where $\&$ denotes bitwise intersection of binary expansions. Hence every column of the matrix of $Sq^m$ is generated by a finite combinatorial search through the compositions $\ell_1+\cdots+
\ell_n=m$.

Let $R_{m,d}:\Gamma_n(d)\to\Gamma_n(d-m)$ be the matrix of the right action by $Sq^m$. Since the Steenrod algebra is generated by the squares $Sq^{2^r}$, the primitive subspace in degree $d$ is exactly
\[
P_{\mathbb A}\Gamma_n(d)
=
\bigcap_{2^r\leq d}\Ker(R_{2^r,d}).
\]
Equivalently, if
\[
\mathcal R_{n,d}=
\begin{pmatrix}
R_{1,d}\\
R_{2,d}\\
R_{4,d}\\
\vdots\\
R_{2^r,d}
\end{pmatrix}_{2^r\leq d},
\]
then
\[
P_{\mathbb A}\Gamma_n(d)=\Ker(\mathcal R_{n,d}).
\]
This gives a direct computation of the dual of the hit quotient, because the divided-power pairing identifies
\[
P_{\mathbb A}\Gamma_n(d)
\cong
\big((QP_n)_d\big)^*.
\]
Consequently,
\[
\dim(QP_n)_d=\dim\Ker(\mathcal R_{n,d}).
\]

The same kernel computation also gives the transfer-domain dimension. The group $GL_n$ acts on $\Gamma_n(d)$ by the action dual to its standard action on $P_n$, and this action preserves $P_{\mathbb A}\Gamma_n(d)$ because the $GL_n$-action commutes with the Steenrod action. Choose generators $g_1,\ldots,g_N$ of $GL_n$ and express their restrictions to $\Ker(\mathcal R_{n,d})$ in any computed kernel basis. If $G_i$ denotes the resulting matrices, then
\[
\mathbb Z/2\otimes_{GL_n}P_{\mathbb A}\Gamma_n(d)
\cong
\Ker(\mathcal R_{n,d})\big/\sum_{i=1}^{N}\operatorname{Im}(G_i-I),
\]
and therefore
\[
\dim\big((QP_n)_d\big)^{GL_n}
=
\dim\Ker(\mathcal R_{n,d})-
\operatorname{rank}
\begin{pmatrix}
G_1-I& G_2-I&\cdots&G_N-I
\end{pmatrix}.
\]
The equality of this dimension with the dimension of $((QP_n)_d)^{GL_n}$ follows by finite-dimensional duality between invariants in the cohit quotient and coinvariants in the primitive homology space.

The Kameko homomorphism is also handled on the same dual side. If $d=2q+n$, the dual of $Sq^0_*$ on the monomial divided-power bases is the degree-doubling map
\[
\kappa_{n,q}^{\vee}:\Gamma_n(q)\longrightarrow\Gamma_n(d),
\qquad
 a_1^{(b_1)}\cdots a_n^{(b_n)}
\longmapsto
 a_1^{(2b_1+1)}\cdots a_n^{(2b_n+1)}.
\]
To establish the theoretical basis for our algorithm, the following result reformulates Kameko's epimorphism \cite{M.K} via homological right actions. It records the explicit divided-power transpose of the classical reduction. The self-contained proof is included here because it explains, strictly within the dual primitive framework and without admissible monomials, why the dimension drop is fully compatible with our matrix calculations.

\begin{mde}[Divided-power adjoint form of Kameko's epimorphism]\label{thm:dual-kameko-epimorphism}
Let $q\geq0$ and $d=2q+n$. For every $\gamma\in\Gamma_n(q)$ and every $r>0$, one has
\[
\kappa_{n,q}^{\vee}(\gamma)Sq^r=
\begin{cases}
0,& r\ \text{odd},\\[2mm]
\kappa_{n,q-r/2}^{\vee}(\gamma Sq^{r/2}),& r\ \text{even},
\end{cases}
\]
where the second term is interpreted as $0$ when $r/2>q$. Consequently,
\[
\kappa_{n,q}^{\vee}\big(P_{\mathbb A}\Gamma_n(q)\big)
\subseteq
P_{\mathbb A}\Gamma_n(d),
\]
and this restricted map is injective. Under the divided-power pairing, it is the dual of Kameko's induced homomorphism
\[
(\widetilde{Sq^0_*})_d:(QP_n)_d\longrightarrow (QP_n)_q.
\]
In particular, $(\widetilde{Sq^0_*})_d$ is surjective.
\end{mde}

\begin{proof}
It suffices to verify the displayed formula on a divided-power monomial
\[
\gamma=a_1^{(b_1)}\cdots a_n^{(b_n)}.
\]
Using the Cartan formula for the homological right action, we have
\[
\kappa_{n,q}^{\vee}(\gamma)Sq^r
=
\sum_{\ell_1+\cdots+\ell_n=r}
\left(\prod_{j=1}^{n}
\binom{2b_j+1-\ell_j}{\ell_j}\right)
 a_1^{(2b_1+1-\ell_1)}\cdots a_n^{(2b_n+1-\ell_n)},
\]
with the usual convention that terms with negative exponents are omitted. If $r$ is odd, then at least one $\ell_j$ is odd in every summand. For such an index $j$, the lower entry $\ell_j$ has last binary digit $1$, while the upper entry $2b_j+1-\ell_j$ is even. Lucas' theorem therefore gives
\[
\binom{2b_j+1-\ell_j}{\ell_j}=0\pmod 2,
\]
so every summand is zero.

Now suppose $r=2s$. Any summand with some odd $\ell_j$ vanishes by the same argument. Hence the only possible nonzero summands have $\ell_j=2m_j$ for all $j$, with $m_1+\cdots+m_n=s$. For such a summand, Lucas' theorem gives
\[
\binom{2b_j+1-2m_j}{2m_j}
=
\binom{b_j-m_j}{m_j}
\pmod 2.
\]
Substituting these identities into the Cartan formula yields exactly
\[
\kappa_{n,q}^{\vee}(\gamma)Sq^{2s}
=
\kappa_{n,q-s}^{\vee}(\gamma Sq^s).
\]
This proves the right-action formula.

If $\gamma\in P_{\mathbb A}\Gamma_n(q)$, then $\gamma Sq^s=0$ for every $s>0$. The formula above therefore implies that $\kappa_{n,q}^{\vee}(\gamma)$ is also annihilated by every positive Steenrod square, and hence is primitive. The map $\kappa_{n,q}^{\vee}$ is injective on the divided-power monomial basis, so its restriction to the primitive subspace is injective.

Finally, for a monomial $u^c=u_1^{c_1}\cdots u_n^{c_n}\in(P_n)_d$, one has
\[
\big\langle \kappa_{n,q}^{\vee}(a_1^{(b_1)}\cdots a_n^{(b_n)}),u^c\big\rangle=1
\]
if and only if $c_j=2b_j+1$ for every $j$. This is equivalent to
\[
\big\langle a_1^{(b_1)}\cdots a_n^{(b_n)}, Sq^0_*(u^c)\big\rangle=1.
\]
Thus $\kappa_{n,q}^{\vee}$ is the transpose of $Sq^0_*$. Since it sends primitive homology into primitive homology, $Sq^0_*$ sends hit elements to hit elements and induces the displayed homomorphism on cohit quotients. The induced dual map is injective, and finite-dimensional duality implies that $(\widetilde{Sq^0_*})_d$ is surjective.
\end{proof}

We can see that Proposition \ref{thm:dual-kameko-epimorphism} improves the algorithmic formula above to the unconditional identity
\[
\dim\Ker\big((\widetilde{Sq^0_*})_d\big)
=
\dim P_{\mathbb A}\Gamma_n(d)-
\dim P_{\mathbb A}\Gamma_n(q),
\qquad d=2q+n.
\]
The following corollary connects Kameko lifts with the Milnor-operation orthogonality tests introduced earlier in this section.

\begin{hqu}[Kameko lifts and Milnor-operation tests]\label{cor:kameko-milnor-tests}
For every $s\geq1$ and every $q\geq0$,
\[
\kappa_{n,q}^{\vee}(\Gamma_n(q))
\subseteq
\Ker\big(R(P_s^0):\Gamma_n(2q+n)\to \Gamma_n(2q+n-2^s+1)\big),
\]
where $R(P_s^0)$ denotes the homological right action of $P_s^0$. Consequently every iterated Kameko lift of a primitive homology class is primitive and is annihilated by the right action of every $P_s^0$.
\end{hqu}

\begin{proof}
Since $P_s^0$ is primitive, its right action on a divided-power monomial is the sum of the one-variable right actions. On one variable,
\[
 a_j^{(m)}P_s^0=
\begin{cases}
 a_j^{(m-(2^s-1))},& m-(2^s-1)\ \text{is odd and }m\geq 2^s,\\[1mm]
0,&\text{otherwise}.
\end{cases}
\]
Every exponent occurring in a monomial of $\kappa_{n,q}^{\vee}(\Gamma_n(q))$ is of the form $2b_j+1$, hence is odd. Since $2^s-1$ is also odd, the difference $(2b_j+1)-(2^s-1)$ is even. Thus no one-variable summand can occur, and the right action of $P_s^0$ is zero. The assertion for iterated Kameko lifts follows from Proposition~\ref{thm:dual-kameko-epimorphism} and from the one-step statement just proved.
\end{proof}

In the rank-five family of Section~\ref{s3}, if $z_{d_1=71}$ is a primitive divided-power representative dual to the invariant class $[\theta_{d_1=71}]$, then the iterated Kameko lift $(\kappa_{5}^{\vee})^{t-1}(z_{d_1})$ is dual to $[\varphi^{t-1}(\theta_{d_1})]$ and is annihilated by the homological right action of $P_{t+5}^{0}$. This gives a second, purely homological explanation for the orthogonality test already forced by the inclusion $H_{d_t}^{\perp}\subseteq\Ker(T_t^*)$.

Now, the dual primitive-kernel algorithm is exact: its output is certified by the displayed matrix equalities over $\mathbb Z/2$. It is also independent of admissible monomials, lower-weight reductions, and the choice of a cohomological admissible basis. Its computational cost is different rather than universally smaller. The matrix $\mathcal R_{n,d}$ has $\binom{n+d-1}{n-1}$ columns, so high degrees still require sparse bit-packed linear algebra. Nevertheless, it gives a genuinely separate verification route. In the rank-five family of this paper, it would recover
\[
\dim(QP_5)_{d_0=33}=\dim\Ker(\mathcal R_{5,d_0}),
\qquad
\dim(QP_5)_{d_1=71}=\dim\Ker(\mathcal R_{5,d_1}),
\]
and, using the dual Kameko map,
\[
\dim\Ker((\widetilde{Sq^0_*})_{d_1})
=
\dim\Ker(\mathcal R_{5,d_1})-
\dim\Ker(\mathcal R_{5,d_0}).
\]
The same procedure computes the one-dimensional $GL_5$-coinvariant line in primitive homology which is dual to the $GL_5$-invariant line in the cohit quotient. In this sense the algorithm supplies a second, homological route to the numerical invariants established in Section~\ref{s3}, while the Milnor-operation test above gives an additional low-cost structural constraint inside this dual computation.

\section{Appendix}\label{s4}

This appendix provides links to our algorithm implemented on the computer algebra system \texttt{OSCAR}, as well as detailed computational output data for $(QP_5)_{d_t}$ with $t = 0$ and $t=1$, together with their invariants. The additional Milnor-operation layer in Theorem~\ref{thm:large-milnor-layer} is independent of the enumeration files: it follows from the derivation formula for $P_{t+5}^{0}$ and the disjoint-support argument given in the proof. The invariant-space computations in the body of the paper are reproduced from the same data by forming the matrices of $\rho_i-1$ on the relevant admissible bases and row reducing them over $\mathbb Z/2$. 

$\bullet$ The \texttt{OSCAR} implementation of our algorithm for computing $(QP_5)_{d_0 = 33}$ and $(QP_5)_{d_1 = 71}$ together with their invariants is available at GitHub:\\ 
\url{https://github.com/phucdv2018/Code-OSCAR-33-and-71/releases/tag/v1.0.0} (archived at DOI: \url{https://doi.org/10.5281/zenodo.17719793}) and included as a supplementary file \texttt{Code\_OSCAR\_33\_and\_71}

\medskip

Note that in this \texttt{OSCAR} code, within the section
\begin{verbatim}

-------- RUN (example) --------
ultra_sparse_bitpacked_main(5, 33)
\end{verbatim}
to obtain results for degree $71,$ one simply replaces the pair $(5, 33)$ by $(5, 71).$

\medskip

$\bullet$ Computational data for $(QP_5)_{d_0=33}$ and its invariants are available at Zenodo:\\ 
\url{https://doi.org/10.5281/zenodo.17601723}.

\medskip

$\bullet$ Computational data for $(QP_5)_{d_1=71}$ and its invariants are available at Zenodo:\\ 
\url{https://doi.org/10.5281/zenodo.17601864}.




\medskip

\section*{Data Availability}
The data supporting the findings of this study are openly available at the GitHub repository and Zenodo archive referenced in Appendix~\ref{s4}.


\end{document}